\documentclass[a4paper]{article}

\usepackage{geometry}
\usepackage{amsmath}
\usepackage{amssymb}
\usepackage{amsthm}
\usepackage{bm}
\usepackage{calc}
\usepackage[british]{babel}
\usepackage[UKenglish]{isodate}
\usepackage{dsfont}
\usepackage{enumitem}
\usepackage[OT2,OT1]{fontenc}
\usepackage{ifthen}
\usepackage[cal=cm,scr=boondoxo]{mathalfa}
\usepackage{pifont}
\usepackage{stmaryrd}
\usepackage{textcomp}
\usepackage[capitalise]{cleveref}
	\crefformat{equation}{#2(#1)#3}
\usepackage{tikz}
	\usetikzlibrary{arrows}
	\usetikzlibrary{patterns}
\usepackage{tikz-cd}
\usepackage[textwidth=3cm,colorinlistoftodos]{todonotes}
\numberwithin{equation}{section}			% Enumeration of equations
\swapnumbers						% Enumeration of theorem-like environments
\newcommand\cyr{%					% cyrillic font
\renewcommand\rmdefault{wncyr}%
\renewcommand\sfdefault{wncyss}%
\renewcommand\encodingdefault{OT2}%
\normalfont
\selectfont}
\DeclareTextFontCommand{\textcyr}{\cyr}

\newcounter{Enum}					% Enumerated list
\newenvironment{Enumerate}{\begin{enumerate}[label={\rm({\roman*})}]}{\end{enumerate}}
\newcommand{\Enumref}[1]{{\setcounter{Enum}{#1}{\rm(\roman{Enum})}}}

\newcounter{Enumalph}					% Enumerated list with letters
\newenvironment{Enumeratealph}{\begin{enumerate}[label={\rm({\alph*})}]}{\end{enumerate}}

\newcommand{\descriptionlabelsave}{}			% Itemized list
\newenvironment{Itemize}{%
	\renewcommand{\descriptionlabelsave}{\descriptionlabel}\renewcommand{\descriptionlabel}{$\triangleright$}%
	\begin{description}[leftmargin=15pt,itemindent=-5.2pt]}{%
	\end{description}\renewcommand{\descriptionlabel}{\descriptionlabelsave}}

\newcounter{StepsCount}					% Enumerated steps (e.g. in proof)
\newenvironment{Steps}{%
	\begin{list}{\ding{\value{StepsCount}}}{\usecounter{StepsCount} \leftmargin=0pt \labelwidth=12pt \itemindent=\labelwidth%
	\itemsep=5pt\listparindent=\parindent} \setcounter{StepsCount}{191}}{\end{list}}

\newenvironment{Ilist}{%			% Itemised list without indentation (not-enumerated steps)
	\begin{list}{$\triangleright$}{\leftmargin=0pt \labelwidth=11pt \itemindent=\labelwidth%
	\itemsep=5pt\listparindent=\parindent}}{\end{list}}

\theoremstyle{plain}
	\newtheorem{lemma}{Lemma}[section]
	\newtheorem{proposition}[lemma]{Proposition}
	\newtheorem{theorem}[lemma]{Theorem}
	
\theoremstyle{definition}
	\newtheorem{definitioN}[lemma]{Definition}
	\newtheorem{assumptioN}[lemma]{Assumption}
\theoremstyle{remark}
	\newtheorem{remarK}[lemma]{Remark}
	\newtheorem{examplE}[lemma]{Example}
\renewcommand{\qedsymbol}{\raisebox{-2pt}{\large\ding{113}}}% Box mit Schatten
\newcommand{\defendsymbol}{$\lozenge$}
\newcommand{\qedsymbolsave}{\qedsymbol}
\newenvironment{definition}{\begin{definitioN}}{
	\renewcommand{\qedsymbolsave}{\qedsymbol}\renewcommand{\qedsymbol}{\defendsymbol}
	\popQED{\qed}\renewcommand{\qedsymbol}{\qedsymbolsave}\end{definitioN}}
\newenvironment{assumption}{\begin{assumptioN}}{
	\renewcommand{\qedsymbolsave}{\qedsymbol}\renewcommand{\qedsymbol}{\defendsymbol}
	\popQED{\qed}\renewcommand{\qedsymbol}{\qedsymbolsave}\end{assumptioN}}
\newenvironment{remark}{\begin{remarK}}{
	\renewcommand{\qedsymbolsave}{\qedsymbol}\renewcommand{\qedsymbol}{\defendsymbol}
	\popQED{\qed}\renewcommand{\qedsymbol}{\qedsymbolsave}\end{remarK}}
\newenvironment{example}{\begin{examplE}}{
	\renewcommand{\qedsymbolsave}{\qedsymbol}\renewcommand{\qedsymbol}{\defendsymbol}
	\popQED{\qed}\renewcommand{\qedsymbol}{\qedsymbolsave}\end{examplE}}

\crefname{assumptioN}{Assumption}{Assumptions}
\Crefname{assumptioN}{Assumption}{Assumptions}
\crefname{remarK}{Remark}{Remarks}
\Crefname{remarK}{Remark}{Remarks}
\crefname{examplE}{Example}{Examples}
\Crefname{examplE}{Example}{Examples}

\newcommand{\mc}[1]{{\mathcal{#1}}}			% --- abbreviation ---
\newcommand{\ms}[1]{{\mathscr{#1}}}			% --- abbreviation ---
\newcommand{\bb}[1]{{\mathbb{#1}}}			% --- abbreviation ---
\newcommand{\qu}{\overline}				% --- abbreviation ---
\newcommand{\mr}{\mathring}				% --- abbreviation ---

\DeclareMathOperator{\RE}{Re}				% real part
\renewcommand{\Re}{\RE}
\DeclareMathOperator{\IM}{Im}				% imaginary part
\renewcommand{\Im}{\IM}

\DeclareMathOperator{\ran}{ran}				% range
\DeclareMathOperator{\Int}{Int}				% interior
\newcommand{\Dummy}{\text{\textvisiblespace\kern1pt}}	% Platzhaltersymbol fuer Funktionsargumente
\newcommand{\BigO}{\mathrm{O}}				% big o

\newcommand{\DS}{\colon\mkern3mu}			% delimiter for set definition
\newcommand{\DP}{\kern2pt{\mathrel{\mathop:}\kern5pt}}	% delimiter for predicate formula
\newcommand{\DF}{\colon}				% delimiter for function domain/codomain
\newcommand{\DE}{\mathrel{\mathop:}=}			% defining equality
\newcommand{\ED}{=\mathrel{\mathop:}}			% defining equality
\newcommand{\D}{\mathrm{d}}				% rm-d for integration differential
\newcommand{\DD}{\mkern4mu\mathrm{d}}			% distance and rm-d for integration differential

\DeclareMathOperator{\Hol}{Hol}				% holomorphic functions
\DeclareMathOperator{\Mer}{Mer}				% meromorphic functions
\DeclareMathOperator{\sgn}{sgn}				% signum function
\DeclareMathOperator{\dist}{dist}			% distance between sets

\newcommand{\RC}{\mathsf{RC}}

\begin{document}

%\begin{flushleft}
	\noindent
	{\Large\textbf{Karamata's theorem for regularised Cauchy transforms}}
	\\[5mm]
	\textsc{
	\hspace*{1ex}
	Matthias Langer
	\,\ $\ast$\,\
	Harald Woracek
		\hspace*{-14pt}
		\renewcommand{\thefootnote}{\fnsymbol{footnote}}
		\setcounter{footnote}{2}
		\footnote{The second author was supported by the joint project I~4600 of the Austrian
			Science Fund (FWF) and the Russian foundation of basic research (RFBR).}
		\renewcommand{\thefootnote}{\arabic{footnote}}
		\setcounter{footnote}{0}
	}
	\\[6mm]
%\end{flushright}
	{\small
	\textbf{Abstract:}
	We prove Abelian and Tauberian theorems for regularised Cauchy transforms
	of positive Borel measures on the real line whose distribution functions
	grow at most polynomially at infinity.  In particular, we relate the asymptotics
	of the distribution functions to the asymptotics of the regularised Cauchy transform.
	}
\begin{flushleft}
	{\small
%	\\[3mm]
	}
\end{flushleft}

%%%%%%%%%%%%%%%%%%%%%%%%%%%%%%%%%%%%%%%%%%%%%%%%%%%%%%%%%%%%%%%%%%%%%%%%%%%%%%%%%%%%%%%%%
%\tableofcontents
%\setcounter{page}{1}
%
%%%%%%%%%%%%%%%%%%%%%%%%%%%%%%%%%%%%%%%%%%%%%%%%%%%%%%%%%%%%%%%%%%%%%%%%%%%%%%%%%%%%%%%%%

%---------
%   TEXTBODY
%---------

%
%
%
\section{Introduction}
\label{L110}

Let $\mu$ be a positive measure on $\bb R$ such that $\int_{\bb R}(1+|t|)^{-1}\DD\mu(t)<\infty$.
The Cauchy transform of $\mu$ is the function 
\begin{equation}\label{L4}
	C[\mu](z)\DE\int_{\bb R}\frac 1{t-z}\DD\mu(t)
\end{equation}
defined and analytic in the open upper half-plane $\bb C^+$. 
It plays an important role in many areas, such as spectral theory, 
moment problems, complex analysis and random matrix theory.
A prominent particular case occurs when $\mu$ is supported on $[0,\infty)$.  
Then we speak of the Stieltjes transform of $\mu$ and write 
\begin{equation}\label{L300}
	S[\mu](z)\DE\int_{[0,\infty)}\frac 1{t-z}\DD\mu(t).
\end{equation}
This function is defined and analytic in the slit plane $\bb C\setminus[0,\infty)$. 

A measure $\mu$ can be reconstructed from its Cauchy (or Stieltjes) transform by 
means of the Stieltjes inversion formula,
\begin{equation}
\label{L61}
	\mu\bigl((\alpha,\beta)\bigr)=
	\lim_{\delta\downarrow0}\lim_{\varepsilon\downarrow 0}
	\frac{1}{\pi}\int_{\alpha+\delta}^{\beta-\delta}\Im C[\mu](x+i\varepsilon)\DD x
\end{equation}
for $-\infty<\alpha<\beta<\infty$. 
This formula can be seen as relating the local behaviour of $\mu$ at a point or on a bounded interval 
in $\bb R$ with the local behaviour of $C[\mu]$ around this point or interval: 
in order to evaluate the right-hand side of \cref{L61} only the values of $C[\mu](z)$ for $z$ in some rectangle 
$(\alpha,\beta)\times(0,\varepsilon)\subseteq\bb C^+$ have to be known. 

One question that has attracted a lot of attention is the relation between the asymptotics of $\mu$ at $\infty$ and the
asymptotics of its transform at $\infty$.  For the case of the Stieltjes transform results were obtained already 
in the early 20$^{\textup{th}}$ century: G.~Valiron \cite{valiron:1913}, E.C.~Titchmarsh \cite{titchmarsh:1927}, 
and G.H.~Hardy and J.E.~Littlewood \cite{hardy.littlewood:1930} proved that, for each $\gamma\in(-1,0)$,
\begin{equation}
\label{L117}
	S[\mu](-x)\sim cx^\gamma\quad\Leftrightarrow\quad\mu([0,t))\sim c't^{\gamma+1},
\end{equation}
where $c,c'$ are related by a certain formula.  Here the symbol $\sim$ means that the quotient of the 
left-hand and right-hand sides tends to $1$, and is understood for $x,t\to +\infty$.  
The asymptotics of $S[\mu]$ along the ray $e^{i\pi}(0,\infty)$ could be
substituted by the asymptotics along any ray contained in the domain of analyticity $\bb C\setminus[0,\infty)$
(allowing the constant $c$ to depend on the angle of the ray), or even non-tangentially. 
This early result about the Stieltjes transform was generalised in several directions, 
and there is a vast literature on that topic.  As examples we mention \cite{karamata:1931a} where J.~Karamata 
generalised \cref{L117} to growth of regular variation instead of power asymptotics, 
or \cite{mcclure.wong:1978} where asymptotic expansions with infinitely many terms instead of a monomial 
on the right-hand side of \cref{L117} is considered. 
In the bilateral case, meaning measures that are not semi-bounded, asymptotics have to be taken along the positive imaginary axis, 
or a ray in $\bb C^+$ or non-tangentially.  In this case there is much less known.  One of the main difficulties is that 
contributions from the positive and negative half-axes can cancel each other.

In many applications, e.g.\ spectral theory of Sturm--Liouville and Schr\"odinger operators, measures are used that 
grow faster at infinity: instead of satisfying $\int_{\bb R}(1+|t|)^{-1}\DD\mu(t)<\infty$ they have only power bounded tails, meaning that
$\int_{\bb R}(1+t^2)^{-\kappa-1}\DD\mu(t)<\infty$ for some $\kappa\in\bb N_0$.  For such measures the Cauchy transform 
\cref{L4} has to be redefined by including appropriate regularising summands in the integrand.  The most common case is that 
$\mu$ is Poisson integrable, i.e.\ $\kappa=0$, and a commonly used regularisation in this case is 
\begin{equation}\label{L5}
	\widetilde C[\mu](z)\DE\int_{\bb R}\Big(\frac 1{t-z}-\frac t{1+t^2}\Big)\DD\mu(t),
	\qquad z\in\bb C^+.
\end{equation}
Some Abelian and Tauberian theorems dealing with polynomial asymptotics in the bilateral case are given in 
\cite{pleijel:1963a,pleijel:1963b}, a Tauberian theorem for the Cauchy transform \cref{L4} can be found in \cite{selander:1963}, 
and an Abelian theorem of somewhat different type (and formulated for integration on the unit circle instead of the real line) is 
\cite{sheremeta:1997}. A Tauberian theorem for the regularised Cauchy transform \cref{L5} is stated in 
\cite{bennewitz:1989}; however, the proof given contains a mistake. 
Fortunately the result itself turns out to be true; see \cref{L16} and the discussion preceding it. 

In the current paper we prove Abelian and Tauberian theorems for higher-order regularised Cauchy transforms and growth of
regular variation in Karamata's sense (see \cref{L115} for this notion). We relate the asymptotics of $\mu([0,t))$ and 
$\mu((-t,0))$ when $t\to+\infty$ to the asymptotics of the higher-order regularised Cauchy transform when $z\to+i\infty$
radially or non-tangentially. 
The main result of the paper is \Cref{L28} where we give a full characterisation (including explicit formulae for constants) in
the generic case. There are some boundary cases, namely when the index of regular variation
is an integer, where only one direction is possible:
either the Abelian direction where we deduce properties
of the regularised Cauchy transform from properties of the measure,
or the Tauberian direction, which is the other way round.
The phenomenon that more complicated behaviour occurs at integer powers was already observed in \cite{pleijel:1963a,pleijel:1963b}. 
We investigate these exceptional cases more closely in \cref{L40}.

For the proof of our results we follow common lines and consider imaginary and real parts of the integral
separately.  The imaginary part can be written as a Stieltjes transform, and thus inherits being well behaved; 
see \cref{L45}.  Contrasting this, the real part is the difference of two Stieltjes transforms, and this is the point 
where cancellation may happen. 

Let us give a brief overview of the contents of the paper.
In \cref{L111} we define higher-order regularised Cauchy transforms and study basic properties. 
In particular, we explore the relation with generalised Nevanlinna functions in the sense of M.G.~Krein and H.~Langer 
\cite{krein.langer:1977}, characterise the range of the transform, and prove an analogue of the classical Grommer--Hamburger 
theorem that relates convergence of a sequence of measures to convergence of their Cauchy transforms.
We use the latter theorem to prove a basic Tauberian theorem in \cref{L112}.
The proofs of Abelian theorems, which are contained in \cref{L113}, 
use different methods: the main ingredients are Karamata's theorems. 
In \cref{L114} we combine the results from \cref{L112,L113} to prove
our main theorems.  We also provide counterexamples for the boundary cases.
Finally, in \cref{L115} we recall and extend some results on
regularly varying functions and Stieltjes transforms of measures supported 
on $[0,\infty)$.

\subsection*{Notation}

Throughout the paper we use the following conventions and notations. 

\begin{Ilist}
\item 
	We set $\bb N\DE\{1,2,\ldots\}$, $\bb N_0\DE\bb N\cup\{0\}$, and let $\bb C$ be the field of complex numbers.
\item 
	We always use the branches of the logarithm and complex powers
	which are analytic on $\bb C\setminus(-\infty,0]$ and take the value $0$ or $1$,
	respectively, at the point $1$.
\item 
	We set $\bb C^+\DE\{z\in\bb C\DS \Im z>0\}$ and $\bb C^-\DE\{z\in\bb C\DS \Im z<0\}$.
\item 
	For a domain $\Omega$ denote by $\Hol(\Omega)$ and $\Mer(\Omega)$ the set
	of holomorphic and of meromorphic functions on $\Omega$ respectively.
\item
	We use the notation $f\sim g$ to express that $\frac fg\to 1$,
	and the notation $f\ll g$ if $\frac fg\to 0$.
	Further, we write $f\lesssim g$ if there exists a constant $c>0$ such that $f\le cg$,
	and we write $f\asymp g$ if $f\lesssim g$ and $g\lesssim f$.
	The domain of validity will be stated or will be clear from the context.
\item 
	When we speak of a ``measure'', we always mean positive Borel measure (unless explicitly specified differently). 
\item
	Throughout the rest of the paper we use the Stieltjes transform $\ms S[\mu]$ as defined in \eqref{L22},
	where we use a different sign convention from the one used in \eqref{L300}.
\end{Ilist}

\section{Regularised Cauchy integrals}
\label{L111}
\subsection{Definition of higher-order regularised Cauchy integrals}

To start with, let us recall the characterisations of the ranges of the transforms
$C$ and $\widetilde C$ introduced in \eqref{L4} and \eqref{L5}.
These are classical results going back to F.~Riesz, G.~Herglotz and R.~Nevanlinna;
for a comprehensive account see, e.g.\ \cite{kac.krein:1968a} or \cite{gesztesy.tsekanovskii:2000}.

\begin{proposition}\label{L6}
	Let $q\in\Hol(\bb C^+)$.
	\begin{Enumerate}
	\item 
		The function $q$ can be represented in the form
		\begin{equation}\label{L20}
			q(z)=a+C[\mu](z)=a+\int_{\bb R}\frac 1{t-z}\DD\mu(t), \qquad z\in\bb C^+,
		\end{equation}
		with some $a\in\bb R$ and a positive measure $\mu$ on $\bb R$ with
		$\int_{\bb R}(1+|t|)^{-1}\D\mu(t)<\infty$ if and only if
		\[
			\forall z\in\bb C^+\DP \Im q(z)\ge 0
			\quad\text{and}\quad
			\int_1^\infty\frac{\Im q(iy)}y\DD y < \infty.
		\]
	\item 
		The function $q$ can be represented in the form
		\begin{equation}\label{L7}
			q(z)=a+bz+\widetilde C[\mu](z)=a+bz+\int_{\bb R}\Big(\frac 1{t-z}-\frac t{1+t^2}\Big)\DD\mu(t), 
			\qquad z\in\bb C^+,
		\end{equation}
		with some $a\in\bb R$, $b\ge 0$ and a positive measure $\mu$ on $\bb R$ with
		$\int_{\bb R}(1+t^2)^{-1}\D\mu(t)<\infty$ if and only if
		\[
			\forall z\in\bb C^+\DP \Im q(z) \ge 0.
		\]
	\end{Enumerate}
\end{proposition}

\begin{remark}\label{L18}
\hfill
\begin{Enumerate}
\item
	Assume that $q$ is represented in the form \eqref{L7}.
	Then the constants $a,b$ are given by
	\begin{equation}\label{L185}
		a=\Re q(i), \qquad b=\lim_{y\to+\infty}\frac 1{iy}q(iy),
	\end{equation}
	and the measure $\mu$ is given by the Stieltjes inversion formula
	\[
		\forall \alpha,\beta\in\bb R,\,\alpha<\beta\DP
		\lim_{\delta\downarrow0}\lim_{\varepsilon\downarrow 0}
		\frac{1}{\pi}\int_{\alpha+\delta}^{\beta-\delta}\Im q(x+i\varepsilon)\DD x
		= \mu\bigl((\alpha,\beta)\bigr);
	\]
	see, e.g.\ \cite{kac.krein:1968a}.
	In particular, the map $(a,b,\mu)\mapsto q$ with $q$ satisfying \eqref{L7}
	is injective.
\item
	Assume that $a=b=0$ in \eqref{L7}, i.e.\ $q(z)=\widetilde C[\mu](z)$.
	Then
	\begin{equation}\label{L183}
		|q(iy)| \ll y, \qquad y\to+\infty,
	\end{equation}
	by \eqref{L185}, and
	\begin{equation}\label{L184}
		\lim_{y\to+\infty}y\Im q(iy) = \sup_{y>0}y\Im(iy) = \mu(\bb R);
	\end{equation}
	see again \cite{kac.krein:1968a}.
\end{Enumerate}
\end{remark}

\medskip\noindent
It is apparent that moving from function \eqref{L20} to function \eqref{L7}
is only the first step on a ladder:
instead of Poisson-integrable measures one may use measures whose tails have
at most power growth, and instead of the term $a+bz$ one may use any polynomial
with real coefficients. In the integral higher-order regularisation will become necessary.

We work with a scale of higher-order regularised Cauchy transforms which is
commonly used in the framework of indefinite inner product spaces;
see, e.g.\ \cite{krein.langer:1977}.

\begin{definition}\label{L8}
	Let $\kappa\in\bb N_0$.
	\begin{Enumerate}
	\item
		We denote by $\bb E_{\leq\kappa}$ the set of all pairs $(\mu,p)$ where
		\begin{Itemize}
		\item
			$\mu$ is a measure on $\bb R$ that satisfies
			\begin{equation}\label{L49}
				\int_{\bb R}\frac{\D\mu(t)}{(1+t^2)^{\kappa+1}}<\infty;
			\end{equation}
		\item
			$p$ is a polynomial with real coefficients whose degree does not exceed $2\kappa+1$;
		\item
			the coefficient of $z^{2\kappa+1}$ in $p$ satisfies
			\begin{equation}\label{L77}
				\frac{1}{(2\kappa+1)!}p^{(2\kappa+1)}(0)
				\ge \int_{\bb R}\frac{\D\mu(t)}{(1+t^2)^{\kappa+1}}.
			\end{equation}
		\end{Itemize}
	\item
		The $\kappa$-regularised Cauchy transform is the
		map $\ms C_\kappa\DF\bb E_{\le\kappa}\to\Hol(\bb C^+)$ defined by
		\begin{equation}\label{L14}
			\ms C_\kappa[\mu,p](z)
			\DE p(z)+(1+z^2)^{\kappa+1}\int_{\bb R}\frac 1{t-z}\cdot
			\frac{\D\mu(t)}{(1+t^2)^{\kappa+1}},
			\qquad z\in\bb C^+.
		\end{equation}
	\end{Enumerate}
\end{definition}

\medskip\noindent
In order to represent polynomials with real coefficients as regularised Cauchy transforms, 
we have to include pairs such as $(0,z^{2\kappa})$ in $\bb E_{\leq\kappa}$. 
For this reason we cannot speak of ``the leading coefficient of $p$'' in \cref{L77}.

These maps can indeed be seen as higher-order regularised Cauchy integrals:
for $k\in\bb N_0$ we have
\begin{align}
	& \frac{1}{t-z}-(t+z)\sum_{j=0}^k\frac{(1+z^2)^j}{(1+t^2)^{j+1}}
	= \frac{1}{t-z}-\frac{t+z}{1+t^2}
	\cdot\frac{1-\bigl(\frac{1+z^2}{1+t^2}\bigr)^{k+1}}{1-\frac{1+z^2}{1+t^2}}
	\nonumber\\[1ex]
	&= \frac{1}{t-z}-(t+z)\cdot\frac{1-\bigl(\frac{1+z^2}{1+t^2}\bigr)^{k+1}}{t^2-z^2}
	= \frac{1}{t-z}\cdot\frac{(1+z^2)^{k+1}}{(1+t^2)^{k+1}};
	\label{L31}
\end{align}
hence we obtain, with $k=\kappa$,
\begin{equation}\label{L38}
	\ms C_\kappa[\mu,p](z) = p(z)
	+ \int_{\bb R}
	\biggl[\frac{1}{t-z}-(t+z)\sum_{j=0}^\kappa\frac{(1+z^2)^j}{(1+t^2)^{j+1}}\biggr]
	\DD\mu(t).
\end{equation}
Note that the regularising terms in the integral on the right-hand side of \eqref{L38}
are the first $\kappa+1$ terms of an expansion of $\frac{1}{t-z}=(t+z)\cdot\frac{1}{t^2-z^2}$
in terms of powers of $\frac{1}{1+t^2}$.
%since the integrand in \eqref{L14} can be written as
%\begin{equation*}%\label{L31}
%	(1+z^2)^{\kappa+1}\cdot\frac 1{t-z}\cdot\frac{1}{(1+t^2)^{\kappa+1}}
%	= \frac{1}{t-z}-(t+z)\sum_{j=0}^\kappa\frac{(1+z^2)^j}{(1+t^2)^{j+1}}
%\end{equation*}
%\todo[color=lime]{include calculation?}
%as a simple calculation shows.
%From this we obtain, for $\kappa=0$, that
For $\kappa=0$ relation \eqref{L38} reads as
\begin{equation}\label{L86}
	\ms C_0[\mu,p](z) = p(z) + \int_{\bb R}\Bigl(\frac{1}{t-z}-\frac{t+z}{1+t^2}\Bigr)\DD\mu(t),
\end{equation}
which yields the following connection with the previously discussed regularised
Cauchy-type integral \eqref{L5}:
\begin{Itemize}
\item
	$a+bz+\widetilde C[\mu](z)=\ms C_0[\mu,p](z)$
	\; with \;
	$\displaystyle p(z)\DE a+\biggl(b+\int_{\bb R}\frac{\D\mu(t)}{1+t^2}\biggr)z$;
\item
	$\displaystyle
	(a,b,\mu)\in
	\bb R\times[0,\infty)\times\biggl\{\mu\DS\text{positive measure with }\int_{\bb R}\frac{\D\mu(t)}{1+t^2}<\infty\biggr\}
	$
	\\[1ex]
	$\displaystyle
	\Longleftrightarrow\quad (\mu,p)\in\bb E_{\le 0}$ \;
	with $p$ related to $a$, $b$ and $\mu$ as above.
\end{Itemize}

\begin{remark}
\label{L39}
	Using again \eqref{L31} with $\kappa=0$ we obtain the following representation
	for $\ms C_\kappa$ for arbitrary $\kappa\in\bb N_0$:
	\begin{align}
		\ms C_\kappa[\mu,p](z)
		&= p(z) + (1+z^2)^\kappa\int_{\bb R}\Bigl(\frac{1}{t-z}-\frac{t+z}{1+t^2}\Bigr)
		\frac{\D\mu(t)}{(1+t^2)^\kappa}
		\nonumber
		\\[1ex]
		&= \bigg(p(z)-z(1+z^2)^\kappa\int_{\bb R}\frac{\D\mu(t)}{(1+t^2)^{\kappa+1}}\bigg)
		\nonumber
		\\[1ex]
		&\quad +(1+z^2)^\kappa\int_{\bb R}\Bigl(\frac{1}{t-z}-\frac{t}{1+t^2}\Bigr)
		\frac{\D\mu(t)}{(1+t^2)^\kappa}.
		\label{L30}
	\end{align}
	If the stronger integrability condition
	$\int_{\bb R}(1+|t|)^{-(2\kappa+1)}\D\mu(t)<\infty$ is satisfied,
	then we can split the second integral on the right-hand side of \eqref{L30}
	and rewrite it as
	\begin{align}
		\ms C_\kappa[\mu,p](z) &= \biggl(p(z)-(1+z^2)^\kappa
		\biggl[z\int_{\bb R}\frac{\D\mu(t)}{(1+t^2)^{\kappa+1}}+
		\int_{\bb R}\frac{t}{1+t^2}\cdot\frac{\D\mu(t)}{(1+t^2)^\kappa}\biggr]\biggr)
		\nonumber
		\\
		&\quad +(1+z^2)^\kappa\int_{\bb R}\frac{1}{t-z}\cdot\frac{\D\mu(t)}{(1+t^2)^\kappa}.
		\label{L42}
	\end{align}
\end{remark}

\medskip\noindent
Before we collect some properties of $\bb E_{\le\kappa}$ and $\ms C_\kappa$,
we recall the Stieltjes--Liv\v{s}ic inversion formula; see, e.g.\
\cite[Corollary~II.1.2]{langer:1982} or \cite[Theorem~1.2.4]{gorbachuk.gorbachuk:1997}).
Let $\sigma$ be a finite measure on $\bb R$, let $\alpha,\beta\in\bb R$ with $\alpha<\beta$,
and let $f$ be an analytic function on a neighbourhood of $[\alpha,\beta]$.
For $\delta,\varepsilon>0$ let $\Gamma_\varepsilon^\delta$ be the path consisting
of the two directed line segments
\begin{equation}\label{L75}
	\alpha+\delta-i\varepsilon\rightsquigarrow \beta-\delta-i\varepsilon
	\quad\text{and}\quad
	\beta-\delta+i\varepsilon\rightsquigarrow \alpha+\delta+i\varepsilon.
\end{equation}
Then
\begin{equation}\label{L74}
	\lim_{\delta\downarrow 0}\lim_{\varepsilon\downarrow 0}
	\frac{-1}{2\pi i}\int_{\Gamma_\varepsilon^\delta}f(z)\int_{\bb R}\frac{1}{t-z}\DD\sigma(t)\DD z
	= \int_{(\alpha,\beta)}f(t)\DD\sigma(t).
\end{equation}

\begin{lemma}\label{L21}
	Let $\kappa\in\bb N_0$.
	\begin{Enumerate}
	\item
		The set $\bb E_{\le\kappa}$ is a positive cone and $\ms C_\kappa$ is a cone map,
		i.e.\ compatible with finite sums and non-negative scalar multiples.
	\item
		The map $\ms C_\kappa$ is injective.  For $q\in\ran\ms C_\kappa$
		the element $(\mu,p)=\ms C_\kappa^{-1}q$ is obtained as follows:
		the polynomial $p$ can be recovered from solving the $2\kappa+2$ equations
		obtained by splitting real and imaginary parts of
		\begin{equation}\label{L15}
			q^{(j)}(i) = p^{(j)}(i), \qquad j\in\{0,\ldots,\kappa\};
		\end{equation}
		the measure $\mu$ can be obtained via the Stieltjes inversion formula:
		for $\alpha,\beta\in\bb R$ with $\alpha<\beta$ we have
		\begin{equation}\label{L73}
			\mu\bigl((\alpha,\beta)\bigr)
			= \lim_{\delta\downarrow 0}\lim_{\varepsilon\downarrow 0}
			\frac{1}{\pi}\int_{\alpha+\delta}^{\beta-\delta} \Im q(t+i\varepsilon)\DD\mu(t).
		\end{equation}
	\item
		Let $\kappa'>\kappa$.
		Then the inclusion $\ran\ms C_\kappa\subseteq\ran\ms C_{\kappa'}$ holds,
		and, for $(\mu,p)\in\bb E_{\le\kappa}$,
		we have $\ms C_\kappa[\mu,p]=\ms C_{\kappa'}[\mu,\widetilde p]$ with
		\[
			\widetilde p(z) = p(z) + \sum_{j=\kappa+1}^{\kappa'}(1+z^2)^j
			\biggl[z\int_{\bb R}\frac{\D\mu(t)}{(1+t^2)^{j+1}}
			+ \int_{\bb R}\frac{t\DD\mu(t)}{(1+t^2)^{j+1}}\biggl].
		\]
	\end{Enumerate}
\end{lemma}

\begin{proof}
	\phantom{}

(i)
	The statements are clear from the definitions of $\bb E_{\le\kappa}$ and $\ms C_\kappa$.

(ii)
	Let $(\mu,p)\in\bb E_{\le\kappa}$ and set $q=\ms C_\kappa[\mu,p]$.
	It follows from the definition of $\ms C_\kappa$ that \eqref{L15} holds,
	which implies that $p$ is uniquely determined by $q$.
	To show \eqref{L73}, let us first extend $q$ to $\bb C\setminus\bb R$ by symmetry:
	$q(z)\DE\qu{q(\qu z)}$ for $z\in\bb C^-$.
	Moreover, let $\alpha,\beta\in\bb R$ with $\alpha<\beta$ and let
	$\Gamma_\varepsilon^\delta$ be the path in \eqref{L75}.
	Then \eqref{L74} implies that
	\begin{align*}
		& \lim_{\delta\downarrow 0}\lim_{\varepsilon\downarrow 0}
		\frac{1}{\pi}\int_{\alpha+\delta}^{\beta-\delta} \Im q(t+i\varepsilon)\DD\mu(t)
		= \lim_{\delta\downarrow 0}\lim_{\varepsilon\downarrow 0}
		\frac{-1}{2\pi i}\int_{\Gamma_\varepsilon^\delta}q(z)\DD z
		\\
		&=
		\lim_{\delta\downarrow 0}\lim_{\varepsilon\downarrow 0}\frac{-1}{2\pi i}
		\int_{\Gamma_\varepsilon^\delta}(1+z^2)^{\kappa+1}
		\bigg[\int_{\bb R}\frac{1}{t-z}\cdot\frac{\D\mu(t)}{(1+t^2)^{\kappa+1}}\bigg]
		\DD z
		\\
		&=
		\int_{(\alpha,\beta)}(1+t^2)^{\kappa+1}\cdot\frac{\D\mu(t)}{(1+t^2)^{\kappa+1}}
		= \mu\big((\alpha,\beta)\big).
	\end{align*}
	The unique determination of $\mu$ and $p$ shows that $\ms C_\kappa$ is injective.

(iii)
	Let $(\mu,p)\in\bb E_{\le\kappa}$ and let $\kappa'>\kappa$.
	It follows from \eqref{L31} that
	\begin{equation}\label{L35}
		\frac{1}{t-z}\cdot\frac{(1+z^2)^{\kappa+1}}{(1+t^2)^{\kappa+1}}
		= \frac{1}{t-z}\cdot\frac{(1+z^2)^{\kappa'+1}}{(1+t^2)^{\kappa'+1}}
		+ (t+z)\sum_{j=\kappa+1}^{\kappa'}\frac{(1+z^2)^j}{(1+t^2)^{j+1}},
	\end{equation}
	which yields the statement in (iii).
\end{proof}

\medskip\noindent
Note that the Stieltjes inversion formula \eqref{L73}
for the recovery of $\mu$ is independent of $\kappa$.

\subsection[Determining the range of $\ms C_\kappa$]{Determining the range of \boldmath{$\ms C_\kappa$}}

An intrinsic characterisation of the range of $\ms C_\kappa$ along the
lines of \cref{L6}\,\Enumref{2} can be given.
This is based on \cite{krein.langer:1977,krein.langer:1981,langer:1986}
and related to \cite[Theorem~3.9]{langer.woracek:ninfrep}
(a predecessor of the latter is \cite[Lemma~3.6]{kaltenbaeck.woracek:polya}).

Let us first recall the definition of generalised Nevanlinna functions
in the sense of \cite{krein.langer:1977}.
We need, in particular, functions from the subclasses $\mc N_\kappa^{(\infty)}$,
which are characterised by a special behaviour at infinity and which were studied in, e.g.\
\cite{dijksma.langer.shondin:2004,dijksma.luger.shondin:2009,dijksma.luger.shondin:2010,hassi.luger:2006,
langer.woracek:gpinf,langer.woracek:ninfrep,langer.woracek:sinham}.

\begin{definition}\label{L19}
	For $q\in\Mer(\bb C^+)$ we denote by $\Omega_q$ its domain of analyticity;
	for the constant $q\equiv\infty$ we set $\Omega_q=\emptyset$.
	\begin{Enumerate}
	\item
		For $q\in\Mer(\bb C^+)\cup\{\infty\}$ we denote by $\kappa_q\in\bb N_0\cup\{\infty\}$
		the number of negative squares of the Hermitian kernel
		\begin{equation}\label{L23}
			K_q(w,z) \DE \frac{q(z)-\qu{q(w)}\,}{z-\qu w}, \qquad z,w\in\Omega_q,
		\end{equation}
		i.e.\ the supremum of all numbers of negative squares of the quadratic forms
		\[
			\sum_{i,j=1}^m K_q(w_j,w_i)\xi_i\qu{\xi_j}
		\]
		with $m\in\bb N$ and $w_1,\ldots,w_m\in\Omega_q$.
	\item
		Let $\kappa\in\bb N_0$.  We denote by $\mc N_\kappa$ the set of all
		functions $q\in\Mer(\bb C^+)$ with $\kappa_q=\kappa$.
		Moreover, we set $\mc N_{\le\kappa}\DE\bigcup\limits_{\kappa'=0}^\kappa \mc N_{\kappa'}$
		and $\mc N_{<\infty}\DE\bigcup\limits_{\kappa'=0}^\infty \mc N_{\kappa'}$.
	\item
		Let $\kappa\in\bb N_0$.  We denote by $\mc N_{\kappa}^{(\infty)}$
		the set of all functions $q\in\mc N_\kappa$ for which
		\[
			\lim_{y\to+\infty}\bigg|\frac{q(iy)}{y^{2\kappa-1}}\bigg| = \infty
			\qquad\text{or}\qquad
			\lim_{y\to+\infty}\frac{q(iy)}{(iy)^{2\kappa-1}}\in(-\infty,0).
		\]
		Moreover, we set $\mc N_{\le\kappa}^{(\infty)}
		\DE\bigcup\limits_{\kappa'=0}^\kappa \mc N_{\kappa'}^{(\infty)}$
		and $\mc N_{<\infty}^{(\infty)}\DE\bigcup\limits_{\kappa'=0}^\infty \mc N_{\kappa'}^{(\infty)}$.
	\end{Enumerate}
\end{definition}

\begin{remark}\label{L68}
	\phantom{}
	\begin{Enumerate}
	\item
		Note that the classes $\mc N_0$ and $\mc N_0^{(\infty)}$ coincide with the set of 
		all Nevanlinna functions,
		i.e.\ those functions $q$ that are analytic on $\bb C^+$ and satisfy $\Im q(z)\ge0$
		for $z\in\bb C^+$.
		Further, functions in $\mc N_\kappa$ have at most $\kappa$ poles and
		at most $\kappa$ zeros in $\bb C^+$.
	\item
		The classes $\mc N_\kappa^{(\infty)}$,
		which have also been denoted by $\mc N_\kappa^\infty$ in the literature,
		can also be characterised differently, namely,
		for $q\in\Mer(\bb C^+)$ the following conditions are equivalent
		(see \cite{hassi.luger:2006}, or also \cite{dijksma.langer.shondin:2004,kaltenbaeck.woracek:polya}):
		\begin{Enumeratealph}
		\item
			$q\in\mc N_{<\infty}^{(\infty)}$;
		\item
			$\infty$ is the only (generalised) pole not of positive type
			(in the sense of \cite[\S3]{krein.langer:1977}),
			i.e.\ $\infty$ is the only (generalised) eigenvalue with a non-positive
			eigenvector of a representing relation in a Pontryagin space
			(see also \cite{langer:1986} for an analytic characterisation of generalised
			poles not of positive type);
		\item
			there exist $m\in\bb N_0$, a real polynomial $p$ and a measure $\sigma$
			on $\bb R$ such that $\int_{\bb R}(1+t^2)^{-1}\D\sigma(t)<\infty$ and
			\begin{equation}\label{L83}
				q(z) = (1+z^2)^m\int_{\bb R}\Bigl(\frac{1}{t-z}-\frac{t}{1+t^2}\Bigr)\DD\sigma(t)
				+ p(z);
			\end{equation}
		\item
			there exist $n\in\bb N_0$, $\beta_1,\ldots,\beta_n\in\bb C^+\cup\bb R$,
			$\rho_j\in\bb N$ and $q_0\in\mc N_0$ such that
			\[
				q(z) = q_0(z)\prod_{j=1}^n(z-\beta_j)^{\rho_j}(z-\qu{\beta_j})^{\rho_j}.
			\]
		\end{Enumeratealph}
		If $q$ is as in \eqref{L83} and $\deg p=l$ with leading coefficient $c_l$,
		then $q\in\mc N_\kappa^{(\infty)}$ with
		\begin{equation}\label{L84}
			\kappa \le \max\{m,\kappa_p\} \qquad\text{where}\quad
			\kappa_p =
			\begin{cases}
				\frac{l}{2}, & l \ \text{even},
				\\[1ex]
				\frac{l-\sgn(c_l)}{2}, & l \ \text{odd};
			\end{cases}
		\end{equation}
		equality holds in \eqref{L84} if $m=0$ or $\sigma$ is an infinite measure;
		see, e.g.\ \cite[(1.16)]{dijksma.langer.shondin:2004}.
		Note that for a real polynomial $p$ one has $\kappa_p\le\kappa'$
		if and only if $\deg p\le2\kappa'+1$ and $p^{(2\kappa'+1)}(0)\ge0$.
	\item
		The representation \eqref{L38} is a special case of the integral representation
		of $\mc N_{<\infty}$-functions given in \cite[Satz~3.1]{krein.langer:1977}.
	\item
		In \cite{langer.woracek:ninfrep} representations of functions
		in $\mc N_{<\infty}^{(\infty)}$ were constructed with distributions
		(more precisely, distributional densities) on the one-point compactification
		$\bb R\cup\{\infty\}$ of $\bb R$ which act like measures on $\bb R$.
	\item
		Functions in $\mc N_{<\infty}^{(\infty)}$
		are analytic in $\bb C^+$.
	\end{Enumerate}
\end{remark}

\begin{theorem}\label{L9}
	For every $\kappa\in\bb N_0$ the equality $\ran\ms C_\kappa=\mc N_{\le\kappa}^{(\infty)}$
	holds.
\end{theorem}

\begin{proof}
	Let $\kappa\in\bb N_0$.
	It follows from \eqref{L30} that $q\in\ran\ms C_\kappa$ if and only if
	it can be written as
	\begin{equation}\label{L85}
		q(z) = \widetilde p(z)
		+ (1+z^2)^\kappa\int_{\bb R}\Bigl(\frac{1}{t-z}-\frac{t}{1+t^2}\Bigr)\DD\sigma(t)
	\end{equation}
	with a real polynomial $\widetilde p$ of degree at most $2\kappa+1$
	with $\widetilde p^{(2\kappa+1)}(0)\ge0$
	and a measure $\sigma$ such that $\int_{\bb R}(1+t^2)^{-1}\D\sigma(t)<\infty$.

	First assume that $q\in\ran\ms C_\kappa$.  Then \eqref{L85} holds
	with $\sigma$ and $\widetilde p$ as above.
	By \cref{L68}\,(ii) we obtain that $q\in\mc N_{\le\kappa}^{(\infty)}$.

	Conversely, assume that $q\in\mc N_{\le\kappa}^{(\infty)}$,
	say $q\in\mc N_{\kappa'}^{(\infty)}$ with $\kappa'\in\{0,\ldots,\kappa\}$.
	Then there exists a representation of $q$ as in \eqref{L83} such that
	$\kappa'=\max\{m,\kappa_p\}$, where $\kappa_p$ is as in \eqref{L84};
	in particular $m\le\kappa$ and $\kappa_p\le\kappa$.
	The function
	\[
		\widehat q(z) \DE (1+z^2)^m
		\int_{\bb R}\Bigl(\frac{1}{t-z}-\frac{t}{1+t^2}\Bigr)\DD\sigma(t)
	\]
	belongs to $\ran\ms C_m$ by the first paragraph of this proof.
	Since $m\le\kappa$, we obtain from \cref{L21}\,(iii) that $\widehat q\in\ran\ms C_\kappa$.
	The relations $\deg p\le2\kappa+1$ and $p^{(2\kappa+1)}(0)\ge0$ show
	that $(0,p)\in\bb E_{\le\kappa}$, and hence $p=\ms C_\kappa[0,p]\in\ran\ms C_\kappa$.
	Now \cref{L21}\,(i) implies that $q=\widehat q+p\in\ran\ms C_\kappa$.
\end{proof}

\subsection[$\ms C_\kappa$ as a homeomorphism: the Grommer--Hamburger theorem]{\boldmath{$\ms C_\kappa$} as a homeomorphism: the Grommer--Hamburger theorem}

Next we discuss a continuity property of $\ms C_\kappa$; see \cref{L71} below.
This result is a variant of a classical theorem of J.~Grommer and H.~Hamburger;
see the discussion in \Cref{L52}.

Before being able to formulate a result, we have to make clear which topologies we use.
On the set $\Hol(\bb C^+)$ we always use the topology of locally uniform convergence.
Topologising $\bb E_{\le\kappa}$ is slightly more subtle.  We proceed as follows.
Fix $\kappa\in\bb N_0$.
The set of all positive measures $\mu$ that satisfies \eqref{L49} is
a subset of the dual space of the weighted $C_0$-space
\[
	C_0(\bb R,\omega_\kappa)
	\DE \Bigl\{f\in C(\bb R)\DS \lim_{|x|\to\infty}|f(x)|\omega_\kappa(x)=0\Bigr\},
	\qquad
	\|f\|\DE\sup_{x\in\bb R}|f(x)|\omega_\kappa(x),
\]
where $\omega_\kappa$ is the weight function $\omega_\kappa(x)\DE(1+x^2)^{\kappa+1}$;
note that
\begin{equation}\label{L109}
	\|\mu\|_{C_0(\bb R,\omega_k)'} = \int_{\bb R}\frac{\D\mu(t)}{(1+t^2)^{\kappa+1}}
\end{equation}
for a positive measure $\mu$ that satisfies \eqref{L49}.
We endow the set of all positive measures $\mu$ that satisfy \eqref{L49}
with the subspace topology of the $w^*$-topology in $C_0(\bb R,\omega_\kappa)'$.
The set of all polynomials of degree at most $2\kappa+1$ is isomorphic to
$\bb R^{2\kappa+2}$ mapping a polynomial to its coefficients,
and we endow polynomials with the Euclidean norm transported via this isomorphism.
The set $\bb E_{\le\kappa}$ is now topologised as a subspace of the product.

This topology has some very nice properties, which are summarised in the following lemma.

\begin{lemma}\label{L11}
	Let $\kappa\in\bb N_0$.
	\begin{Enumerate}
	\item
		$\bb E_{\le\kappa}$ is a closed subset of
		$C_0(\bb R,\omega_\kappa)'\times\bb R^{2\kappa+2}$.
	\item
		A subset $\mc E\subseteq\bb E_{\leq\kappa}$ is relatively compact if and only if
		\begin{equation}\label{L12}
			\sup\big\{\|p\|\DS(\mu,p)\in\mc E\big\} < \infty.
		\end{equation}
	\item
		The sets
		\[
			\mc E_N\DE\big\{(\mu,p)\in\bb E_{\le\kappa}\DS\|p\|\leq N\big\}
		\]
		are compact.  We have $\mc E_N\subseteq\Int\mc E_{N+1}$,
		and $\bigcup_{N\in\bb N}\mc E_N=\bb E_{\leq\kappa}$,
		where the interior $\Int\mc E_{N+1}$ of $\mc E_{N+1}$ is understood 
		within $\bb E_{\leq\kappa}$.
	\end{Enumerate}
\end{lemma}

\begin{proof}
	\phantom{}

	\Enumref{1}
	Assume that $((\mu_i,p_i))_{i\in I}$ is a net in $\bb E_{\leq\kappa}$ that converges
	to some element $(\mu,p)\in C_0(\bb R,\omega_\kappa)'\times\bb R^{2\kappa+2}$.
	Then $\mu$ is again a positive measure, and, by \eqref{L109},
	\begin{align*}
		\frac 1{(2\kappa+1)!}p^{(2\kappa+1)}(0)
		&= \lim_{i\in I}\frac 1{(2\kappa+1)!}p_i^{(2\kappa+1)}(0)
		\\
		&\ge \limsup_{i\in I}\int_{\bb R}\frac{\D\mu_i(t)}{(1+t^2)^{\kappa+1}}
		\ge \int_{\bb R}\frac{\D\mu(t)}{(1+t^2)^{\kappa+1}}.
	\end{align*}
	Thus $(\mu,p)\in\bb E_{\leq\kappa}$, and we see that $\bb E_{\leq\kappa}$ is indeed closed.

	\Enumref{2}
	Let $\pi_2$ be the projection onto the second component
	of $C_0(\bb R,\omega_\kappa)'\times\bb R^{2\kappa+2}$. Then $\pi_2$ is continuous.
	This already shows the implication ``$\Rightarrow$'' in \Enumref{2}.
	Conversely, assume that \eqref{L12} holds.
	It follows from \eqref{L109} that
	\[
		\sup\Bigl\{\|\mu\|_{C_0(\bb R,\omega_\kappa)'}\DS (\mu,p)\in\mc E\Bigr\}
		\le \sup\bigl\{\|p\|\DS(\mu,p)\in\mc E\bigr\}
		\ED c < \infty.
	\]
	We see that
	\[
		\mc E
		\subseteq \{\mu\in C_0(\bb R,\omega_\kappa)'\DS\|\mu\|_{C_0(\bb R,\omega_\kappa)'}\le c\}
		\times \{p\in\bb R^{2\kappa+2}\DS\|p\|\le c\},
	\]
	and hence $\mc E$ is relatively compact in $C_0(\bb R,\omega_\kappa)'\times\bb R^{2\kappa+2}$
	by the Banach--Alaoglu Theorem.
	Since $\bb E_{\le\kappa}$ is closed in this product space, $\mc E$ is
	also relatively compact in $\bb E_{\leq\kappa}$.

	\Enumref{3} The statement follows from what we have shown so far and from the fact
	that the continuity of $\pi_2$ implies that $\mc E_N$ is closed and that
	\[
		\Int\mc E_{N+1}=\big\{(\mu,p)\in\bb E_{\leq\kappa}\DS\|p\|<N+1\big\}.
		\qedhere
	\]
\end{proof}

\begin{proposition}\label{L13}
	The range of $\ms C_\kappa$ is closed in $\Hol(\bb C^+)$, and $\ms C_\kappa$
	is a homeomorphism onto its range.
\end{proposition}

\begin{proof}
	Continuity of $\ms C_\kappa$ is clear from our choice of topology.
	Let $(q_i)_{i\in I}$ be a net in $\ran\ms C_\kappa$, and assume
	that $\lim_{i\in I}q_i=q$ in $\Hol(\bb C^+)$.
	Remembering \eqref{L15} we find $i_0\in I$ and $N\in\bb N$
	such that $\ms C_\kappa^{-1}(q_i)\in\mc E_N$ for all $i\ge i_0$.
	Since $\mc E_N$ is compact (and $\ms C_\kappa$ is continuous and injective),
	it follows that the limit $\lim_{i\in I}\ms C_\kappa^{-1}(q_i)$
	exists in $\mc E_N\subseteq\bb E_{\leq\kappa}$.
\end{proof}

\medskip\noindent
The following result is used in the proofs of \cref{L71,L91}.

\begin{proposition}\label{L78}
	Let $\kappa\in\bb N_0$ and $q_n\in\mc N_{\le\kappa}$, $n\in\bb N_0$.
	Assume that
	\begin{Enumerate}
	\item
		for each compact $K\subseteq\bb C^+$ with non-empty interior $O$
		there exists $m_K\in\bb N$ such that $q_n$ is analytic on $O$
		for all $n\ge m_K$;
	\item
		there exists $M\subseteq\bb C^+$ with accumulation point in $\bb C^+$
		such that $\lim\limits_{n\to\infty}q_n(z)$ exists for all $z\in M$.
	\end{Enumerate}
	Then there exists $\mr q\in\mc N_{\le\kappa}\cap\Hol(\bb C^+)$
	such that $\lim\limits_{n\to\infty}q_n=\mr q$ locally uniformly in $\bb C^+$.
	Here we understand locally uniform convergence in the space of meromorphic functions considered as analytic functions
	into the Riemann sphere. 
\end{proposition}

\begin{proof}
	Assumptions (i) and (ii) imply, in particular, that there
	exist $\kappa+1$ points $z_0,\ldots,z_\kappa\in\bb C^+$ such
	that $|q_n(z_i)|\le c$ for all $n\in\bb N$ and $i\in\{0,\ldots,\kappa\}$
	and some $c>0$.
	By \cite[Theorem~3.2]{langer.luger.matsaev:2011} there exist
	a subsequence $(q_{n_k})_{k\in\bb N}$, a set $P\subseteq\bb C^+$
	with $|P|\le\kappa$, and $\mr q\in\mc N_{\le\kappa}$ such that
	\[
		\lim_{k\to\infty}q_{n_k} = \mr q
		\qquad\text{locally uniformly on} \ \bb C^+\setminus P.
	\]
	Note that $\mr q$ is meromorphic on $\bb C^+$ because $\mr q\in\mc N_{\le\kappa}$.
	Let $w\in P$.  There exists a closed disc $K\subseteq\bb C^+$ around $w$
	with interior $O$ such that $\mr q$ is zero-free on $O\setminus\{w\}$.
	By assumption (i), $q_{n_k}$ is analytic on $O$ for all $k$ with $n_k\ge m_K$.
	The convergence of the logarithmic residue implies that $\mr q$ is analytic at $w$,
	and hence $\lim_{k\to\infty}q_{n_k}=\mr q$ locally uniformly on $O$.
	Since $w$ was arbitrary in $P$, this shows that $\mr q$ is analytic on $\bb C^+$
	and that $\lim_{k\to\infty}q_{n_k}=\mr q$ locally uniformly on $\bb C^+$.

	The above considerations can be done for every subsequence of $(q_n)$
	instead of $(q_n)$ itself.  Now assumption (ii) implies
	that $\lim_{n\to\infty}q_n=\mr q$ locally uniformly on $\bb C^+$.
\end{proof}

\medskip\noindent
We can now prove an analogue of the classical Grommer--Hamburger theorem for
regularised Cauchy transforms.

\begin{theorem}\label{L71}
	Let $\kappa\in\bb N_0$, let $(\mu_n,p_n)\in\bb E_{\le\kappa}$ for $n\in\bb N$,
	set $q_n\DE\ms C_\kappa[\mu_n,p_n]$, and let $\mr q\in\Hol(\bb C^+)$.
	Then the following three statements are equivalent:
	\begin{Enumerate}
	\item
		$\exists M\subseteq\bb C^+$ such that $M$ has an accumulation point in $\bb C^+$
		and that $\lim\limits_{n\to\infty}q_n(z)=\mr q(z)$ for all $z\in M$;
	\item
		$\lim\limits_{n\to\infty}q_n=\mr q$ locally uniformly on $\bb C^+$;
	\item
		$\exists(\mr\mu,\mr p)\in\bb E_{\le\kappa}$ such
		that $\mr q=\ms C_\kappa[\mr\mu,\mr p]$, $\lim\limits_{n\to\infty}p_n=\mr p$ and
	\begin{equation}\label{L72}
		\forall a,b\in\bb R\DP a<b,\, \mr\mu(\{a\})=\mr\mu(\{b\})=0
		\quad\implies\quad
		\lim_{n\to\infty}\mu_n((a,b)) = \mr\mu((a,b)).
	\end{equation}
	\end{Enumerate}
\end{theorem}

\begin{proof}
	The equivalence of (i) and (ii) follows directly from \cref{L9,L78}
	since $q_n$ is analytic on $\bb C^+$ for every $n\in\bb N$.

	Let us now prove the equivalence of (ii) and (iii).
	It follows from \cref{L13} that (ii) is equivalent to
	\begin{equation}\label{L76}
		\exists(\mr\mu,\mr p)\in\bb E_{\le\kappa}\DP
		\mr q = \ms C_\kappa[\mr\mu,\mr p],\,\lim_{n\to\infty}p_n=\mr p
	\end{equation}
	together with $\lim_{n\to\infty}\mu_n=\mr\mu$ w.r.t.\ $w^*$ in $C_0(\bb R,\omega_\kappa)'$.
	Now assume that \eqref{L76} holds. Then the convergence of $(p_n)$ implies that
	\[
		\|\mu_n\|_{C_0(\bb R,\omega_\kappa)'}
		= \int_{\bb R}\frac{\D\mu_n(t)}{(1+t^2)^{\kappa+1}}
		\le \frac{1}{(2\kappa+1)!}p_n^{(2\kappa+1)}(0)
		\le c
	\]
	for some $c>0$.
	Since $C_{00}(\bb R)$, which denotes the set of
	compactly supported continuous functions on $\bb R$, is dense
	in $C_0(\bb R,\omega_\kappa)$,
	the relation $\mu_n\to\mr\mu$ w.r.t.\ $w^*$ in $C_0(\bb R,\omega_\kappa)'$
	is equivalent to $\mu_n\to\mr\mu$ w.r.t.\ $w^*$ in $C_{00}(\bb R)'$.
	The latter relation is equivalent to \eqref{L72}
	by the portmanteau-type theorem \cite[Theorem~1]{barczy.pap:2006}.
\end{proof}

\begin{remark}\label{L52}
\rule{0ex}{1ex}
\begin{Enumerate}
\item
	Let us make the connection with the original formulation of the Grommer--Hamburger theorem;
	see, e.g.\ \cite[\S48]{wintner:1929}.
	The latter is about Cauchy transforms, \eqref{L4}, of finite measures,
	and states the following:
	let $(\mu_n)_{n\in\bb N}$ be a sequence of finite measures on $\bb R$ whose total variations
	are uniformly bounded and let $\mr q\in\Hol(\bb C^+)$;
	then the following statements are equivalent:
	\begin{Enumeratealph}
	\item
		$\lim\limits_{n\to\infty}C[\mu_n](z) = \mr q(z)$ for all $z\in\bb C^+$;
	\item
		there exists a finite measure $\mr\mu$ such that $\mr q=C[\mr\mu]$ and \eqref{L72} holds.
	\end{Enumeratealph}
\item
	There has been some confusion about the formulation of the Grommer--Hamburger theorem.
	The condition in \eqref{L72} says that $(\mu_n)_{n\in\bb N}$ converges vaguely,
	i.e.\ w.r.t.\ $w^*$ in $C_{00}(\bb R)'$.
	However, at some places in the literature
	it is claimed that (a) in item (i) implies that $(\mu_n)_{n\in\bb N}$ converges weakly,
	i.e.\ w.r.t.\ $w^*$ in $C_{\rm b}(\bb R)'$,
	where $C_{\rm b}(\bb R)$ is the space of bounded continuous functions.
	The example $\mu_n=\delta_n$, where $\delta_n$ denotes the Dirac measure at $n$,
	shows that this is not true: $\lim_{n\to\infty}C[\delta_n]=0$
	locally uniformly, but $\lim_{n\to\infty}\delta_n=0$ only vaguely and not weakly;
	in particular, mass is lost.
	Note that, by the portmanteau theorem (see, e.g.\ \cite[Theorem~13.16]{klenke:2020}),
	a sequence of uniformly bounded measures $(\mu_n)$ on $\bb R$ converges weakly
	to a measure $\mr\mu$ if and only if it converges vaguely to $\mr\mu$
	and $\lim_{n\to\infty}\mu_n(\bb R)=\mr\mu(\bb R)$.
	See also the discussion in \cite{geronimo.hill:2003}.
\item
	In the original Grommer--Hamburger theorem one needs the a priori assumption
	that the total variations are uniformly bounded.
	On the other hand, in \cref{L71} the integrals
	$\int_{\bb R}(1+t^2)^{-(\kappa+1)}\D\mu_n(t)$ are automatically
	bounded by \eqref{L77} and the convergence of the polynomials $(p_n)$.
	Consider also the following example: let $\mu_n=n^2\delta_n$ and $p_n(z)=\frac{n^2}{1+n^2}z$.
	By \eqref{L86} we have
	\begin{align*}
		\ms C_0[\mu_n,p_n](z) &= p_n(z) + \int_{\bb R}\Bigl(\frac{1}{t-z}-\frac{t+z}{1+t^2}\Bigr)
		\DD\mu_n(t)
		\\
		&= \widetilde C[\mu_n](z)
		= n^2\Bigl(\frac{1}{n-z}-\frac{n}{1+n^2}\Bigr)
		\to z
	\end{align*}
	locally uniformly as $n\to\infty$.  Note that the limit function
	belongs to $\ran\ms C_0$ but is not the Cauchy transform of a finite measure.
\item
	It follows from \cite[Theorem~1]{barczy.pap:2006} that \eqref{L72}
	is equivalent to the following condition:
	\[
		\text{for every bounded Borel set $A$ with} \ \mr\mu(\partial A)=0:
		\quad
		\lim_{n\to\infty}\mu_n(A) = \mr\mu(A).
	\]
\item
	Under the additional assumption that $\mr q\in\mc N_\kappa^{(\infty)}$, i.e.\
	$\mr q$ has the same number of negative squares as $q_n$,
	the implication (ii)\,$\Rightarrow$\,\eqref{L72} in \cref{L71}
	can also be deduced from \cite[Corollary~3.1]{langer.luger.matsaev:2011};
	see also \cite[Lemma~3.7]{langer.woracek:sinham}.
\end{Enumerate}
\end{remark}

\section{A Tauberian theorem}
\label{L112}

In this section we prove a Tauberian theorem for the transform $\ms C_\kappa$
where the asymptotic behaviour of the measure $\mu$ towards infinity can be derived
from the asymptotic behaviour of the function $q=\ms C_\kappa[\mu,p]$ at infinity.
As mentioned in the Introduction, there is a wide range of Tauberian theorems 
for Stieltjes transforms, where the measure is only supported on the positive half-line.
Surprisingly, it seems there is much less known for Cauchy integrals,
where the measure is allowed to be supported on the whole real line.
One result, which has been frequently cited, is claimed in \cite[Theorem~7.5]{bennewitz:1989}.
The proof given in that paper contains a mistake\footnote{%
In particular, in that paper it is claimed that the relation $\int_{\bb R}(t-z)^{-3}\D\tau(t)=0$
for all $z\in\bb C^+$ for a non-decreasing function $\tau$ on $\bb R$ implies that $\tau$ 
is constant at all points of continuity, which is not true as the example $\tau(t)=t$ shows.}.
Fortunately, the result itself turns out to be true.
In this section we provide a simple and conclusive argument which allows us, at the same time,
to drop one assumption made in \cite{bennewitz:1989} and to generalise it
to higher-order regularised Cauchy transforms.

\Cref{L16} contains the above mentioned Tauberian theorem for the
transform $q=\ms C_\kappa[\mu,p]$.
In most cases the asymptotic behaviour of the measure $\mu$ can be determined
independently for the positive and the negative real axis; see \eqref{L92} and \eqref{L93}.
The assumption about the asymptotic behaviour of $q$ at infinity can be formulated
in different ways.
\Cref{L91} shows the equivalence of these assumptions, where conditions (i) and (ii)
are relatively minimal assumptions.  In particular, (ii) says that, along one ray
towards infinity, $q$ behaves like a constant times a regularly varying function;
for the latter notion see Appendix~\ref{L115}.
Note that we prove \cref{L91} for the larger class $\mc N_{\le\kappa}$ instead
of the class $\mc N_{\le\kappa}^{(\infty)}=\ran\ms C_\kappa$.

The following functions play an important role in the current section,
namely as limiting functions of rescalings of a given $q$:
\begin{equation}\label{L87}
	Q_{\alpha,\omega}(z) \DE i\omega\Bigl(\frac zi\Bigr)^\alpha,
	\qquad z\in\bb C^+,
\end{equation}
where $\alpha\in\bb R$ and $\omega\in\bb C\setminus\{0\}$.
Note that, by the normalisation of the power function 
at the end of the Introduction, we have $Q_{\alpha,\omega}(i)=i\omega$.

\begin{theorem}\label{L91}
	Let $\kappa\in\bb N_0$ and $q\in\mc N_{\le\kappa}$,
	and let $f:[x_0,\infty)\to(0,\infty)$ with $x_0>0$ be measurable.
	Then the following statements are equivalent:
	\begin{Enumerate}
	\item
		there exists $M\subseteq\bb C^+$ with an accumulation point in $\bb C^+$
		such that
		\begin{equation}\label{L96}
			\forall z\in M\DP
			\lim_{r\to\infty}\frac{q(rz)}{f(r)} \quad\text{exists and is non-zero;}
		\end{equation}
	\item
		$f$ is regularly varying and there exists $z_0\in\bb C^+$ such that
		\begin{equation}\label{L63}
			\lim_{r\to\infty}\frac{q(rz_0)}{f(r)} \quad\text{exists and is non-zero;}
		\end{equation}
	\item
		$f$ is regularly varying with index $\alpha\in[-2\kappa-1,2\kappa+1]$
		and there exists $\omega\in\bb C\setminus\{0\}$ such that
		\[
			\lim_{r\to\infty}\frac{q(rz)}{f(r)} = Q_{\alpha,\omega}(z)
			\qquad\text{locally uniformly for} \ z\in\bb C^+;
		\]
		moreover, there exists $\kappa'\in\{0,\ldots,\kappa\}$
		such that $Q_{\alpha,\omega}\in\mc N_{\kappa'}$ and hence
		\begin{equation}\label{L102}
			\big|\arg\bigl((-1)^{\kappa'}\omega\bigr)\big|
			\le \frac\pi 2\Bigl(1-\big||\alpha|-2\kappa'\big|\Bigr).
		\end{equation}
	\end{Enumerate}
	If \textup{(i)--(iii)} are satisfied and $\xi$ denotes the limit
	in \eqref{L63}, then $\omega=\frac{\xi}{i}\bigl(\frac{z_0}{i}\bigr)^{-\alpha}$.
\end{theorem}

\medskip\noindent
Note that in the case when $\kappa=0$ we must have $\kappa'=0$ in (iii) and hence \eqref{L102}
reduces to $|\arg\omega|\le\frac{\pi}{2}(1-|\alpha|)$.

\begin{theorem}\label{L16}
	Let $\kappa\in\bb N_0$ and $(\mu,p)\in\bb E_{\le\kappa}$,
	set $q\DE\ms C_\kappa[\mu,p]$, and let $f:[x_0,\infty)\to(0,\infty)$ with $x_0>0$ be measurable.
	Assume that the equivalent conditions \textup{(i)}--\textup{(iii)} in \cref{L91} hold.
	Then $\alpha\in[-1,2\kappa+1]$ and $Q_{\alpha,\omega}\in\ran\ms C_\kappa$.
	Let $\mu_{\alpha,\omega}$ be the measure component
	of $\ms C_\kappa^{-1}Q_{\alpha,\omega}$.  For all $a,b\in\bb R$
	with $a<b$ and $\mu_{\alpha,\omega}(\{a\})=\mu_{\alpha,\omega}(\{b\})=0$,
	we have
	\begin{equation}\label{L98}
		\lim_{r\to\infty}\frac{1}{rf(r)}\mu\bigl((ra,rb)\bigr)
		= \mu_{\alpha,\omega}\bigl((a,b)\bigr).
	\end{equation}
	In particular, if $\alpha>-1$, then
	\begin{align}
		\lim_{r\to\infty}\frac{1}{rf(r)}\mu\bigl((0,r)\bigr)
		&= \frac{1}{\pi}\cdot\frac{|\omega|}{\alpha+1}
		\cos\Bigl(\frac{\alpha\pi}{2}-\arg\omega\Bigr),
		\label{L92}
		\\[1ex]
		\lim_{r\to\infty}\frac{1}{rf(r)}\mu\bigl((-r,0)\bigr)
		&= \frac{1}{\pi}\cdot\frac{|\omega|}{\alpha+1}
		\cos\Bigl(\frac{\alpha\pi}{2}+\arg\omega\Bigr);
		\label{L93}
	\end{align}
	if $\alpha=-1$, then $\omega>0$ and
	\begin{equation}\label{L101}
		\lim_{r\to\infty}\frac{1}{rf(r)}\mu\bigl((-r,r)\bigr) = \omega.
	\end{equation}
\end{theorem}

\begin{remark}\label{L108}
\rule{0ex}{1ex}
\begin{Enumerate}
\item
	Any asymmetry of the limits on the right-hand sides of \eqref{L92} and \eqref{L93}
	can be seen from $\omega$.  To this end, assume that $\alpha>-1$ and
	write $\alpha=2m+\alpha_0$ with $m\in\bb N_0$ and $|\alpha_0|\le1$;
	if $\alpha\in2\bb N_0+1$, choose $m$ such that $(-1)^m\omega>0$.
	The limits on the right-hand sides of \eqref{L92} and \eqref{L93} can be
	rewritten as follows:
	\begin{align*}
		c_\pm \hspace*{0.72ex}&\hspace*{-0.72ex}\DE \frac{1}{\pi}\cdot\frac{|\omega|}{\alpha+1}
		\cos\Bigl(\frac{\alpha\pi}{2}\mp\arg\omega\Bigr)
		\\
		&= \frac{1}{\pi}\cdot\frac{|\omega|}{\alpha+1}
		\begin{cases}
			\displaystyle \cos\Bigl(\frac{\alpha_0\pi}{2}\mp\arg\omega\Bigr), & m \ \text{even},
			\\[2ex]
			\displaystyle \cos\Bigl(\frac{\alpha_0\pi}{2}+\pi\mp\arg\omega\Bigr), & m \ \text{odd},
		\end{cases}
		\\
		&= \frac{1}{\pi}\cdot\frac{|\omega|}{\alpha+1}
		\cos\Bigl(\frac{\alpha_0\pi}{2}\mp\arg\bigl((-1)^m\omega\bigr)\Bigr).
	\end{align*}
	It follows from \cref{L27} below and its proof that $m=\kappa'$ and hence,
	by \eqref{L102}, that
	$\big|\frac{\alpha_0\pi}{2}\mp\arg\bigl((-1)^m\omega\bigr)\big|\le\frac{\pi}{2}$.
	From this the following equivalences follow easily,
	where we set $\psi\DE\arg((-1)^m\omega)$,
	\begin{alignat*}{3}
		c_+ &> c_- \quad &&\iff\quad && \alpha_0,\psi\ne0 \;\;\wedge\;\; \sgn\alpha_0 = \sgn\psi,
		\\[1ex]
		c_+ &= c_- \quad &&\iff\quad && \alpha_0=0 \;\;\vee\;\; \psi=0,
		\\[1ex]
		c_+ &< c_- \quad &&\iff\quad && \alpha_0,\psi\ne0 \;\;\wedge\;\; \sgn\alpha_0 = -\sgn\psi,
		\\[1ex]
		c_\pm &= 0 \quad &&\iff\quad &&
			\big|\tfrac{\alpha_0\pi}{2}\mp\psi\big| = \tfrac{\pi}{2},
		\\[1ex]
		c_+ &= c_- = 0 \quad &&\iff\quad && |\alpha_0|=1
			\;\;\vee\;\; \bigl(\alpha_0=0 \;\wedge\; |\psi|=\tfrac{\pi}{2}\bigr).
	\end{alignat*}
	If $c_+\ne0$, then the function $r\mapsto\mu((0,r))$ is regularly varying
	with index $\alpha+1$; if $c_-\ne0$, then $r\mapsto\mu((-r,0))$ is
	regularly varying with index $\alpha+1$.
\item
	\Cref{L208} below shows that there are situations where (i)--(iii) in \cref{L91}
	are satisfied but none of $t\mapsto\mu((-t,t))$, $t\mapsto\mu((0,t))$,
	$t\mapsto\mu((-t,0))$ is regularly varying.
	See also \cref{L212}.
\end{Enumerate}
\end{remark}

\medskip\noindent
Before we prove Theorems~\ref{L91} and \ref{L16}, we need some lemmas.

\begin{lemma}\label{L27}
	Let $\alpha\in\bb R$ and $\omega\in\bb C\setminus\{0\}$ and let $Q_{\alpha,\omega}$
	be as in \eqref{L87}.
	\begin{Enumerate}
	\item
		We have $Q_{\alpha,\omega}\in\mc N_\kappa$ if and only if
		\begin{equation}\label{L65}
			\big|\arg\bigl((-1)^\kappa\omega\bigr)\big|
			\le \frac\pi 2\Bigl(1-\big||\alpha|-2\kappa\big|\Bigr).
		\end{equation}
		This is possible only when
		\begin{equation}\label{L81}
			\kappa =
			\begin{cases}
				\bigl\lfloor\frac{|\alpha|+1}{2}\bigr\rfloor
				&
				\text{if} \ \alpha\notin 2\bb Z+1,
				\\[1ex]
				\frac{|\alpha|+1}{2}
				&
				\text{if} \ \alpha\in 2\bb Z+1 \ \text{and} \ (-1)^{\frac{|\alpha|+1}{2}}\omega>0,
				\\[1ex]
				\frac{|\alpha|-1}{2}
				&
				\text{if} \ \alpha\in 2\bb Z+1 \ \text{and} \ (-1)^{\frac{|\alpha|+1}{2}}\omega<0.
			\end{cases}
		\end{equation}
		In particular, $\frac{|\alpha|-1}{2}\le\kappa\le\frac{|\alpha|+1}{2}$.
	\item
		Assume that $Q_{\alpha,\omega}\in\mc N_\kappa$, i.e.\ that \eqref{L65} is satisfied.
		Then $Q_{\alpha,\omega}\in\mc N_\kappa^{(\infty)}$
		if and only if $\alpha\ge-1$ and, in addition, $\omega>0$ in the case when $\alpha=-1$.
	\item
		Assume that $Q_{\alpha,\omega}\in\mc N_\kappa^{(\infty)}$ and let $\mu_{\alpha,\omega}$
		be the measure component of $\ms C_\kappa^{-1} Q_{\alpha,\omega}$.
		If $\alpha=-1$, then $\mu_{\alpha,\omega}=\omega\delta_0$, where $\delta_0$ is the
		Dirac measure at $0$.
		If $\alpha>-1$, then $\mu_{\alpha,\omega}$ is absolutely continuous w.r.t.\
		the Lebesgue measure and has density
		\[
			\frac{\D\mu_{\alpha,\omega}}{\D t}(t) = \frac{|\omega|}{\pi}|t|^\alpha
			\cos\Bigl(\frac{\alpha\pi}{2}-(\sgn t)\arg\omega\Bigr),
			\qquad \text{a.e.} \ t\in\bb R.
		\]
	\end{Enumerate}
\end{lemma}

\begin{proof}
	\phantom{}

(i)
	Write $\alpha=2m+\alpha_0$ with $m\in\bb Z$ and $|\alpha_0|\le1$
	(note that if $\alpha$ is an odd integer, then $m$ and $\alpha_0$ are not unique).
	Then
	\begin{equation}\label{L70}
		Q_{\alpha,\omega}(z) = i\omega\Bigl(\frac{z}{i}\Bigr)^{2m}\Bigl(\frac{z}{i}\Bigr)^{\alpha_0}
		= z^{2m}q_0(z)
	\end{equation}
	with
	\[
		q_0(z) = i(-1)^m\omega\Bigl(\frac{z}{i}\Bigr)^{\alpha_0}.
	\]
	Since the only generalised poles and zeros not of positive type
	(in the sense of \cite{krein.langer:1977}; see also \cite{langer:1986})
	can be $0$ and $\infty$,
	it follows from \cite[Corollary]{dijksma.langer.luger.shondin:2000}
	or \cite[Proposition~3.2 and Theorem~3.3]{derkach.hassi.snoo:1999}
	that $Q_{\alpha,\omega}\in\mc N_\kappa$ if and only if $|m|=\kappa$ and $q_0\in\mc N_0$.
	Determining the sector onto which $\bb C^+$ is mapped under $q_0$
	one can easily show that $q_0\in\mc N_0$ if and only if
	\[
		\big|\arg\bigl((-1)^m\omega\bigr)\big| \le \frac{\pi}{2}\bigl(1-|\alpha_0|\bigr).
	\]
	Since $|\alpha_0|=\big||\alpha|-|2m|\big|=\big||\alpha|-2\kappa\big|$, the equivalence
	of $Q_{\alpha,\omega}\in\mc N_\kappa$ and \eqref{L65} follows.
	The formula for $\kappa$ can be derived easily from \eqref{L65}.

(ii)
	It follows from the factorisation \eqref{L70}
	and \cite[Corollary]{dijksma.langer.luger.shondin:2000} that $q\in\mc N_\kappa^{(\infty)}$
	if and only if $m\ge0$.

(iii)
	If $\alpha=-1$ and $\omega>0$, then $Q_{\alpha,\omega}(z)=-\frac{\omega}{z}$ and
	hence $Q_{\alpha,\omega}\in\mc N_0=\mc N_0^{(\infty)}$
	and $\mu_{\alpha,\omega}=\omega\delta_0$.
	Assume that $\alpha>-1$.  The function $Q_{\alpha,\omega}$ can be extended
	to a continuous function on $(\bb C^+\cup\bb R)\setminus\{0\}$,
	and, for $t>0$, we have
	\begin{align*}
		\Im Q_{\alpha,\omega}(\pm t)
		&= \Im\bigl[i\omega(\mp it)^\alpha\bigr]
		= \Im\Bigl[|\omega|t^\alpha e^{i(\frac{\pi}{2}+\arg\omega\mp\alpha\frac{\pi}{2})}\Bigr]
		\\
		&= |\omega|t^\alpha\sin\Bigl(\frac{\pi}{2}+\arg\omega\mp\frac{\alpha\pi}{2}\Bigr)
		= |\omega|t^\alpha\cos\Bigl(\arg\omega\mp\frac{\alpha\pi}{2}\Bigr).
	\end{align*}
	Now the assertion follows from the Stieltjes inversion formula \eqref{L73}.
	Note that there is no point mass at $0$ if $\alpha>-1$.
\end{proof}

\begin{remark}\label{L82}
	It follows from \eqref{L65} and \eqref{L81} that $Q_{\alpha,\omega}\in\mc N_{<\infty}$
	if and only if
	\begin{alignat}{2}
		& \alpha\notin 2\bb Z+1
		\quad\text{and}\quad &
		& \big|\arg\bigl((-1)^{\lfloor\frac{|\alpha|+1}{2}\rfloor}\omega\bigr)\big|
		\le \frac{\pi}{2}\dist(\alpha,2\bb Z+1)
		\label{L89}
		\\[-1ex]
		& \hspace*{-10ex} \text{or}
		\nonumber
		\\[-1ex]
		& \alpha\in 2\bb Z+1
		\quad\text{and}\quad &
		& \omega\in\bb R.
		\label{L90}
	\end{alignat}
\end{remark}

\medskip\noindent
The significance of regular variation is that having a regularly varying asymptotics
for $q(rz)$ for one single point $z$ already suffices to get locally uniform asymptotics
depending on $z$ as a power.  The reason for this is the multiplicative nature
of the argument in $q(rz)$.
In the following lemma we use the standard notation $z_0M=\{z_0z\DS z\in M\}$
with $M\subseteq\bb R$.

\begin{lemma}\label{L24}
	Let $z_0\in\bb C^+$, $r_0>0$ and let $q:z_0(r_0,\infty)\to\bb C\setminus\{0\}$
	be a continuous function.
	Further, let $f:[x_0,\infty)\to(0,\infty)$ with $x_0>0$ be a measurable function
	and let $B\subseteq(r_0,\infty)$ be a set with positive Lebesgue measure.
	Assume that
	\[
		\forall z\in z_0B\DP
		\lim_{r\to\infty}\frac{q(rz)}{f(r)} \quad\text{exists and is non-zero}.
	\]
	Then there exist $\alpha\in\bb R$ and $\omega\in\bb C\setminus\{0\}$
	such that $f$ is regularly varying with index $\alpha$ and
	\begin{equation}\label{L88}
		\forall z\in z_0(0,\infty)\DP
		\lim_{r\to\infty}\frac{q(rz)}{f(r)} = Q_{\alpha,\omega}(z),
	\end{equation}
	where $Q_{\alpha,\omega}$ is as in \eqref{L87}.
\end{lemma}

\begin{proof}
	Set $\mr q(z)\DE\lim_{r\to\infty}\frac{q(rz)}{f(r)}$, $z\in z_0B$.
	Choose $s_0\in B$.  For every $\lambda\in\frac{1}{s_0}B$
	we have $\lambda s_0z_0\in z_0B$ and $s_0z_0\in z_0B$ and hence
	\[
		\lim_{r\to\infty}\frac{q(r\cdot\lambda s_0z_0)}{f(r)}
		= \mr q(\lambda s_0z_0),
		\qquad
		\lim_{r\to\infty}\frac{q(\lambda r\cdot s_0z_0)}{f(\lambda r)}
		= \mr q(s_0z_0).
	\]
	Taking quotients of these equations we obtain
	\[
		\lim_{r\to\infty}\frac{f(\lambda r)}{f(r)}
		= \frac{\mr q(\lambda s_0z_0)}{\mr q(s_0z_0)}.
	\]
	Since the set $\frac{1}{s_0}B$ has positive measure,
	the Characterisation Theorem \cite[Theorem~1.4.1]{bingham.goldie.teugels:1989}
	yields that $f$ is regularly varying with index, say, $\alpha\in\bb R$.
%	the right-hand side of \eqref{tmp1} equals $\lambda^\alpha$
	Hence, for every $\lambda>0$ we have
	\[
		\lim_{r\to\infty}\frac{q(r\cdot\lambda s_0z_0)}{f(r)}
		= \lim_{r\to\infty}\frac{q(\lambda r\cdot s_0z_0)}{f(\lambda r)}
		\cdot\lim_{r\to\infty}\frac{f(\lambda r)}{f(r)}
		= \mr q(s_0z_0)\lambda^\alpha.
	\]
	Replacing $\lambda s_0z_0$ by $z$ we obtain that, for every $z\in z_0(0,\infty)$,
	\[
		\lim_{r\to\infty}\frac{q(rz)}{f(r)} = \mr q(s_0z_0)\Bigl(\frac{z}{s_0z_0}\Bigr)^\alpha,
	\]
	which implies \eqref{L88}.
\end{proof}

\begin{lemma}\label{L17}
	Let $q\in\ran\ms C_\kappa$, $r>0$, and consider the
	function $q_{\langle r\rangle}(z)\DE q(rz)$.
	Then the following statements hold.
	\begin{Enumerate}
	\item
		$q_{\langle r\rangle}\in\ran\ms C_\kappa$.
	\item
		The measure components $\mu$ and $\mu_{\langle r\rangle}$ of $\ms C_\kappa^{-1}(q)$
		and $\ms C_\kappa^{-1}(q_{\langle r\rangle})$, respectively, are related by
		\[
			\mu_{\langle r\rangle}=\tfrac{1}{r}\Sigma^r_*\mu,
		\]
		where $\Sigma^r_*\mu$ is the push-forward of $\mu$ under the
		map $\Sigma^r\DF t\mapsto \frac{1}{r}t$,
		i.e.\ $\mu_{\langle r\rangle}(M)=\frac{1}{r}\mu(rM)$ for a measurable
		set $M\subseteq\bb R$.
	\end{Enumerate}
\end{lemma}

\begin{proof}
	The statement in (i) is obvious from \cref{L9}.
	Write $q=\ms C_\kappa[\mu,p]$ and
	$q_{\langle r\rangle}=\ms C_\kappa[\mu_{\langle r\rangle},p_{\langle r\rangle}]$.
	Making a change of variable ($t=rs$) we obtain
	\begin{align*}
		q_{\langle r\rangle}(z) &= q(rz)
		= p(rz)+\bigl(1+(rz)^2\bigr)^{\kappa+1}
		\int_{\bb R}\frac{1}{t-rz}\cdot\frac{\D\mu(t)}{(1+t^2)^{\kappa+1}}
		\\[1ex]
		&= p(rz)+\bigl(1+(rz)^2\bigr)^{\kappa+1}\cdot\int_{\bb R}\frac{1}{s-z}
		\cdot\frac{\frac{1}{r}\DD(\Sigma^r_*\mu)(s)}{(1+(rs)^2)^{\kappa+1}}\,.
	\end{align*}
	Extend $q_{\langle r\rangle}$ to $\bb C^+\cup\bb C^-$ by symmetry:
	$q_{\langle r\rangle}(z)\DE\qu{q_{\langle r\rangle}(\qu z)}$,
	$z\in\bb C^-$, and let $\Gamma_{\varepsilon}^\delta$ be the path in \eqref{L75}.
	The Stieltjes inversion formula \eqref{L73} and
	the Stieltjes--Liv\v{s}ic inversion formula \eqref{L74} yield that,
	for all $\alpha,\beta\in\bb R$ with $\alpha<\beta$, we have
	\begin{align*}
		\mu_{\langle r\rangle}\bigl((\alpha,\beta)\bigr)
		&= \lim_{\delta\downarrow 0}\lim_{\varepsilon\downarrow 0}
		\frac{1}{\pi}\int_{\alpha+\delta}^{\beta-\delta} \Im q_{\langle r\rangle}(t+i\varepsilon)\DD t
		= \lim_{\delta\downarrow 0}\lim_{\varepsilon\downarrow 0}\frac{-1}{2\pi i}
		\int_{\Gamma_\varepsilon^\delta}q_{\langle r\rangle}(z)\DD z
		\\[1ex]
		&= \lim_{\delta\downarrow 0}\lim_{\varepsilon\downarrow 0}\frac{-1}{2\pi i}
		\int_{\Gamma_\varepsilon^\delta}\bigl(1+(rz)^2\bigr)^{\kappa+1}
		\biggl[\int_{\bb R}\frac 1{s-z}\cdot
		\frac{\frac 1r\DD(\Sigma^r_*\mu)(s)}{(1+(rs)^2)^{\kappa+1}}\biggr]\DD z
		\\[1ex]
		&= \int_{(\alpha,\beta)}\big(1+(rs)^2\big)^{\kappa+1}\cdot
		\frac{\frac 1r\DD(\Sigma^r_*\mu)(s)}{(1+(rs)^2)^{\kappa+1}}
		= \tfrac 1r(\Sigma^r_*\mu)\bigl((a,b)\bigr).
		\qedhere
	\end{align*}
\end{proof}

\begin{proof}[Proof of \cref{L91}.]
	The implication (iii)\,$\Rightarrow$\,(ii) is trivial.

	Next let us show the implication (ii)\,$\Rightarrow$\,(i).
	For each $s>0$ the limit
	\[
		\lim_{r\to\infty}\frac{q(r\cdot sz_0)}{f(r)}
		= \lim_{r\to\infty}\biggl(\frac{q(rs\cdot z_0)}{f(rs)}\cdot\frac{f(rs)}{f(r)}\biggr)
		= \lim_{t\to\infty}\frac{q(tz_0)}{f(t)}\cdot\lim_{r\to\infty}\frac{f(rs)}{f(r)}
	\]
	exists and is non-zero.  Hence (i) is satisfied with $M=\{sz_0\DS s\in(0,\infty)\}$.

	Finally, we prove the implication (i)\,$\Rightarrow$\,(iii).
	Let $(r_n)_{n\in\bb N}$ be an arbitrary sequence of positive numbers
	with $r_n\to\infty$, and set
	\begin{equation}\label{L97}
		q_n(z) \DE \frac{q(r_nz)}{f(r_n)}, \qquad z\in\bb C^+.
	\end{equation}
	It is easy to see that $q_n\in\mc N_{\le\kappa}$.
	Since $q$, as a function from $\mc N_{\le\kappa}$, has only finitely many
	poles in $\bb C^+$, assumption (i) in \cref{L78} is satisfied.
	By \eqref{L96} also assumption (ii) in \cref{L78} is fulfilled.
	Hence, the latter proposition implies that there
	exists $\mr q\in\mc N_{\le\kappa}\cap\Hol(\bb C^+)$
	such that $\lim_{n\to\infty}q_n=\mr q$ locally uniformly in $\bb C^+$.
	Since the sequence $(r_n)_{n\in\bb N}$ was arbitrary, it follows again from \eqref{L96} that
	\[
		\lim_{r\to\infty}\frac{q(rz)}{f(r)} = \mr q(z)
	\]
	locally uniformly for $z\in\bb C^+$.  Now \cref{L24} implies that
	there exist $\alpha\in\bb R$ and $\omega\in\bb C\setminus\{0\}$ such that
	$f$ is regularly varying with index $\alpha$ and that $\mr q=Q_{\alpha,\omega}$.
	Since $\mc N_{\le\kappa}$ is closed under locally uniform convergence,
	we have $Q_{\alpha,\omega}\in\mc N_{\kappa'}$ with some $\kappa'\le\kappa$.
	By \cref{L27} this shows that $|\alpha|\le2\kappa'+1\le2\kappa+1$
	and that \eqref{L102} holds.
\end{proof}

\begin{proof}[Proof of \cref{L16}.]
	Let $(r_n)_{n\in\bb N}$ be an arbitrary sequence of positive numbers
	with $r_n\to\infty$, and define $q_n$ as in \eqref{L97}.
	It follows from \cref{L17} that $q_n\in\ran\ms C_\kappa$ and that
	the measure component $\mu_n$ of $\ms C_\kappa^{-1}q_n$ satisfies
	\begin{equation}\label{L99}
		\mu_n\bigl((a,b)\bigr) = \frac{1}{r_nf(r_n)}\mu\bigl((r_na,r_nb)\bigr)
	\end{equation}
	for all $a,b\in\bb R$ with $a<b$.
	\Cref{L71} implies that $Q_{\alpha,\omega}\in\ran\ms C_\kappa=\mc N_{\le\kappa}^{(\infty)}$,
	which, by \cref{L27}\,(ii) shows that $\alpha\ge-1$.
	Let $\mu_{\alpha,\omega}$ be the measure component of $\ms C_\kappa^{-1}Q_{\alpha,\omega}$.
	Further, let $a,b\in\bb R$ be such that $a<b$ and
	$\mu_{\alpha,\omega}(\{a\})=\mu_{\alpha,\omega}(\{b\})=0$.
	It follows from \cref{L71} that
	\begin{equation}\label{L100}
		\lim_{n\to\infty}\mu_n\bigl((a,b)\bigr) = \mu_{\alpha,\omega}\bigl((a,b)\bigr).
	\end{equation}
	Since the sequence $(r_n)_{n\in\bb N}$ was arbitrary, relations \eqref{L99} and \eqref{L100}
	imply \eqref{L98}.
	We obtain from \cref{L27}\,(iii) that
	\[
		\mu_{\alpha,\omega}\bigl(\pm(0,1)\bigr) = \frac{1}{\pi}\cdot\frac{|\omega|}{\alpha+1}
		\cos\Bigl(\frac{\alpha\pi}{2}\mp\arg\omega\Bigr)
	\]
	if $\alpha>-1$ and $\mu_{\alpha,\omega}((-1,1))=\omega$ if $\alpha=-1$.
	This, combined with \eqref{L98}, yields \eqref{L92}--\eqref{L101}.
\end{proof}

\section{The Abelian direction}
\label{L113}

In this section we consider Abelian theorems, i.e.\ we study the asymptotic behaviour
of $q$ at infinity using some knowledge about the asymptotic behaviour of the
distribution functions $t\mapsto\mu([0,t))$ and $t\mapsto\mu((-t,0))$.
Some of the theorems contain also a Tauberian direction,
which complement the results in Section~\ref{L112}.
The main tool is Karamata's theorem about Stieltjes transforms of measures
that are supported on the half-axis $[0,\infty)$.
We follow the common lines to pass from unilateral to bilateral theorems,
e.g.\ \cite{pleijel:1963a, selander:1963},
and represent imaginary and real parts as Stieltjes transforms.
To this end we use the push-forward measure $\mu_*$ of $\mu$
under the map $t\mapsto t^2$, which satisfies
\begin{equation}\label{L2}
	\mu_*((-\infty,0))=0 \qquad\text{and}\qquad
	\mu_*([0,t^2))=\mu((-t,t)), \quad t>0.
\end{equation}
We shall often use a substitution to change between $\mu$ and $\mu_*$;
let us note that, for a non-negative, measurable function $h$ on $[0,\infty)$, we have 
$\int_{[0,\infty)}h(s)\DD\mu_*(s)=\int_{\bb R}h(t^2)\DD\mu(t)$.

In order to apply Karamata's theorem in an effective way, 
we need a finer classification of the growth properties of the positive measure $\mu$,
namely, let us set
\begin{equation}\label{L215}
	\ms p(\mu) \DE \inf\biggl\{n\in\bb N\DS
	\int_{\bb R}\frac{\D\mu(t)}{(1+|t|)^n}<\infty\biggr\}
	\in \bb N\cup\{\infty\}.
\end{equation}

In the Abelian theorems we often assume that the symmetrised distribution function
$t\mapsto\mu((-t,t))$ is regularly varying.
The following lemma can be used to give a different characterisation of $\ms p(\mu)$ 
and its relation to the index of the regularly varying distribution function.

\begin{lemma}\label{L159}
	Let $\mu$ be a measure on $\bb R$.  For $\gamma>0$ we have
	\[
		\int_{\bb R}\frac{\D\mu(t)}{(1+|t|)^\gamma} < \infty
		\quad\Leftrightarrow\quad
		\int_1^\infty\frac{\mu((-t,t))}{t^{\gamma+1}}\DD t < \infty.
	\]
	If $t\mapsto\mu((-t,t))$ is regularly varying with index $\beta$,
	then $\ms p(\mu)$ is finite and $\beta\in[\ms p(\mu)-1,\ms p(\mu)]$.
\end{lemma}

\begin{proof}
	Let the measure $\mu_*$ be defined as in \eqref{L2}, let $\gamma>0$
	and define the measure $\nu$ on $[1,\infty)$ 
	such that $\nu((t,\infty))=t^{-\frac{\gamma}{2}}$.  It follows from \cref{L34} that 
	the following equivalences hold:
	\begin{align}
		\int_{\bb R}\frac{\D\mu(t)}{(1+|t|)^\gamma} < \infty
		\quad&\Leftrightarrow\quad
		\int_{[0,\infty)}\frac{\D\mu_*(s)}{(1+\sqrt{s})^\gamma} < \infty
		\quad\Leftrightarrow\quad
		\int_{[1,\infty)}s^{-\frac{\gamma}{2}}\DD\mu_*(s) < \infty
		\nonumber\\[1ex]
		&\Leftrightarrow\quad
		\int_1^\infty s^{-\frac{\gamma}{2}-1}\mu_*([1,s))\DD s < \infty
		\label{L174}
		\\[1ex]
		&\Leftrightarrow\quad
		\int_1^\infty\frac{\mu((-t,t))}{t^{\gamma+1}}\DD t < \infty,
		\nonumber
	\end{align}
	where in the last step we used the substitution $s=t^2$.
	
	Now assume that $t\mapsto\mu((-t,t))$ is regularly varying with index $\beta$.
	Then $s\mapsto\mu_*{[1,s)}$ is regularly varying with index $\frac{\beta}{2}$.
	It is clear that the integral in \eqref{L174} is finite if $\gamma$ is large enough,
	which shows that $\ms p(\mu)$ is finite.
	Further, the fact that the integral in \eqref{L174} is finite for $\gamma=\ms p(\mu)$
	and infinite for $\gamma=\ms p(\mu)-1$ (unless $\mu$ itself is finite) implies 
	that $\frac{\beta}{2}-\frac{\ms p(\mu)}{2}\le0$
	and $\frac{\beta}{2}-\frac{\ms p(\mu)-1}{2}\ge0$, which finishes the proof;
	cf.\ \cref{L33}.
	Note that, when $\mu$ is finite, then $\beta=0$ and $\ms p(\mu)=1$.
\end{proof}

\medskip\noindent
For the Abelian theorems we treat real and imaginary parts of $q(iy)$ separately
as they have different representations in terms of Stieltjes transforms.
This is done in the following two subsections.

\subsection{The imaginary part}

The imaginary part of $q(iy)$ is relatively well behaved as it can be written
in terms of one Stieltjes transform.
In order to apply Karamata's theorem, we choose $\kappa$ minimal in \eqref{L49}
for a given measure $\mu$.  
To this end, let us define
\begin{equation}\label{L60}
	\kappa(\mu) \DE \inf\biggl\{n\in\bb N_0\DS \int_{\bb R}\frac{\D\mu(t)}{(1+|t|)^{2n+2}}<\infty\biggr\} < \infty.
\end{equation}
Comparing \eqref{L60} with \eqref{L215} we can easily deduce that
\begin{equation}\label{L216}
	\kappa(\mu) = \Bigl\lfloor\frac{\ms p(\mu)-1}{2}\Bigr\rfloor.
\end{equation}
Throughout this section we suppose that the following assumption is satisfied.

\begin{assumption}\label{L171}
	Let $\mu$ be a measure on $\bb R$ such that $\ms p(\mu)<\infty$
	and set $\kappa\DE\kappa(\mu)$.
	Further, let $p\in\bb R[z]$ with $p(z)=c_{2\kappa+1}z^{2\kappa+1}+\ldots+c_0$
	such that $(\mu,p)\in\bb E_{\le\kappa}$, and set $q\DE\ms C_\kappa[\mu,p]$.
\end{assumption}

\medskip\noindent
Note that, by \cref{L8}\,(i), we have $c_{2\kappa+1}\ge\int_{\bb R}(1+t^2)^{-(\kappa+1)}\D\mu(t)$.
The two cases, equality and strict inequality, lead to different asymptotic behaviour
of $\Im q(iy)$ as $y\to\infty$, as the next proposition shows.
Naturally, the case when the polynomial dominates the integral is the simpler one.

\begin{proposition}\label{L172}
	Let $\mu$, $\kappa$, $p$ and $q$ be as in \cref{L171}.
	\begin{Enumerate}
	\item
		If $c_{2\kappa+1}>\int_{\bb R}(1+t^2)^{-(\kappa+1)}\D\mu(t)$, then
		\[
			q(iy) \sim i(-1)^\kappa\biggl(c_{2\kappa+1}-\int_{\bb R}\frac{\D\mu(t)}{(1+t^2)^{\kappa+1}}\biggr)y^{2\kappa+1},
			\qquad y\to\infty.
		\]
	\item
		If $c_{2\kappa+1}=\int_{\bb R}(1+t^2)^{-(\kappa+1)}\D\mu(t)$, then
		\begin{equation}\label{L210}
			|q(iy)| \ll y^{2\kappa+1}
			\qquad\text{and}\qquad
			\Im q(iy) \sim (-1)^\kappa y^{2\kappa+1}
			\int_{\bb R}\frac{1}{t^2+y^2}\cdot\frac{\D\mu(t)}{(1+t^2)^\kappa},
		\end{equation}
		as $y\to\infty$.
	\item
		If $\mu$ is an infinite measure, then
		\begin{equation}\label{L211}
			(-1)^\kappa\Im q(iy) \gg y^{2\kappa-1},
			\qquad y\to\infty.
		\end{equation}
	\end{Enumerate}
\end{proposition}

\begin{proof}
	\phantom{}
\begin{Steps}
\item
	From \eqref{L30} we obtain
	\begin{align}
		q(iy) &= p(iy) + iy(1-y^2)^\kappa\int_{\bb R}\frac{\D\mu(t)}{(1+t^2)^{\kappa+1}}
		+ (1-y^2)^\kappa\int_{\bb R}\biggl(\frac{1}{t-iy}-\frac{t}{1+t^2}\biggr)\frac{\D\mu(t)}{(1+t^2)^\kappa}
		\nonumber\\[1ex]
		&= i(-1)^\kappa\biggl(c_{2\kappa+1}-\int_{\bb R}\frac{\D\mu(t)}{(1+t^2)^{\kappa+1}}\biggr)y^{2\kappa+1}
		+ c_{2\kappa}(-1)^\kappa y^{2\kappa} + \BigO\bigl(y^{2\kappa-1}\bigr)
		\nonumber\\[1ex]
		&\quad + (-1)^\kappa\Bigl(y^{2\kappa}+\BigO\bigl(y^{2\kappa-2}\bigr)\Bigr)
		\int_{\bb R}\biggl(\frac{1}{t-iy}-\frac{t}{1+t^2}\biggr)\frac{\D\mu(t)}{(1+t^2)^\kappa}.
		\label{L186}
	\end{align}
	Together with \eqref{L183}, this proves the assertion in (i), relation \eqref{L211}
	when $c_{2\kappa+1}>\int_{\bb R}(1+t^2)^{-(\kappa+1)}\D\mu(t)$, 
	and the first relation in \eqref{L210} 
	when $c_{2\kappa+1}=\int_{\bb R}(1+t^2)^{-(\kappa+1)}\D\mu(t)$.

\item
	For the rest of the proof assume 
	that $c_{2\kappa+1}=\int_{\bb R}(1+t^2)^{-(\kappa+1)}\D\mu(t)$.
	If the measure $\mu$ is finite, then $\kappa=0$ and
	\[
		\Im q(iy) = \Im\int_{\bb R}\frac{1}{t-iy}\DD\mu(t) 
		= y\int_{\bb R}\frac{1}{t^2+y^2}\DD\mu(t).
	\]

\item	
	Let us now consider the case when $\mu$ is infinite.
	We can use \eqref{L186} to write the imaginary part as
	\[
		\Im q(iy) = (-1)^\kappa\Bigl(y^{2\kappa}+\BigO\bigl(y^{2\kappa-2}\bigr)\Bigr)
		\Im\biggl[\int_{\bb R}\biggl(\frac{1}{t-iy}-\frac{t}{1+t^2}\biggr)
		\frac{\D\mu(t)}{(1+t^2)^\kappa}\biggr]
		+ \BigO\bigl(y^{2\kappa-1}\bigr).
	\]
	By the definition of $\kappa$ we have 
	$\int_{\bb R}(1+t^2)^{-\kappa}\D\mu(t)=\infty$.
	Hence \eqref{L184} implies that
	\[
		\Im\biggl[\int_{\bb R}\biggl(\frac{1}{t-iy}-\frac{t}{1+t^2}\biggr)
		\frac{\D\mu(t)}{(1+t^2)^\kappa}\biggr] \gg \frac{1}{y},
	\]
	from which the second relation in \eqref{L210} follows.
\end{Steps}
\end{proof}

\medskip\noindent
In the following we assume that the leading asymptotics of $q(iy)$ is not given 
by a polynomial term but the measure $\mu$.
More precisely, we suppose that the following assumption is satisfied.

\begin{assumption}\label{L173}
	Let $\mu$, $\kappa$, $p$ and $q$ satisfy the conditions in \cref{L171}.
	Further, assume that
	\begin{equation}\label{L213}
		c_{2\kappa+1} = \int_{\bb R}\frac{\D\mu(t)}{(1+t^2)^{\kappa+1}}.
	\end{equation}
\end{assumption}

\begin{remark}\label{L204}
	Suppose that \cref{L171} holds and that $\kappa=\kappa(\mu)=0$, 
	with $\kappa$ defined in \eqref{L60}.  According to \eqref{L86} we can write
	\begin{equation}\label{L206}
		q(z) = \ms C_0[\mu,p] = \biggl(c_1-\int_{\bb R}\frac{\D\mu(t)}{1+t^2}\biggr)z
		+ c_0 + \int_{\bb R}\biggl(\frac{1}{t-z}-\frac{t}{1+t^2}\biggr)\DD\mu(t).
	\end{equation}
	\Cref{L173} is equivalent to the coefficient of $z$ in \eqref{L206}
	vanishing.  Hence with \cref{L173} being satisfied we have
	\begin{equation}\label{L207}
		q(z) = c_0 + \int_{\bb R}\biggl(\frac{1}{t-z}-\frac{t}{1+t^2}\biggr)\DD\mu(t)
		= c_0 + \widetilde C[\mu](z),
	\end{equation}
	where $\widetilde C[\mu]$ is defined in \eqref{L5}.
\end{remark}

\medskip\noindent
The next lemma shows that the imaginary part of $q(iy)$ can be written---at least
asymptotically---in terms of a Stieltjes transform.
For the definition of the Stieltjes transform see \eqref{L22}.

\begin{lemma}\label{L217}
	Let $p$, $\kappa$, $p$ and $q$ satisfy the conditions in \cref{L171} and \cref{L173}.
	Further, let $\mu_*$ be the push-forward measure of $\mu$ as in \eqref{L2}
	and define the measure $\tau_\kappa$ on $[0,\infty)$ by
	\begin{equation}\label{L218}
		\D\tau_\kappa(s) = \frac{\D\mu_*(s)}{(1+s)^\kappa},
		\qquad s\in[0,\infty).
	\end{equation}
	Then the Stieltjes transform $\ms S[\tau_\kappa]$ is well defined and
	\begin{equation}\label{L175}
		\Im q(iy) \sim (-1)^\kappa y^{2\kappa+1}\ms S[\tau_\kappa](y^2),
		\qquad y\to\infty.
	\end{equation}
\end{lemma}

\begin{proof}
	Since
	\[
		\int_{[0,\infty)}\frac{\D\tau_\kappa(s)}{1+s} 
		= \int_{[0,\infty)}\frac{\D\mu_*(s)}{(1+s)^{\kappa+1}}
		= \int_{\bb R}\frac{\D\mu(t)}{(1+t^2)^{\kappa+1}} < \infty,
	\]
	the Stieltjes transform $\ms S[\tau_\kappa]$ is well defined.
	It follows from \cref{L172}\,(ii) that
	\begin{align*}
		\Im q(iy) &\sim (-1)^\kappa y^{2\kappa+1}
		\int_{\bb R}\frac{1}{t^2+y^2}\cdot\frac{\D\mu(t)}{(1+t^2)^\kappa}
		\\[1ex]
		&= (-1)^\kappa y^{2\kappa+1}
		\int_{[0,\infty)}\frac{1}{s+y^2}\cdot\frac{\D\mu_*(s)}{(1+s)^\kappa}
		\\[1ex]
		&= (-1)^\kappa y^{2\kappa+1}
		\int_{[0,\infty)}\frac{1}{s+y^2}\DD\tau_\kappa(s)
		= (-1)^\kappa y^{2\kappa+1}\ms S[\tau_\kappa](y^2),
	\end{align*}
	which proves \eqref{L175}.
\end{proof}

\medskip\noindent
In the following theorem we prove that the imaginary part of $q(iy)$ is related
to the symmetrised distribution function $\mu((-t,t))$ of the measure $\mu$.
In most cases $|\Im q(iy)|$ is regularly varying if and only if $t\mapsto\mu((-t,t))$
is regularly varying.

\begin{theorem}\label{L45}
	Let $\mu$, $\kappa$, $p$ and $q$ satisfy the conditions in \cref{L171,L173}.
	Further, let $\beta\ge0$ and consider the following two statements:
	\begin{Enumeratealph}
	\item
		the symmetrised distribution function $t\mapsto\mu((-t,t))$ is regularly varying 
		with index $\beta$;
	\item
		the function $y\mapsto(-1)^\kappa\Im q(iy)$ is regularly varying with index $\beta-1$.
	\end{Enumeratealph}
	Then we have the following relations.
	\begin{Enumerate}
	\item
		The implication \textup{(a)}\,$\Rightarrow$\,\textup{(b)} holds.
	\item
		Unless $\kappa>0$ and $\beta=2\kappa$, also \textup{(b)}\,$\Rightarrow$\,\textup{(a)} holds.
	\item
		Assume that \textup{(a)} and \textup{(b)} are satisfied.  Then $\beta\in[2\kappa,2\kappa+2]$ and
		\begin{equation}\label{L46}
			\Im q(iy) \sim
			\begin{cases}
				\displaystyle \frac{\frac{\pi\beta}{2}}{\sin\frac{\pi\beta}{2}}\cdot\frac{\mu((-y,y))}{y},
				& \beta\in[0,\infty)\setminus\{2,4,\ldots\},
				\\[3ex]
				\displaystyle (-1)^\kappa\beta y^{\beta-1}\int_1^y\frac{\mu((-t,t))}{t^{\beta+1}}\DD t,
				& \beta=2\kappa \;\wedge\; \kappa>0,
				\\[3ex]
				\displaystyle (-1)^\kappa\beta y^{\beta-1}\int_y^\infty\frac{\mu((-t,t))}{t^{\beta+1}}\DD t,
				& \beta=2\kappa+2,
			\end{cases}
		\end{equation}
		as $y\to\infty$, where the first fraction in the first case on the right-hand side is 
		understood as $1$ when $\beta=0$.  In particular, if $\beta\in 2\bb N$, then
		\begin{equation}\label{L47}
			(-1)^\kappa\Im q(iy) \gg \frac{\mu((-y,y))}{y}.
		\end{equation}
	\end{Enumerate}
\end{theorem}

\begin{remark}\label{L48}
	In the situation of \cref{L45} assume that (a) is satisfied and that $\beta\in2\bb N$.  
	It follows from the definition of $\kappa$ in \eqref{L60} and from \cref{L159} that
	\begin{alignat*}{5}
		\beta &= 2\kappa 
		\quad &&\Leftrightarrow\quad & 
		\int_{\bb R}\frac{\D\mu(t)}{(1+|t|)^\beta} &= \infty
		\quad &&\Leftrightarrow\quad & 
		\int_1^\infty\frac{\mu((-t,t))}{t^{\beta+1}}\DD t &= \infty,
		\\[2ex]
		\beta &= 2\kappa+2 
		\quad &&\Leftrightarrow\quad & 
		\int_{\bb R}\frac{\D\mu(t)}{(1+|t|)^\beta} &< \infty
		\quad &&\Leftrightarrow\quad & 
		\int_1^\infty\frac{\mu((-t,t))}{t^{\beta+1}}\DD t &< \infty.
	\end{alignat*}
\end{remark}

\begin{proof}[Proof of \cref{L45}]
	Let $\mu_*$ be the push-forward measure of $\mu$ as in \eqref{L2}
	and define $\tau_\kappa$ as in \eqref{L218}.
	Then \eqref{L175} holds.
	We prove the theorem in several steps.

\begin{Steps}
\item
	Let us first consider the case when $\kappa=0$.  Then $\tau_0=\mu_*$ and,
	by \cref{L1}, we have the following equivalences:
	\begin{align*}
		\text{(a)} \quad &\Leftrightarrow\quad
		s\mapsto\mu_*([0,s)) \ \text{is regularly varying with index} \ \tfrac{\beta}{2}
		\\[1ex]
		&\Leftrightarrow\quad
		\ms S[\mu_*] \ \text{is regularly varying with index} \ \tfrac{\beta}{2}-1
		\\[1ex]
		&\Leftrightarrow\quad \text{(b)}.
	\end{align*}
	Assume now that (a) and (b) hold.  Another application of \cref{L1} implies
	that $\frac{\beta}{2}\in[0,1]$.  When $\beta\in[0,2)$, \cref{L3} yields
	\[
		\Im q(iy) = y\,\ms S[\mu_*](y^2)
		\sim y\frac{\frac{\pi\beta}{2}}{\sin\frac{\pi\beta}{2}}\cdot\frac{\mu_*([0,y^2))}{y^2}
		\sim \frac{\frac{\pi\beta}{2}}{\sin\frac{\pi\beta}{2}}\cdot\frac{\mu((-y,y))}{y}.
	\]
	When $\beta=2$, we obtain from \cref{L1} and the substitution $s=t^2$ that
	\[
		\Im q(iy) = y\,\ms S[\mu_*](y^2)
		\sim y\int_{y^2}^\infty\frac{\mu_*([0,s))}{s^2}\DD s
		= 2y\int_y^\infty\frac{\mu((-t,t))}{t^3}\DD t,
	\]
	which proves \eqref{L46} when $\kappa=0$.

\item
	Now let us consider the case when $\kappa>0$.  Set $h(s)\DE\frac{1}{(1+s)^\kappa}\sim s^{-\kappa}$,
	$s\to\infty$.  By the definition of $\kappa$ we have
	\begin{align*}
		\int_{[0,\infty)}h(s)\DD\mu_*(s) 
		&= \int_{[0,\infty)}\frac{\D\mu_*(s)}{(1+s)^\kappa} 
		= \int_{\bb R}\frac{\D\mu(t)}{(1+t^2)^\kappa} = \infty,
		\\[1ex]
		\int_{[0,\infty)}\frac{h(s)}{1+s}\DD\mu_*(s) 
		&= \int_{[0,\infty)}\frac{\D\mu_*(s)}{(1+s)^{\kappa+1}}
		= \int_{\bb R}\frac{\D\mu(t)}{(1+t^2)^{\kappa+1}} < \infty,
	\end{align*}
	which shows that \eqref{L118} with $\nu=\mu_*$ is satisfied.

\item
	Assume that $\kappa>0$ and (a) holds.  Then $s\mapsto\mu_*([0,s))$ 
	is regularly varying with index $\frac{\beta}{2}$.
	By \cref{L122}\,(i) with $\alpha=\frac{\beta}{2}$ and $\gamma=-\kappa$ 
	we have $\alpha+\gamma=\frac{\beta}{2}-\kappa\in[0,1]$ and hence $\beta>0$.  
	Therefore we can apply \cref{L122}\,(i) again to 
	obtain that $y\mapsto y^{2\kappa+1}\ms S[\tau_\kappa](y^2)$ is 
	regularly varying with index $2\kappa+1+2\bigl(\frac{\beta}{2}-\kappa-1\bigr)=\beta-1$,
	i.e.\ (b) holds.

	If $\beta\notin2\bb N$, then $\alpha+\gamma=\frac{\beta}{2}-\kappa\in(0,1)$ 
	and \eqref{L127} implies that
	\begin{align*}
		\Im q(iy) &\sim (-1)^\kappa y^{2\kappa+1}\ms S[\tau_\kappa](y^2)
		\\[1ex]
		&\sim (-1)^\kappa y^{2\kappa+1}
		\frac{\pi\frac{\beta}{2}}{\sin\bigl(\pi\bigl(\frac{\beta}{2}-\kappa\bigr)\bigr)}
		\cdot\bigl(y^2\bigr)^{-\kappa+1}\mu_*([0,y))
		\\[1ex]
		&= \frac{\frac{\pi\beta}{2}}{\sin\frac{\pi\beta}{2}}\cdot\frac{\mu((-y,y))}{y}.
	\end{align*}

	Now assume that $\beta\in2\bb N$.  The relation $\frac{\beta}{2}-\kappa\in[0,1]$ implies that
	either $\kappa=\frac{\beta}{2}$ or $\kappa=\frac{\beta}{2}-1$.
	Let us consider the former case; the other case is similar.
	It follows again from \eqref{L127} that, with the substitution $s=t^2$,
	\begin{align*}
		\Im q(iy) &\sim (-1)^\kappa y^{2\kappa+1}\ms S[\tau_\kappa](y^2)
		\\[1ex]
		&\sim (-1)^\kappa y^{2\kappa+1}\cdot\frac{\beta}{2}\cdot\frac{1}{y^2}
		\int_1^{y^2} s^{-\kappa-1}\mu_*([0,s))\DD s
		\\[1ex]
		&= (-1)^\kappa \beta y^{2\kappa-1}\int_1^y \frac{\mu_*([0,t^2))}{t^{2\kappa+1}}\DD t
		= (-1)^\kappa \beta y^{\beta-1}\int_1^y \frac{\mu((-t,t))}{t^{\beta+1}}\DD t,
	\end{align*}
	which proves \eqref{L46} also in the case $\kappa>0$.

\item
	Now assume that $\kappa>0$, $\beta\ne2\kappa$ and (b) holds.
	It follows from \eqref{L175} that $\ms S[\tau_\kappa]$ is regularly varying.
	Since $\alpha+\gamma=\frac{\beta}{2}-\kappa>0$, we can use \cref{L122}\,(ii)
	to deduce that $s\mapsto\mu_*([0,s))$ is regularly varying,
	which implies that (a) holds.

\item
	Finally, assume that $\beta\in2\bb N$.  
	Relation \eqref{L47} follows from \eqref{L46} and \cref{L36}.
\end{Steps}
\end{proof}

\medskip\noindent
The following example shows that the implication \textup{(b)}\,$\Rightarrow$\,\textup{(a)}
is, in general, not valid when $\beta=2\kappa>0$.

\begin{example}\label{L208}
	Let $\kappa\in\bb N$ and set $h(s)=\frac{1}{(1+s)^\kappa}$, $s\in[0,\infty)$.
	Choose the measures $\sigma$ and $\nu$ as in \cref{L202}, and let $\mu$
	be the symmetric measure on $\bb R$ such that $\mu_*=\nu$, 
	i.e.\ $\mu$ is the discrete measure with point masses
	\[
		\mu(\{e^{n/2}\}) = \mu(\{-e^{n/2}\}) = \nu(\{e^n\})
		= (1+e^n)^\kappa,
		\qquad n\in\bb N.
	\]
	According to \cref{L202} the distribution function $t\mapsto\mu((-t,t))=\nu([0,t^2))$
	is not regularly varying, which means that \textup{(a)} does not hold.
	On the other hand, $t\mapsto\tau_\kappa([0,t))=\sigma([0,t))$ is slowly varying, 
	again by \cref{L202}.  It follows from \cref{L1} that $\ms S[\tau_\kappa]$
	is regularly varying with index $-1$, and hence \textup{(b)} holds
	with $\beta=2\kappa$; see \eqref{L175}.
	Note that we have 
	\[
		\Im q(iy) \sim 2(-1)^\kappa y^{2\kappa-1}\log y,
		\qquad y\to\infty,
	\]
	by \eqref{L209}.
	
	Since $\mu$ is symmetric, we have $\Re q(iy)=0$ for $y>0$.
	Hence \eqref{L63} is satisfied with $z_0=i$ and $f(r)=r^{2\kappa-1}\log r$.
	\Cref{L16} implies that also \eqref{L92} and \eqref{L93} are satisfied.
	However, since $\omega=2(-1)^\kappa\in\bb R$ and $\alpha=2\kappa-1$ is odd,
	the right-hand sides of \eqref{L92} and \eqref{L93} vanish.
	This shows that, in some cases, one cannot use \cref{L16} to 
	deduce from the validity of (i)--(iii) in \cref{L91} 
	that $t\mapsto\mu((-t,t))$ is regularly varying.
\end{example}

\begin{example}\label{L64}
	Let $a,b\ge0$ with $a\ne b$, and consider the function
	\[
		q(z) = a\log z-b\log(-z),
	\]
	which belongs to the Nevanlinna class $\mc N_0$.
	Since $q(ri)=(a-b)\log r + i(a+b)\frac{\pi}{2}$, conditions (i)--(iii)
	in \cref{L91} are satisfied with $f(r)=\log r$,
	$\alpha=0$ and $\omega=i(b-a)$.
	Let $\mu$ be the measure in the representation $q=\ms C_0[\mu,p]$.
	\Cref{L16} only yields that
	\[
		\lim_{r\to\infty}\frac{\mu((0,r))}{r\log r} = 0, \qquad
		\lim_{r\to\infty}\frac{\mu((-r,0))}{r\log r} = 0.
	\]
	On the other hand, since $\Im q(ri)=(a+b)\frac{\pi}{2}$, we can apply \cref{L45}
	to obtain that $r\mapsto\mu((-r,r))$ is regularly varying and that
	\[
		\mu((-r,r)) \sim \frac{2}{\pi}r\Im q(ir) = (a+b)r.
	\]
	Note that, actually, $\mu((0,r))=br$ and $\mu((-r,0))=ar$.
\end{example}

\medskip\noindent
The next proposition shows that, in the case $\kappa=0$, the validity
of the first asymptotic relation in \eqref{L46} implies already (a) and (b)
in \cref{L45}.

\begin{proposition}\label{L50}
	Let $\mu$ be a measure on $\bb R$ such that $\int_{\bb R}(1+t^2)^{-1}\D\mu(t)<\infty$,
	let $c_0\in\bb R$ and set $q\DE c_0+\widetilde C[\mu]$.
	Assume that the limit
	\[
		\lim_{y\to\infty}\biggl(\raisebox{0.5ex}{$\displaystyle\Im q(iy)$}\bigg/
		\raisebox{-0.5ex}{$\displaystyle\frac{\mu((-y,y))}{y}$}\biggr)
	\]
	exists and is positive.  Then \textup{(a)} and \textup{(b)} in \cref{L45} are satisfied.
\end{proposition}

\begin{proof}
	Let $\mu_*$ be as in \eqref{L2}.  By \eqref{L175}, the following limit
	\[
		\lim_{t\to\infty}\frac{\mu_*([0,t))}{t\ms S[\mu_*](t)}
		= \lim_{y\to\infty}\frac{\mu_*([0,y^2))}{y^2\ms S[\mu_*](y^2)}
		= \lim_{y\to\infty}\frac{\mu((-y,y))}{y\Im q(iy)}
	\]
	exists and is positive.  Hence \cite[Theorem~B]{shea:1969} implies that
	$t\mapsto\mu_*([0,t))$ is regularly varying, and hence (a) and (b) hold.
\end{proof}

\subsection{The real part}

The real part of $q(iy)$ is more subtle since it can be written only in terms of 
a difference of two Stieltjes transforms and cancellations can arise.
We introduce the following notation for the main part of $\Re q(iy)$.
Let $\mu$ be a measure on $\bb R$ such that there exists $\ell\in\bb N_0$ with
\begin{equation}\label{L131}
	\int_{\bb R}\frac{\D\mu(t)}{(1+|t|)^{2\ell+1}}<\infty,
\end{equation}
and define 
\begin{equation}\label{L187}
	\RC_\ell[\mu](y) \DE (1-y^2)^\ell\int_{\bb R}\frac{t}{t^2+y^2}\cdot\frac{\D\mu(t)}{(1+t^2)^\ell},
	\qquad y>0.
\end{equation}

\begin{lemma}\label{L177}
	Let $\mu$, $\kappa$, $p$ and $q$ be as in \cref{L171},
	and let $\ms p(\mu)$ and $\RC_\ell$ be as in \eqref{L215} and \eqref{L187} respectively.
	\begin{Enumerate}
	\item
		If $\ms p(\mu)$ is odd, then there exists a real, even polynomial $\widetilde p$
		of degree at most $2\kappa-2$ such that
		\begin{equation}\label{L181}
			\Re q(iy) = (-1)^\kappa
			\biggl(c_{2\kappa}-\int_{\bb R}\frac{t}{(1+t^2)^{\kappa+1}}\DD\mu(t)\biggr)y^{2\kappa}
			+ \RC_\kappa[\mu](y) + \widetilde p(y).
		\end{equation}
	\item
		If $\ms p(\mu)$ is even, then there exists a real, even polynomial $\widetilde p$
		of degree at most $2\kappa$ such that
		\begin{equation}\label{L182}
			\Re q(iy) = \RC_{\kappa+1}[\mu](y) + \widetilde p(y)
		\end{equation}
		with $\widetilde p(y) = (-1)^\kappa c_{2\kappa}y^{2\kappa}+\BigO(y^{2\kappa-2})$.
	\end{Enumerate}
\end{lemma}

\begin{proof}
	\phantom{}

	(i)
	When $\ms p(\mu)$ is odd, then $\int_{\bb R}(1+|t|)^{-(2\kappa+1)}\D\mu(t)<\infty$.
	Hence from \eqref{L42} we obtain
	\[
		\Re q(iy) = \Re(p(iy))
		-(1-y^2)^\kappa\int_{\bb R}\frac{t}{(1+t^2)^{\kappa+1}}\DD\mu(t)
		+(1-y^2)^\kappa\int_{\bb R}\frac{t}{t^2+y^2}\cdot\frac{\D\mu(t)}{(1+t^2)^\kappa},
	\]
	which yields \eqref{L181}.
	
	(ii)
	When $\ms p(\mu)$ is even, then we use \eqref{L14} to write
	\[
		\Re q(iy) = \Re(p(iy))
		+ (1-y^2)^{\kappa+1}\int_{\bb R}\frac{t}{t^2+y^2}\cdot\frac{\D\mu(t)}{(1+t^2)^{\kappa+1}},
	\]
	which gives \eqref{L182}.
\end{proof}

\begin{assumption}\label{L176}
	Let $\mu$, $\kappa$, $p$ and $q$ be as in \cref{L171}.
	Assume that, if $\ms p(\mu)$ is odd, then
	\[
		c_{2\kappa} = \int_{\bb R}\frac{t}{(1+t^2)^{\kappa+1}}\DD\mu(t).
	\]
\end{assumption}

\begin{remark}\label{L205}
	Assume that $\ms p(\mu)=1$ and that \cref{L171,L173}
	are satisfied.  We then have $\kappa=\kappa(\mu)=0$ and, by \cref{L204}, we can write
	$q$ as in \eqref{L207}, i.e.\ $q(z)=c_0+\widetilde C[\mu](z)$.
	If, in addition, \cref{L176} is satisfied, then $c_0=\int_{\bb R}\frac{t}{1+t^2}\DD\mu(t)$
	and hence $q$ is the Cauchy transform of $\mu$, i.e.\
	\[
		q(z) = \int_{\bb R}\frac{1}{t-z}\DD\mu(t) = C[\mu](z),
	\]
	where $C[\mu]$ is defined in \eqref{L4}.
\end{remark}

\medskip\noindent
To investigate the behaviour of the real part, we often choose $\ell$ minimal so that \eqref{L187}
makes sense, i.e.\ let us set
\begin{equation}\label{L161}
	\ell(\mu) \DE \inf\biggl\{n\in\bb N_0\DS \int_{\bb R}\frac{\D\mu(t)}{(1+|t|)^{2n+1}}<\infty\biggr\}.
\end{equation}
It follows easily that $\ell(\mu)=\bigl\lfloor\frac{\ms p(\mu)}{2}\bigr\rfloor$,
where $\ms p(\mu)$ is as in \eqref{L215}.

\begin{remark}\label{L154}
	Let $\mu$ be a measure on $\bb R$ and define $\ms p(\mu)$, $\kappa(\mu)$ 
	and $\ell(\mu)$ as in \eqref{L215}, \eqref{L60} and \eqref{L161} respectively.
	If $\ms p(\mu)$ is odd, then $\kappa(\mu)=\ell(\mu)=\frac{\ms p(\mu)-1}{2}$.
	If $\ms p(\mu)$ is even, then $\kappa(\mu)+1=\ell(\mu)=\frac{\ms p(\mu)}{2}$.
	In both cases we have $\kappa(\mu)+\ell(\mu)=\ms p(\mu)-1$.
\end{remark}

\begin{lemma}\label{L178}
	Let $\mu$, $\kappa$, $p$ and $q$ be as in \cref{L171} and $\ell=\ell(\mu)$.
	Then there exists a real, even polynomial $\widetilde p$ of degree 
	at most $2\ell-2$ such that
	\[
		\Re q(iy) = \RC_\ell[\mu](y) + \widetilde p(y).
	\]
\end{lemma}

\begin{proof}
	If $\ms p(\mu)$ is odd, then $\ell(\mu)=\kappa(\mu)$ and the statement follows from \cref{L177}\,(i).
	If $\ms p(\mu)$ is even, then $\ell(\mu)=\kappa(\mu)+1$ and we can apply \cref{L177}\,(ii).
\end{proof}

\begin{definition}\label{L160}
	Let $\mu$ be a measure on $\bb R$, let $\mu_*$
	be the push-forward measure as in \eqref{L2}, and let $\ell\in\bb N_0$ 
	be such that \eqref{L131} is satisfied, i.e.\ $\ell\ge\ell(\mu)$.
	Define the measure $\sigma_\ell$ on $[0,\infty)$ by
	\begin{equation}\label{L162}
		\D\sigma_\ell(s) = \frac{\sqrt{s}}{(1+s)^\ell}\DD\mu_*(s), \qquad s\in[0,\infty),
	\end{equation}
	and set
	\[
		\mc F_\ell[\mu](y) \DE (1-y^2)^\ell\ms S[\sigma_\ell](y^2), \qquad y>0,
	\]
	where $\ms S$ is the Stieltjes transform defined in \eqref{L22}.
\end{definition}

\begin{lemma}\label{L163}
	Let $\mu$ be a measure on $\bb R$ and define the
	measures $\mu^+$ and $\mu^-$ on $\bb R$ by
	\begin{equation}\label{L164}
		\D\mu^+(t) \DE \mathds{1}_{[0,\infty)}(t)\DD\mu(t), \qquad
		\D\mu^-(t) \DE \mathds{1}_{(-\infty,0)}(t)\DD\mu(t), \qquad
		t\in\bb R.
	\end{equation}
	Let $\ell\in\bb N_0$ such that \eqref{L131} holds and
	let $\mc F_\ell$ be as in \cref{L160}.  Then
	\begin{align}
		\RC_\ell[\mu](y)
		&= \mc F_\ell[\mu^+](y) - \mc F_\ell[\mu^-](y),
		\label{L165}
		\\[1ex]
		\big|\RC_\ell[\mu](y)\big|
		&\le \big|\mc F_\ell[\mu](y)\big|
		\label{L166}
	\end{align}
	for $y>0$.
\end{lemma}

\begin{proof}
	Let $\mu_*$ and $\mu_*^\pm$ be the push-forward measures of $\mu$ and $\mu^\pm$
	respectively as in \eqref{L2},
	and let $\sigma_\ell$, $\sigma_\ell^\pm$ be as in \eqref{L162}.

	The definition of $\ell$ in \eqref{L161} implies that
	\[
		\int_{[0,\infty)}\frac{\D\sigma_\ell^\pm(s)}{1+s}
		= \int_{[0,\infty)}\frac{\sqrt{s}}{(1+s)^{\ell+1}}\DD\mu_*^\pm(s)
		= \int_{\bb R}\frac{|t|}{(1+t^2)^{\ell+1}}\DD\mu^\pm(t) < \infty,
	\]
	and hence $\ms S[\sigma_\ell^+]$ and $\ms S[\sigma_\ell^-]$ are well defined.
	Moreover,
	\begin{align*}
		\RC_\ell[\mu](y) &= (1-y^2)^\ell\biggl[\int_{\bb R}\frac{1}{t^2+y^2}\cdot\frac{|t|}{(1+t^2)^\ell}\DD\mu^+(t)
		-\int_{\bb R}\frac{1}{t^2+y^2}\cdot\frac{|t|}{(1+t^2)^\ell}\DD\mu^-(t)\biggr]
		\\[1ex]
		&= (1-y^2)^\ell
		\biggl[\int_{[0,\infty)}\frac{1}{s+y^2}\cdot\frac{\sqrt{s}}{(1+s)^\ell}\DD\mu_*^+(s)
		-\int_{[0,\infty)}\frac{1}{s+y^2}\cdot\frac{\sqrt{s}}{(1+s)^\ell}\DD\mu_*^-(s)\biggr]
		\\[1ex]
		&= (1-y^2)^\ell\Bigl[\ms S[\sigma_\ell^+](y^2)-\ms S[\sigma_\ell^-](y^2)\Bigr],
	\end{align*}
	which yields \eqref{L165}.
	In a similar way as above one shows that $\ms S[\sigma_\ell]$ is well defined.
	Further, we have
	\begin{align*}
		\big|\RC_\ell[\mu](y)\big| 
		&\le \big|(1-y^2)^\ell\big|\int_{\bb R}\frac{1}{t^2+y^2}\cdot\frac{|t|}{(1+t^2)^\ell}\D\mu(t)
		\\[1ex]
		&= \big|(1-y^2)^\ell\big|\int_{[0,\infty)}\frac{1}{s+y^2}\cdot\frac{\sqrt{s}}{(1+s)^\ell}\DD\mu_*(s)
		= \big|\ms S[\sigma_\ell](y^2)\big|,
	\end{align*}
	which proves \eqref{L166}.
\end{proof}

\medskip\noindent
In the following key proposition the asymptotic behaviour of $\mc F_\ell[\mu]$ is determined.
It is used in \cref{L56} and in the Abelian implications in \cref{L28,L40}.

\begin{proposition}\label{L167}
	Let $\mu$ be a measure on $\bb R$, set $\ell\DE\ell(\mu)$
	and define $\mc F_\ell[\mu]$ as in \cref{L160}.
	Further, assume that $t\mapsto\mu((-t,t))$ is regularly varying with index $\beta$.
	Then $\beta\in[2\ell-1,2\ell+1]$. 
	% $\ell\in\bigl[\frac{\beta-1}{2},\frac{\beta+1}{2}\bigr]$.
	Moreover, if $\beta=0$, then
	\begin{equation}\label{L168}
		\mc F_\ell[\mu](y) \ll \frac{\mu((-y,y))}{y},
		\qquad y\to\infty.
	\end{equation}
	If $\beta>0$, then $(-1)^\ell\mc F_\ell[\mu]$ is regularly varying 
	with index $\beta-1$ and satisfies
	\begin{equation}\label{L169}
		\mc F_\ell[\mu](y) \sim
		\begin{cases}
			\displaystyle \frac{\frac{\pi\beta}{2}}{\cos\frac{\pi\beta}{2}}\cdot\frac{\mu((-y,y))}{y},
			& \beta\notin2\bb N-1,
			\\[3ex]
			\displaystyle (-1)^\ell\beta y^{\beta-1}\int_1^y\frac{\mu((-t,t))}{t^{\beta+1}}\DD t,
			& \beta=2\ell-1,
			\\[3ex]
			\displaystyle (-1)^\ell\beta y^{\beta-1}\int_y^\infty\frac{\mu((-t,t))}{t^{\beta+1}}\DD t,
			& \beta=2\ell+1,
		\end{cases}
	\end{equation}
	and
	\begin{equation}\label{L156}
		y^{2\ell-2} \ll \big|\mc F_\ell[\mu](y)\big| \ll y^{2\ell}
	\end{equation}
	as $y\to\infty$.
	In particular, if $\beta\in2\bb N-1$, then
	\begin{equation}\label{L170}
		\big|\mc F_\ell[\mu](y)\big| \gg \frac{\mu((-y,y))}{y}.
	\end{equation}
\end{proposition}

\begin{remark}\label{L158}
\rule{0ex}{1ex}
\begin{Enumerate}
\item
	In the situation of \cref{L167} assume that $\beta\in2\bb N-1$ and
	that $t\mapsto\mu((-t,t))$ is regularly varying with index $\beta$.
	It follows from the definition of $\ell(\mu)$ in \eqref{L161} and
	from \cref{L159} that, with $\ell=\ell(\mu)$,
	\begin{alignat*}{5}
		\beta &= 2\ell-1 
		\quad &&\Leftrightarrow\quad & 
		\int_{\bb R}\frac{\D\mu(t)}{(1+|t|)^\beta} &= \infty
		\quad &&\Leftrightarrow\quad & 
		\int_1^\infty\frac{\mu((-t,t))}{t^{\beta+1}}\DD t &= \infty,
		\\[2ex]
		\beta &= 2\ell+1 
		\quad &&\Leftrightarrow\quad & 
		\int_{\bb R}\frac{\D\mu(t)}{(1+|t|)^\beta} &< \infty
		\quad &&\Leftrightarrow\quad &
		\int_1^\infty\frac{\mu((-t,t))}{t^{\beta+1}}\DD t &< \infty.
	\end{alignat*}
\item
	Instead of the assumption $\ell=\ell(\mu)$ in the proposition, 
	let us consider the case when $\ell>\ell(\mu)$.  Then
	\[
		\sigma_\ell([0,\infty)) = \int_{[0,\infty)}\frac{\sqrt{s}}{(1+s)^\ell}\DD s
		= \int_{\bb R}\frac{|t|}{(1+t^2)^\ell}\DD\mu(t) < \infty,
	\]
	and, by \cref{L26}, we have
	\[
		\mc F_\ell[\mu](y) \sim (-1)^\ell y^2\ms S[\sigma_\ell](y^2)
		\sim (-1)^\ell\sigma_\ell([0,\infty))y^{2\ell-2}
	\]
	as $y\to\infty$.
\end{Enumerate}
\end{remark}

\begin{proof}[Proof of \cref{L167}]
	Let $\mu_*$ be the push-forward measure of $\mu$ as \eqref{L2},
	and let $\sigma_\ell$ be as in \eqref{L162}.
	Then $s\mapsto\mu_*([0,s))=\mu((-\sqrt{s},\sqrt{s}))$ is regularly varying
	with index $\frac\beta2$.
	Set $h(s) \DE \frac{\sqrt{s}}{(1+s)^\ell} \sim s^{\frac12-\ell}$ as $s\to\infty$.
	Then
	\begin{equation}\label{L188}
		\int_{[0,\infty)}h(s)\DD\mu_*(s) = \int_{\bb R}\frac{|t|}{(1+t^2)^\ell}\DD\mu(t).
	\end{equation}

	Let us first consider the case when the integrals in \eqref{L188} are finite. 
	Then $\ell=0$ by the definition of $\ell$, and $\sigma_0$ and $\mu$ are finite measures, 
	which implies that $\beta=0$.  From \cref{L26} we obtain that
	\[
		\mc F_0[\mu](y) = \ms S[\sigma_0](y^2) \sim \frac{\sigma_0([0,\infty))}{y^2}
		\ll \frac{\mu(\bb R)}{y} \sim \frac{\mu((-y,y))}{y},
	\]
	which proves \eqref{L168} in this case.
	
	For the rest of the proof assume that the integrals in \eqref{L188} are infinite.
	The definition of $\ell$ also implies that
	\[
		\int_{[0,\infty)}\frac{h(s)}{1+s}\DD\mu_*(s) = \int_{\bb R}\frac{|t|}{(1+t^2)^{\ell+1}}\DD\mu(t) 
		< \infty,
	\]
	which shows that \eqref{L118} with $\nu=\mu_*$ is satisfied.
	Hence we can apply \cref{L122} with $\alpha=\frac{\beta}{2}$ and $\gamma=\frac12-\ell$,
	which yields $\alpha+\gamma=\frac\beta2+\frac12-\ell\in[0,1]$, i.e.\
	$\beta\in[2\ell-1,2\ell+1]$.
	\Cref{L122} also implies that, if $\beta>0$, then
	$(-1)^\ell\mc F_\ell[\mu](y) \sim y^{2\ell}\ms S[\mu_\ell](y^2)$
	is regularly varying with index $2\ell+2\bigl(\frac{\beta}{2}+\frac12-\ell-1\bigr)=\beta-1$.

	Let us first consider the case when $\beta\notin2\bb N-1$ and $\beta>0$.  
	Then $\alpha+\gamma=\frac{\beta}{2}+\frac12-\ell\in(0,1)$, and from \cref{L122} we obtain
	\begin{align*}
		\mc F_\ell[\mu](y) &\sim (-1)^\ell y^{2\ell}\ms S[\sigma_\ell](y^2)
		\sim (-1)^\ell y^{2\ell}\frac{\pi\frac{\beta}{2}}{\sin\bigl(\pi\bigl(\frac{\beta}{2}+\frac12-\ell\bigr)\bigr)}
		\bigl(y^2\bigr)^{-\frac12-\ell}\mu_*([0,y^2))
		\\[1ex]
		&= \frac{\frac{\pi\beta}{2}}{\sin\bigl(\frac{\pi\beta}{2}+\frac{\pi}{2}\bigr)}
		\cdot\frac{\mu((-y,y))}{y} 
		= \frac{\frac{\pi\beta}{2}}{\cos\frac{\pi\beta}{2}}\cdot\frac{\mu((-y,y))}{y}
	\end{align*}
	as $y\to\infty$.  

	Next assume that $\beta=2\ell-1$.
	Then $\alpha+\gamma=\frac{\beta}{2}+\frac12-\ell=0$, and \cref{L122} and the substitution $s=t^2$ yield
	\begin{align*}
		\mc F_\ell[\mu](y) &\sim (-1)^\ell y^{2\ell}\ms S[\sigma_\ell](y^2)
		\sim (-1)^\ell y^{2\ell}\frac{\beta}{2}\cdot\frac{1}{y^2}
		\int_1^{y^2} s^{-\frac12-\ell}\mu_*([0,s))\DD s
		\\[1ex]
		&= (-1)^\ell\frac{\beta}{2} y^{2\ell-2}\int_1^y t^{-1-2\ell}\mu_*([0,t^2))\,2t\DD t
		= (-1)^\ell\beta y^{2\ell-2}\int_1^y \frac{\mu((-t,t))}{t^{2\ell}}\DD t,
	\end{align*}
	which shows \eqref{L169} in this case.

	The proofs of \eqref{L168} when $\beta=0$ and the integrals in \eqref{L188} are infinite 
	and of \eqref{L169} in the remaining case $\beta=2\ell+1$ are similar.

	We can use \eqref{L189} to obtain
	\[
		\big|\mc F_\ell[\mu](y)\big| \sim y^{2\ell}\ms S[\sigma_\ell](y^2)
		\qquad\text{and}\qquad
		\frac{1}{y^2} \ll \ms S[\sigma_\ell](y^2) \ll 1,
	\]
	which yields \eqref{L156}.

	Finally, the relation in \eqref{L170} follows easily 
	from \eqref{L169} and \cref{L36}.
\end{proof}

\medskip\noindent
We also need the following comparison result.

\begin{lemma}\label{L155}
	Let $\mu_1$ and $\mu_2$ be measures on $\bb R$ and set $\ell\DE\ell(\mu_2)$.
	Further, assume that $t\mapsto\mu_2((-t,t))$ is regularly varying 
	with index $\beta>0$ and that the limit
	\[
		\lim_{t\to\infty}\frac{\mu_1((-t,t))}{\mu_2((-t,t))}
	\]
	exists in $[0,\infty)$.  Then $\mc F_\ell[\mu_1]$ is well defined, 
	i.e.\ $\ell\ge\ell(\mu_1)$, and
	\[
		\lim_{y\to\infty}\frac{\mc F_\ell[\mu_1](y)}{\mc F_\ell[\mu_2](y)}
		= \lim_{t\to\infty}\frac{\mu_1((-t,t))}{\mu_2((-t,t))}.	
	\]
\end{lemma}

\begin{proof}
	Since $\beta>0$, it follows as in the proof of \cref{L167}
	that \eqref{L192} with $\nu_2=(\mu_2)_*$ and $h(s)=\frac{\sqrt{s}}{(1+s)^\ell}$
	is satisfied.  Now the claim follows from \cref{L121}.
\end{proof}

\medskip\noindent
When $t\mapsto\mu((-t,t))$ is regularly varying with an index that is not an odd integer,
the real part of $q(iy)$ is dominated by the imaginary part, as the following proposition shows.

\begin{proposition}\label{L56}
	Suppose that $\mu$, $p$ and $q$ satisfy \cref{L171,L173,L176}.
	Further, assume that $t\mapsto\mu((-t,t))$ is regularly varying with index $\beta$
	such that $\beta\notin2\bb N-1$.  Then
	\begin{equation}\label{L180}
		\limsup_{y\to\infty}\bigg|\frac{\Re q(iy)}{\Im q(iy)}\bigg|
		\le \Big|\tan\frac{\pi\beta}{2}\Big|.
	\end{equation}
\end{proposition}

\begin{proof}
	Set $\ell\DE\ell(\mu)$.  It follows from \cref{L178,L163} that
	\begin{equation}\label{L179}
		\bigg|\frac{\Re q(iy)}{\Im q(iy)}\bigg|
		= \bigg|\frac{\RC_\ell[\mu](y)+\widetilde p(y)}{\Im q(iy)}\bigg|
		\le \frac{|\mc F_\ell[\mu](y)|}{|\Im q(iy)|} 
		+ \BigO\biggl(\frac{y^{2\ell-2}}{|\Im q(iy)|}\biggr).
	\end{equation}
	\Cref{L45} implies that $y\mapsto|\Im q(iy)|$ is regularly varying with 
	index $\beta-1$.  Further, we obtain from \cref{L167} that $\beta\in[2\ell-1,2\ell+1]$.
	Since, by assumption, $\beta\notin2\bb N-1$, we have $\beta-1>2\ell-2$
	and hence $|\Im q(iy)|\gg y^{2\ell-2}$, i.e.\ the $\BigO$-term on the right-hand side
	of \eqref{L179} converges to $0$ as $y\to\infty$.
	
	From \cref{L167,L45} we obtain
	\begin{align*}
		|\mc F_\ell[\mu](y)|\ &
		\begin{cases}
			\displaystyle \ll \frac{\mu((-y,y))}{y},
			& \beta=0,
			\\[3ex]
			\displaystyle \sim \frac{\frac{\pi\beta}{2}}{|\cos\frac{\pi\beta}{2}|}
			\cdot\frac{\mu((-y,y))}{y},
			& \beta\in(0,\infty)\setminus(2\bb N-1),
		\end{cases}
		\\[1ex]
		|\Im q(iy)|\ &
		\begin{cases}
			\displaystyle \sim \frac{\frac{\pi\beta}{2}}{|\sin\frac{\pi\beta}{2}|}
			\cdot\frac{\mu((-y,y))}{y},
			& \beta\in[0,\infty)\setminus2\bb N,
			\\[3ex]
			\displaystyle \gg \frac{\mu((-y,y))}{y},
			& \beta\in2\bb N,
		\end{cases}
	\end{align*}
	which yields
	\[
		\lim_{y\to\infty}\frac{|\mc F_\ell[\mu](y)|}{|\Im q(iy)|} 
		= \Big|\tan\frac{\pi\beta}{2}\Big|
	\]
	and hence \eqref{L180}.
\end{proof}

\section{Abelian--Tauberian theorems}
\label{L114}

In the following theorem, one of the main results of this paper,
we combine the Tauberian and the Abelian theorems from the previous sections.
In most cases we can give a full characterisation when the asymptotic behaviour
of the regularised Cauchy transform is described by a regularly varying function.
Thereby, we assume that the behaviour of $q=\ms C_\kappa[\mu,p]$ is not governed 
by the polynomial summand. 

\subsection{The generic situation}

Recall the notation $\ms p(\mu)$ from \cref{L215} and $\kappa(\mu)$ from \cref{L60}.

\begin{theorem}\label{L28}
	Let $\mu$ be a measure on $\bb R$ such that $\ms p(\mu)<\infty$
	and set $\kappa\DE\kappa(\mu)$.
	Further, let $p\in\bb R[z]$ with $p(z)=c_{2\kappa+1}z^{2\kappa+1}+\ldots+c_0$
	such that $(\mu,p)\in\bb E_{\le\kappa}$, and assume that 
	\[
		c_{2\kappa+1} = \int_{\bb R}\frac{\D\mu(t)}{(1+t^2)^{\kappa+1}}.
	\]
	and, if $\ms p(\mu)$ is odd, that 
	\[
		c_{2\kappa} = \int_{\bb R}\frac{t}{(1+t^2)^{\kappa+1}}\DD\mu(t).
	\]
	Set $q\DE\ms C_\kappa[\mu,p]$.

	Let $\beta\in[0,\infty)$ and consider the following statements:
	\begin{Enumeratealph}
	\item
		the symmetrised distribution function $t\mapsto\mu((-t,t))$ is
		regularly varying with index~$\beta$;
	\item[\textup{(a)$'$}\hspace*{-0.6ex}]
		the limit
		\begin{equation}\label{L80}
			\zeta\DE\lim_{t\to\infty}\frac{\mu((-t,0))}{\mu([0,t))}
		\end{equation}
		exists in $[0,\infty]$;
	\item
		there exist a regularly varying function $f\DF[x_0,\infty)\to(0,\infty)$ with $x_0>0$ 
		and a constant $\omega\in\bb C\setminus\{0\}$ such that
		\begin{equation}\label{L32}
			\lim_{r\to\infty}\frac{q(rz)}{f(r)}
			= i\omega\Bigl(\frac{z}{i}\Bigr)^{\beta-1}
		\end{equation}
		holds locally uniformly for $z\in\bb C^+$.
	\end{Enumeratealph}
	Then the following relations hold.
	\begin{Enumerate}
	\item
		If $\beta\notin\bb N_0$, then
		\[
			\textup{(a)} \;\wedge\; \textup{(a)$'$} \quad\Leftrightarrow\quad \textup{(b)}.
		\]
	\item
		Assume that $\beta\in2\bb N_0$.  Then
		\[
			\textup{(a)} \;\;\Rightarrow\;\; \textup{(b)}.
		\]
		If $\beta=0$ or $\int_{\bb R}(1+|t|)^{-\beta}\D\mu(t)<\infty$, then
		\[
			\textup{(b)} \;\;\Rightarrow\;\; \textup{(a)}.
		\]
	\item
		If $\beta\in2\bb N-1$, then
		\[
			\textup{(a)} \;\wedge\; \textup{(a)$'$} \;\wedge\; \zeta\ne1 
			\quad\Rightarrow\quad \textup{(b)}.
		\]
%		and
%		\[
%			\textup{(b)} \;\wedge\; ???
%			\quad\Rightarrow\quad \textup{(a)} \;\wedge\; \textup{(a)$'$}.
%		\]
	\end{Enumerate}
	Further, assume that either \textup{(a)} and \textup{(b)} hold and $\beta\in2\bb N_0$,
	or that \textup{(a)}, \textup{(a)$'$} and \textup{(b)} hold and $\beta\notin\bb N_0$,
	or that \textup{(a)}, \textup{(a)$'$} with $\zeta\ne1$ and \textup{(b)} hold and $\beta\in2\bb N-1$;
	then $\omega$ and $f$ can be chosen as
	\begin{equation}\label{L190}
		\omega = (-1)^{\ms p(\mu)-1}\biggl(\cos\frac{\pi\beta}{2}
		+i\frac{\zeta-1}{\zeta+1}\cdot\sin\frac{\pi\beta}{2}\biggr),
	\end{equation}
	where $\frac{\zeta-1}{\zeta+1}$ is understood as $1$ when $\zeta=\infty$,
	where $\ms p(\mu)$ is defined in \eqref{L215}, and
	\begin{equation}\label{L191}
		f(r) =
		\begin{cases}
			\dfrac{\pi\beta}{|\sin(\pi\beta)|}\cdot\dfrac{\mu((-r,r))}{r},
			& \beta\in[0,\infty)\setminus\bb N,
			\\[3ex]
			\displaystyle \beta r^{\beta-1}\int_1^r\frac{\mu((-t,t))}{t^{\beta+1}}\DD t,
			& \displaystyle \beta\in\bb N \;\wedge\; 
			\int_{\bb R}\frac{\D\mu(t)}{(1+|t|)^\beta}=\infty,
			\\[3ex]
			\displaystyle \beta r^{\beta-1}\int_r^\infty\frac{\mu((-t,t))}{t^{\beta+1}}\DD t, \quad
			& \displaystyle \beta\in\bb N \;\wedge\; 
			\int_{\bb R}\frac{\D\mu(t)}{(1+|t|)^\beta}<\infty;
		\end{cases}
	\end{equation}
	here $\frac{\pi\beta}{|\sin(\pi\beta)|}$ is understood as $1$ when $\beta=0$.
\end{theorem}

\medskip\noindent
Before we prove the theorem, let us add a couple of comments.

\begin{remark}\label{L157}
\rule{0ex}{1ex}
\begin{Enumerate}
\item
	Note that, in the case $\beta\in\bb N$, one has $f(r)\gg\frac{\mu((-r,r))}{r}$ by \cref{L36}.
\item
	When $\beta\in2\bb N_0$, then $\omega\in\bb R$ and hence $|\Im q(iy)|\gg|\Re q(iy)|$.
	When $\beta\in2\bb N-1$ in the theorem, then $\omega\in i\bb R$ and 
	hence $|\Re q(iy)|\gg|\Im q(iy)|$.
\item
	It follows from the proof of \cref{L28} (see \eqref{L220}) that, if $\beta\in2\bb N-1$ and (a) and (a)$'$ 
	are satisfied with $\zeta=1$, then $|q(iy)|\ll f(y)$ with $f$ as in \eqref{L191}.
	A more detailed discussion of some cases in this situation are contained
	in \cref{L40} below.
\item
	In the case when $\beta\notin\bb N_0$, condition (a) is not sufficient to guarantee (b)
	since one can easily construct measures $\mu$ that satisfy (a) but not (a)$'$, e.g.\
	by distributing the mass in an alternating way on the positive and negative half-axes.
\item
	In \cref{L212} below we show that the converse implication in \cref{L28}\,(iii) does not hold.
	However, see \cref{L40}\,(ii) below for a Tauberian implication in 
	the case when $\omega\notin i\bb R$.
\item
	\Cref{L53} below deals with the situation when $\beta=1$ and $\zeta=1$,
	where \cref{L28} is not applicable.  This example shows that, in general, 
	the asymptotic behaviour of $q(iy)$ is not determined by the leading asymptotic 
	behaviour of $\mu([0,t))$ and $\mu((-t,0))$.
	However, in this example statement (b) in \cref{L28} still holds.
	We do not know whether there exist a measure $\mu$ and a function $q$ such that
	$\beta\in2\bb N-1$ and (a) and (a)$'$ hold with $\zeta=1$ but (b) does not hold.
\end{Enumerate}
\end{remark}

\begin{proof}[Proof of \cref{L28}]
Let $\kappa=\kappa(\mu)$ and $\ell=\ell(\mu)$.
\begin{Steps}
\item 
	It follows from \cref{L91} that \textup{(b)} is equivalent to
	\begin{Enumeratealph}
	\item[\textup(c)]
		there exist a regularly varying function $f\DF[x_0,\infty)\to(0,\infty)$ with $x_0>0$ 
		and a constant $\omega\in\bb C\setminus\{0\}$ such that
		\begin{equation}\label{L195}
			\lim_{y\to\infty}\frac{q(iy)}{f(y)} = i\omega.
		\end{equation}
	\end{Enumeratealph}
\item
	Let us first consider the case when $\beta\in2\bb N_0$.
	It follows from \cref{L56} that $|\Re q(iy)|\ll|\Im q(iy)|$ and hence, \eqref{L195}
	is equivalent to
	\[
		\lim_{y\to\infty}\frac{\Im q(iy)}{f(y)} = \omega.
	\]
	From \cref{L45,L48} we therefore obtain the implications in (ii).
	Assume now that \textup{(a)} and \textup{(b)} hold.
	It follows from \cref{L167} that $\beta\in[2\ell-1,2\ell+1]$ and 
	hence $\ell=\frac{\beta}{2}$.
	With $f$ as in \eqref{L191} we obtain from \eqref{L46} and \cref{L48} that
	\begin{equation}\label{L196}
		\lim_{y\to\infty}\frac{\Im q(iy)}{f(y)}
		= (-1)^\kappa.
	\end{equation}
	On the other hand, by \cref{L154}, $\omega$ from \eqref{L190} equals
	$\omega = (-1)^{\ms p(\mu)-1}\cos\frac{\pi\beta}{2} = (-1)^{\kappa+\ell}(-1)^\ell = (-1)^\kappa$,
	which coincides with the limit in \eqref{L196}.
\item
	Next we assume that $\beta\notin\bb N_0$ and (b) holds.
	It follows from \cref{L91} that the index of the regularly varying function $f$
	is $\alpha\DE\beta-1$, and hence $r\mapsto rf(r)$ is regularly varying
	with index $\beta$.
	Now \cref{L16} implies that \eqref{L92} and \eqref{L93} hold.
	Since $\alpha\notin\bb N_0$, at least one of the right-hand sides of
	\eqref{L92}, \eqref{L93} is non-zero.
	Taking the quotient of \eqref{L93} and \eqref{L92} we obtain that the
	limit in \eqref{L80} exists; further, $t\mapsto\mu((-t,t))$ is regularly varying
	with index $\beta$, i.e.\ (a) and (a)$'$ hold.
\item
	For the rest of the proof we assume that $\beta\notin2\bb N_0$,
	that (a) and (a)$'$ hold, and if $\beta\in2\bb N-1$, then also $\zeta\ne1$.  
	Let $\mu^+$ and $\mu^-$ be as in \eqref{L164} and $f$ as in \eqref{L191}.
	From \cref{L178,L163} we obtain that
	\begin{align}
		\frac{1}{i}\cdot\frac{q(iy)}{f(y)}
		&= \frac{\mc F_\ell[\mu](y)}{f(y)}\cdot\frac{\Im q(iy)-i\Re q(iy)}{\mc F_\ell[\mu](y)}
		\nonumber\\[1ex]
		&= \frac{\mc F_\ell[\mu](y)}{f(y)}\biggl[\frac{\Im q(iy)}{\mc F_\ell[\mu](y)}
		-i\biggl(\frac{\mc F_\ell[\mu^+](y)}{\mc F_\ell[\mu](y)}
		- \frac{\mc F_\ell[\mu^+](y)}{\mc F_\ell[\mu](y)} 
		+ \frac{\widetilde p(y)}{\mc F_\ell[\mu](y)}\biggr)\biggr],
		\label{L199}
	\end{align}
	where $\widetilde p$ is a real, even polynomial of degree at most $2\ell-2$.
	In the next couple of steps we evaluate the limits of parts of this expression.
\item
	First we show that
	\begin{equation}\label{L201}
		\lim_{y\to\infty}\frac{\mc F_\ell[\mu](y)}{f(y)} 
		= (-1)^{\ms p(\mu)-1}\sin\frac{\pi\beta}{2}.
	\end{equation}
	We start with the case when $\beta\notin2\bb N-1$.
	It follows from \cref{L159} that $\beta\in(\ms p(\mu)-1,\ms p(\mu))$ and 
	hence $\sgn(\sin(\pi\beta))=(-1)^{\ms p(\mu)-1}$.  
	Further, we obtain from \cref{L167} and \eqref{L191} that
	\[
		\lim_{y\to\infty}\frac{\mc F_\ell[\mu](y)}{f(y)} 
		= \frac{|\sin(\pi\beta)|}{2\cos\frac{\pi\beta}{2}}
		= \frac{(-1)^{\ms p(\mu)-1}\sin(\pi\beta)}{2\cos\frac{\pi\beta}{2}}
		= (-1)^{\ms p(\mu)-1}\sin\frac{\pi\beta}{2}.
	\]
	Now assume that $\beta\in2\bb N-1$.  
	Then $\lim_{y\to\infty}\frac{\mc F_\ell[\mu](y)}{f(y)}=(-1)^\ell$.
	On the other hand, by \cref{L45} we have $\beta\in[2\kappa,2\kappa+2]$
	and hence $\beta=2\kappa+1$, which, by \cref{L154}, implies that
	\[
		(-1)^{\ms p(\mu)-1}\sin\frac{\pi\beta}{2}
		= (-1)^{\kappa+\ell}\sin\Bigl(\pi\kappa+\frac{\pi}{2}\Bigr)
		= (-1)^{\kappa+\ell}\cos(\pi\kappa) = (-1)^\ell,
	\]
	which proves \eqref{L201} also in this case.
\item
	Next we show that
	\begin{equation}\label{L200}
		\lim_{y\to\infty}\frac{\Im q(iy)}{\mc F_\ell[\mu](y)} = \cot\frac{\pi\beta}{2}.
	\end{equation}
	When $\beta\notin2\bb N-1$, we can use \cref{L45,L167} 
	to obtain \eqref{L200}.
	When $\beta\in2\bb N-1$, \cref{L45} and \eqref{L170} imply
	that $\lim_{y\to\infty}\frac{\Im q(iy)}{\mc F_\ell[\mu](y)}=0$,
	and again \eqref{L200} holds.

\item
	Let us consider the expressions within the round brackets 
	on the right-hand side of \eqref{L199}.
	Assume first that $\mu^+$ is not the zero measure.
	The relation $\mu=\mu^++\mu^-$ implies that
	\[
		\frac{\mu((-t,t))}{\mu^+((-t,t))}
		= 1+\frac{\mu^-((-t,t))}{\mu^+((-t,t))}
		= 1+\frac{\mu((-t,0))}{\mu([0,t))}
		\to 1+\zeta,
		\qquad t\to\infty,
	\]
	and hence
	\[
		\lim_{t\to\infty}\frac{\mu^+((-t,t))}{\mu((-t,t))} = \frac{1}{1+\zeta},
		\qquad
		\lim_{t\to\infty}\frac{\mu^-((-t,t))}{\mu((-t,t))} = \frac{\zeta}{1+\zeta}.
	\]
	If $\mu^+$ is the zero measure, then these relations also hold with $\zeta=\infty$.
	Together with \cref{L155} this yields
	\begin{equation}\label{L197}
		\lim_{y\to\infty}\frac{\mc F_\ell[\mu^+](y)}{\mc F_\ell[\mu](y)} = \frac{1}{1+\zeta},
		\qquad
		\lim_{y\to\infty}\frac{\mc F_\ell[\mu^-](y)}{\mc F_\ell[\mu](y)} = \frac{\zeta}{1+\zeta}.
	\end{equation}
	Further, \eqref{L156} implies that
	\begin{equation}\label{L198}
		\lim_{y\to\infty}\frac{\widetilde p(y)}{\mc F_\ell[\mu](y)} = 0.
	\end{equation}
\item
	Combining \eqref{L199}, \eqref{L201}, \eqref{L200}, \eqref{L197} and \eqref{L198}
	we arrive at
	\begin{equation}\label{L220}
		\frac{1}{i}\lim_{y\to\infty}\frac{q(iy)}{f(y)}
		= (-1)^{\ms p(\mu)-1}\sin\frac{\pi\beta}{2}
		\biggl(\cot\frac{\pi\beta}{2}-i\frac{1-\zeta}{1+\zeta}\biggr),
	\end{equation}
	which is equal to the expression on the right-hand side of \eqref{L190};
	this proves \eqref{L195}.
	Since, by assumption, either $\beta\notin2\bb N-1$ or $\zeta\ne1$,
	we have $\omega\ne0$.
	It follows from \eqref{L201} that 
	\[
		f(y)\sim\big|\sin\tfrac{\pi\beta}{2}\big|\cdot\big|\mc F_\ell[\mu](y)\big|.
	\]
	By \cref{L167} the right-hand side, and hence also $f$, is regularly varying.
	This shows that (c) is satisfied, which, in turn, implies (b).
\end{Steps}
\end{proof}

\medskip\noindent
The following example shows that the converse of the implication in \cref{L28}\,(iii)
does not hold.

\begin{example}\label{L212}
	Let $\ell\in\bb N$.  With $h(s)=\frac{\sqrt{s}}{(1+s)^\ell}$, $s\in[0,\infty)$,
	choose $\sigma$ and $\nu$ as in \cref{L202}.
	Moreover, let $\mu$ be the measure on $\bb R$ such that $\mu((-\infty,0))=0$
	and $\mu_*=\mu^+_*=\nu$.  Then
	\[
		\mu\bigl(\bigl\{e^{\frac{k}{2}}\bigr\}\bigr) = \mu_*^+(\{e^k\}) 
		= \frac{(1+e^k)^\ell}{e^{\frac{k}{2}}},
		\qquad k\in\bb N.
	\]
	From the relation
	\[
		\int_{\bb R}\frac{\D\mu(t)}{(1+|t|)^n}
		= \sum_{k=1}^\infty \frac{1}{\bigl(1+e^{\frac{k}{2}}\bigr)^n}\cdot\frac{(1+e^k)^\ell}{e^{\frac{k}{2}}}
	\]
	we can easily deduce that $\ms p(\mu)=2\ell$ and hence $\ell=\ell(\mu)$.
	Now set $q=\ms C_\kappa[\mu,p]$ with $\kappa=\kappa(\mu)$ and
	a polynomial $p$ such that \eqref{L213} holds.
	Then \cref{L171,L173,L176} are satisfied
	since $\ms p(\mu)$ is even.
	It follows from \cref{L202} that $t\mapsto\mu((-t,t))=\mu^+([0,t))=\nu([0,t^2))$
	is not regularly varying, i.e.\ (a) in \cref{L28} is not satisfied.

	Let us now consider the behaviour of $q(iy)$ as $y\to\infty$.
	It is clear that $\sigma_\ell=\sigma$ with $\sigma_\ell$ as in \eqref{L162}.
	From \cref{L163} and \eqref{L209} we obtain
	\[
		\RC_\ell[\mu](y) = \mc F_\ell[\mu^+](y)
		\sim (-1)^\ell y^{2\ell}\ms S[\sigma](y^2)
		\sim 2(-1)^\ell y^{2\ell-2}\log y,
	\]
	and hence
	\begin{equation}\label{L214}
		\Re q(iy) \sim 2(-1)^\ell y^{2\ell-2}\log y, \qquad y\to\infty,
	\end{equation}
	by \cref{L178}.
	Next we show that the imaginary part is dominated by the real part.
	Let $\tau_\kappa$ be as in \eqref{L218}.
	Since $\kappa(\mu)=\ell(\mu)-1$ by \cref{L154}, we have
	\[
		\tau_\kappa(\{e^k\}) 
		= \frac{1}{(1+e^k)^\kappa}\cdot\frac{(1+e^k)^\ell}{e^{\frac{k}{2}}} 
		= \frac{1+e^k}{e^{\frac{k}{2}}},
		\qquad k\in\bb N.
	\]
	Let $\tilde\tau$ be the measure on $[0,\infty)$ 
	with $\tilde\tau([0,s))=\frac{2\sqrt{e}}{\sqrt{e}-1}\bigl(s^{\frac12}-1\bigr)$
	for $s\in(1,\infty)$ and $\tilde\tau([0,1])=0$.
	For $t\in(e^n,e^{n+1}]$, $n\in\bb N$, we have
	\[
		\tau_\kappa([0,t)) = \sum_{k=1}^n\tau_k(\{e^k\})
		\le 2\sum_{k=1}^n e^{\frac{k}{2}}
		= \frac{2\sqrt{e}}{\sqrt{e}-1}\bigl(e^{\frac{n}{2}}-1\bigr)
		\le \frac{2\sqrt{e}}{\sqrt{e}-1}\bigl(t^{\frac12}-1\bigr)
		= \tilde\tau([0,t)).
	\]
	By \cref{L54}\,(i) and \cref{L3} we have
	\begin{align}
		|\Im q(iy)| &\sim y^{2\kappa+1}\ms S[\tau_\kappa](y^2)
		\le y^{2\kappa+1}\ms S[\tilde\tau](y^2)
		\sim y^{2\kappa+1}\frac{\pi}{2}\cdot\frac{\tilde\tau([0,y^2))}{y^2}
		\nonumber\\[1ex]
		&\sim \frac{\pi\sqrt{e}}{\sqrt{e}-1}y^{2\kappa}
		= \frac{\pi\sqrt{e}}{\sqrt{e}-1}y^{2\ell-2},
		\label{L221}
	\end{align}
	which shows that $|\Re q(iy)|\gg|\Im q(iy)|$.
	This, together with \eqref{L214}, implies that \eqref{L96}
	with $f(r)=2r^{2\ell-2}\log r$ is satisfied and therefore also
	(b) in \cref{L28} with $\beta=2\ell-1$ and $\omega=i(-1)^{\ell+1}$.
	Note also that (a)$'$ is satisfied with $\zeta=0$.
	
	For $\ell=1$ we obtain a function $q\in\mc N_0$ since $\kappa=0$ in this case.
	We can choose $p$ such that
	\[
		q(z) = \widetilde C[\mu](z) 
		= \sum_{k=1}^\infty\biggl(\frac{1}{e^{\frac{k}{2}}-z}-\frac{e^{\frac{k}{2}}}{1+e^k}\biggr)
		\cdot\frac{1+e^k}{e^{\frac{k}{2}}}.
	\]
	According to \eqref{L214} and \eqref{L221} it satisfies $q(iy)\sim2\log y$
	as $y\to\infty$.
\end{example}

\begin{example}\label{L53}
	Let $\gamma\in(0,1)$ and define the measure $\mu=\lambda+\nu$ on $\bb R$
	where $\lambda$ is the Lebesgue measure and $\nu$ is the measure such that
	$\nu((-\infty,0))=0$ and
	\[
		\nu([0,t)) = 
		\begin{cases}
			t, & t\in(0,e],
			\\[1ex]
			\dfrac{t}{(\log t)^\gamma}, & t\in(e,\infty).
		\end{cases}
	\]
	Further, set $q(z)\DE\widetilde C[\mu](z)$.
	Clearly, $\mu([0,t))\sim t$ and $\mu((-t,0))=t$ and hence $\mu((-t,t))\sim 2t$,
	which shows that (a) and (a)$'$ in \cref{L28} are satisfied with $\beta=1$ and $\zeta=1$.
	Further, \cref{L45} implies that $\Im q(iy)\to\pi$ as $y\to\infty$.
	For the real part we obtain from \cref{L163,L167} that
	\begin{align*}
		\Re q(iy) &= \RC_1[\mu](y)
		= \mc F_1[\lambda^++\nu](y)-\mc F_1[\lambda^-](y)
		= \mc F_1[\nu](y)
		\\[1ex]
		&\sim -\int_1^y\frac{\nu((-t,t))}{t^2}\DD t
		\sim -\int_1^y\frac{1}{t(\log t)^\gamma}\DD t
		= -\frac{1}{1-\gamma}(\log y)^{1-\gamma}.
	\end{align*}
	This shows that $q(iy)\sim -\frac{1}{1-\gamma}(\log y)^{1-\gamma}$ as $y\to\infty$,
	and hence (b) in \cref{L28} is satisfied.
	However, the asymptotics of $q(iy)$ at infinity is not determined by the leading term
	of the asymptotics of $\mu([0,t))$ or $\mu((-t,0))$
	as the former depends on $\gamma$ whereas the latter does not.
\end{example}

\subsection{The exceptional case}

In the following theorem we consider certain situations when $\beta\in2\bb N-1$ and $\zeta=1$,
a case that is not covered by \cref{L28}.
When $\beta\in2\bb N-1$ and $\zeta\ne1$, the real part of $q(iy)$ dominates
the imaginary part and is strictly larger than $\frac{\mu((-y,y))}{y}$;
see \cref{L28} and \cref{L157}\,(ii). 
In the next theorem we consider cases when there is cancellation between contributions 
from the measure $\mu$ on the positive and negative axes to the real part and 
where the growth of $q(iy)$ is determined by the imaginary part.

\begin{theorem}\label{L40}
	Let $\mu$, $\kappa$, $p$ and $q$ be as in \cref{L171,L173} and let $\beta\in2\bb N-1$.
	Then the following statements are true.
	\begin{Enumerate}
	\item
		Assume that $\mu=\mu_0+\mu_1+\mu_2$ with measures $\mu_0$, $\mu_1$ and $\mu_2$ 
		such that the following conditions are satisfied:
		\begin{Enumeratealph}
		\item
			$\mu_0$ is symmetric and $t\mapsto\mu_0((-t,t))$ is regularly varying with index $\beta$;
		\item
			either $\mu_1([0,\infty))=0$ or $t\mapsto\mu_1([0,t))$ 
			is regularly varying with index $\beta$;
			\\[0.3ex]
			either $\mu_1((-\infty,0))=0$ or $t\mapsto\mu_1((-t,0))$ 
			is regularly varying with index $\beta$;
		\item
			$\mu_2((-t,t))\lesssim t^\gamma$ as $t\to\infty$ with some $\gamma<\beta$;
		\item
			if $\int_{\bb R}(1+|t|)^{-\beta}\D\mu_1(t)<\infty$, then
			\[
				c_{2\kappa} = \int_{\bb R}\frac{t}{(1+t^2)^{\kappa+1}}\DD(\mu_1+\mu_2)(t);
			\]
		\item
			with the notation in \eqref{L164} the following two limits
			\textup{(}for $+$ and $-$\textup{)} exist in $\bb R$:
			\[
				\eta_\pm \DE
				\begin{cases}
					\displaystyle \lim_{t\to\infty}
					\biggl(\raisebox{0.5ex}{$\displaystyle\int_1^t\frac{\mu_1^\pm((-s,s))}{s^{\beta+1}}\DD s$}
					\bigg/\raisebox{-0.5ex}{$\displaystyle\frac{\mu_0((-t,t))}{t^\beta}$}\biggr)
					& \displaystyle \text{if} \ \int_{\bb R}\frac{\D\mu_1(t)}{(1+|t|)^\beta}=\infty,
					\\[4ex]
					\displaystyle -\lim_{t\to\infty}
					\biggl(\raisebox{0.5ex}{$\displaystyle\int_t^\infty\frac{\mu_1^\pm((-s,s))}{s^{\beta+1}}\DD s$}
					\bigg/\raisebox{-0.5ex}{$\displaystyle\frac{\mu_0((-t,t))}{t^\beta}$}\biggr)
					& \displaystyle \text{if} \ \int_{\bb R}\frac{\D\mu_1(t)}{(1+|t|)^\beta}<\infty.
				\end{cases}
			\]
		\end{Enumeratealph}
		Then $\mu((-t,t))\sim\mu_0((-t,t))$, $\beta=2\kappa+1$, and \eqref{L32} 
		holds locally uniformly for $z\in\bb C^+$ with
		\begin{equation}\label{L41}
			f(r) = \beta\frac{\mu((-r,r))}{r}, \qquad
			\omega = (-1)^\kappa\Bigl(\frac{\pi}{2}+i\bigl(\eta_+-\eta_-\bigr)\Bigr).
		\end{equation}
	\item
		Assume that there exist a regularly varying function $f:[x_0,\infty)\to(0,\infty)$ 
		with $x_0>0$ and a constant $\omega\notin i\bb R$ such that \eqref{L32}
		holds uniformly for $z\in\bb C^+$.
		Then $t\mapsto\mu((-t,t))$ is regularly varying with index $\beta$
		and \eqref{L80} holds with $\zeta=1$.
		One can choose $f$ as in \eqref{L41}, in which case $\omega$
		satisfies $\Re\omega=(-1)^\kappa\frac{\pi}{2}$.
	\end{Enumerate}
\end{theorem}

\begin{remark}\label{L57}
	\Cref{L212} shows that the condition $\omega\notin i\bb R$ in (ii) is essential.
\end{remark}

\begin{proof}[Proof of \cref{L40}]
We split the proof into a couple of steps; the first five steps deal with the proof of (i).
\begin{Steps}
\item
	First we note that either $\mu_1=0$ or $t\mapsto\mu_1((-t,t))$ is regularly varying
	with index $\beta$.  In the latter case we obtain from \cref{L36}\,(i) 
	and assumptions (b) and (e) that, 
	when $\int_{\bb R}(1+|t|)^{-\beta}\D\mu_1(t)=\infty$,
	\[
		\frac{\mu_1((-t,t))}{t^\beta} 
		\ll \int_1^t\frac{\mu_1((-s,s))}{s^{\beta+1}}\DD s
		= \int_1^t\frac{\mu_1^+((-s,s))}{s^{\beta+1}}\DD s
		+ \int_1^t\frac{\mu_1^-((-s,s))}{s^{\beta+1}}\DD s
		\lesssim \frac{\mu_0((-t,t))}{t^\beta}
	\]
	as $t\to\infty$.  A similar calculation and the use of \cref{L36}\,(ii)
	show that also in the case when $\int_{[0,\infty)}(1+|t|)^{-\beta}\D\mu_1(t)<\infty$,
	we have $\mu_1((-t,t))\ll\mu_0((-t,t))$.  
	Together with assumption (c), this implies that $\mu((-t,t))\sim\mu_0((-t,t))$.
	
	It follows from \cref{L45} that $\beta=2\kappa+1$ and
	\begin{equation}\label{L134}
		\Im q(iy) \sim \frac{\frac{\pi\beta}{2}}{\sin\bigl(\pi\kappa+\frac{\pi}{2}\bigr)}
		\cdot\frac{\mu((-y,y))}{y}
		= (-1)^\kappa\frac{\pi}{2}f(y).
	\end{equation}
\item
	Set $\ell\DE\ell(\mu)$ and $\widehat\mu\DE\mu_1+\mu_2$.
	By \cref{L154} we have either $\ell=\kappa$ or $\ell=\kappa+1$.
	It follows from assumption (c) and \cref{L159} 
	that $\int_{\bb R}(1+|t|)^{-\beta}\D\mu_2(t)<\infty$
	and hence $\ell(\mu_2)\le\kappa$.
	If $\ms p(\mu)$ is odd (i.e.\ $\ell=\kappa$), then
	\[
		\int_{\bb R}\frac{\D\mu_1(t)}{(1+|t|)^\beta} 
		\le \int_{\bb R}\frac{\D\mu(t)}{(1+|t|)^\beta} < \infty
	\]
	and hence, by the symmetry of $\mu_0$ and assumption (d),
	\[
		\int_{\bb R}\frac{t}{(1+t^2)^{\kappa+1}}\DD\mu(t)
		= \int_{\bb R}\frac{t}{(1+t^2)^{\kappa+1}}\DD\widehat\mu(t)
		= c_{2\kappa}.
	\]
	This shows that \cref{L176} is satisfied.
	It follows from \cref{L178} and the symmetry of $\mu_0$ that
	\begin{equation}\label{L107}
		\Re q(iy) = \RC_\ell[\mu](y) + \widetilde p(y)
		= \RC_\ell[\widehat\mu](y) + \widetilde p(y)
	\end{equation}
	with $\widetilde p$ as in \cref{L178}.
	
	Set $\hat\ell\DE\kappa$ when $\mu_1=0$ and $\hat\ell\DE\ell(\widehat\mu)$ otherwise.
	It follows in a similar way as above that, by assumption (b), either $\mu_1^+=0$
	or $\ell(\mu_1^+)=\kappa$ or $\ell(\mu_1^+)=\kappa+1$,
	and similarly for $\mu_1^-$.
	Since $\ell(\mu_2)\le\kappa$, we have $\kappa\le\hat\ell\le\ell\le\kappa+1$.
	
	Let us rewrite the expression in \eqref{L107} in the case when $\hat\ell=\kappa$
	and $\ell=\kappa+1$.
	We obtain from \eqref{L35} with $\kappa$ and $\kappa'$ there replaced by $\hat\ell-1$
	and $\ell-1$ respectively that
	\begin{align*}
		\frac{t}{t^2+y^2}\cdot\frac{(1-y^2)^\ell}{(1+t^2)^\ell}
		&= \Re\biggl[\frac{1}{t-iy}\cdot\frac{(1+(iy)^2)^\ell}{(1+t^2)^\ell}\biggr]
		\\[1ex]
		&= \Re\biggl[\frac{1}{t-iy}\cdot\frac{(1+(iy)^2)^{\hat\ell}}{(1+t^2)^{\hat\ell}}\biggr]
		- \Re\biggl[(t+iy)\cdot\frac{(1+(iy)^2)^{\hat\ell}}{(1+t^2)^{\hat\ell+1}}\biggr]
		\\[1ex]
		&= \frac{t}{t^2+y^2}\cdot\frac{(1-y^2)^{\hat\ell}}{(1+t^2)^{\hat\ell}}
		- (1-y^2)^{\hat\ell}\cdot\frac{t}{(1+t^2)^{\hat\ell+1}}.
	\end{align*}
	Hence
	\[
		\Re q(iy) = \RC_{\hat\ell}[\widehat\mu](y)
		- (1-y^2)^\kappa\int_{\bb R}\frac{t}{(1+t^2)^{\kappa+1}}\DD\widehat\mu(t) + \widetilde p(y)
		= \RC_{\hat\ell}[\widehat\mu])(y) + \widehat p(y),
	\]
	where the polynomial $\widehat p$ satisfies
	\[
		\widehat p(y) 
		= (-1)^\kappa\biggl[c_{2\kappa}
		-\int_{\bb R}\frac{t}{(1+t^2)^{\kappa+1}}\DD\widehat\mu(t)\biggr]y^{2\kappa}
		+ \BigO\bigl(y^{2\kappa-2}\bigr)
		= \BigO\bigl(y^{2\kappa-2}\bigr) = \BigO\bigl(y^{2\hat\ell-2}\bigr)
	\]
	by assumption (d).
	
	In all cases we can write
	\[
		\Re q(iy) = \RC_{\hat\ell}[\widehat\mu](y) + \widehat p(y)
		= \mc F_{\hat\ell}[\mu_1^+](y) - \mc F_{\hat\ell}[\mu_1^-](y)
		+ \RC_{\hat\ell}[\mu_2](y) + \widehat p(y)
	\]
	with an even polynomial $\widehat p$ of degree at most $2\hat\ell-2$;
	note that we have used \cref{L163} in the last step and that $\widehat p=\widetilde p$
	when $\hat\ell=\ell$.
	We therefore have
	\begin{equation}\label{L132}
		\frac{1}{i}\cdot\frac{q(iy)}{f(y)}
		= \frac{\Im q(iy)}{f(y)} 
		- i\frac{\mc F_{\hat\ell}[\mu_1^+](y)}{f(y)} + i\frac{\mc F_{\hat\ell}[\mu_1^-](y)}{f(y)}
		- i\frac{\RC_{\hat\ell}[\mu_2](y)}{f(y)} - i\frac{\widehat p(y)}{f(y)}.
	\end{equation}
\item
	We show that
	\begin{equation}\label{L135}
		\lim_{y\to\infty}\frac{\mc F_{\hat\ell}[\mu_1^\pm](y)}{f(y)} = -(-1)^\kappa\eta_\pm.
	\end{equation}
	Let us start with $\mu_1^+$.
	When $\mu_1^+=0$, the equality in \eqref{L135} is obvious.
	Assume now that $\mu_1^+\ne0$.
	Let us first consider the case when $\ell(\mu_1^+)=\hat\ell$.
	Then the assumptions of \cref{L167} are satisfied for $\mu$ replaced by $\mu_1^+$.
	Hence that lemma and assumption (e) yield that, 
	when $\int_{\bb R}(1+|t|)^{-\beta}\D\mu_1(t)=\infty$, i.e.\ $\ell=\kappa+1$,
	\[
		\lim_{y\to\infty}\frac{\mc F_{\hat\ell}[\mu_1^+](y)}{f(y)}
		= \lim_{y\to\infty}\biggl[\raisebox{0.5ex}{$\displaystyle (-1)^{\hat\ell}\beta y^{\beta-1}
		\int_1^y\frac{\mu_1^+((-s,s))}{s^{\beta+1}}\DD s$}
		\bigg/\raisebox{-0.5ex}{$\displaystyle \Bigl(\beta\frac{\mu_0((-y,y))}{y}\Bigr)$}\biggr]
		= -(-1)^\kappa\eta_+;
	\]
	a similar calculation proves relation \eqref{L135} also in the case 
	when $\int_{\bb R}(1+|t|)^{-\beta}\D\mu_1(t)<\infty$.
	The same considerations can be applied to $\mu_1^-$ when $\mu_1^-=0$
	or $\ell(\mu_1^-)=\hat\ell$.

	It remains to prove \eqref{L135} for $\mu_1^+$ 
	when $\ell(\mu_1^+)<\ell(\mu_1^-)=\hat\ell=\kappa+1$
	or for $\mu_1^-$ when $\ell(\mu_1^-)<\ell(\mu_1^+)$.
	We consider only the first case.
	It follows from \cref{L158}\,(ii) applied to $\mu_1^+$,
	the first inequality in \eqref{L156} and the already proved relation \eqref{L135}
	for $\mu_1^-$ (since $\ell(\mu_1^-)=\hat\ell$) that, with some $c>0$,
	\[
		\big|\mc F_{\hat\ell}[\mu_1^+](y)\big|
		\sim cy^{2\hat\ell-2} = cy^{2\kappa}
		\ll \big|\mc F_{\hat\ell}[\mu_1^-](y)\big|
		\lesssim f(y).
	\]
	On the other hand, the relations $\ell(\mu_1^-)=\kappa+1$,
	$\ell(\mu_1^+)<\kappa+1$ and \cref{L159} imply that
	\[
		\int_1^\infty\frac{\mu_1^+((-s,s))}{s^{\beta+1}}\DD s < \infty,
		\qquad
		\int_1^\infty\frac{\mu_1^-((-s,s))}{s^{\beta+1}}\DD s = \infty,
	\]
	and hence
	\[
		\int_1^t\frac{\mu_1^+((-s,s))}{s^{\beta+1}}\DD s 
		\ll \int_1^t\frac{\mu_1^-((-s,s))}{s^{\beta+1}}\DD s.
	\]
	This and the existence of the limit for $\eta_-$ yield $\eta_+=0$,
	which shows that \eqref{L135} holds also for $\mu_1^+$ in this case.
\item
	Next let us show that
	\begin{equation}\label{L136}
		\lim_{y\to\infty}\biggl(\frac{\RC_{\hat\ell}[\mu_2](y)}{f(y)} 
		+ \frac{\widehat p(y)}{f(y)}\biggr) = 0.
	\end{equation}
	Let us consider the two possible cases for $\hat\ell$ separately.
	First assume that $\hat\ell=\kappa$.
	Choose $\rho\in(\max\{\gamma,\beta-1\},\beta)$ and define the measure $\nu$ 
	by $\nu((-\infty,0))=0$ and $\nu([0,t))=t^\rho$ for $t>0$.
	Since $\rho\in(\beta-1,\beta)$, we have $\ell(\nu)=\hat\ell$.
	It follows from \cref{L155} and assumption (c) that
	\[
		\lim_{y\to\infty}\frac{\mc F_{\hat\ell}[\mu_2](y)}{\mc F_{\hat\ell}[\nu](y)}
		= \lim_{t\to\infty}\frac{\mu_2((-t,t))}{\nu((-t,t))} = 0.
	\]
	From \eqref{L166} and \eqref{L169} we obtain
	\[
		\big|\RC_{\hat\ell}[\mu_2](y)\big|
		\le \big|\mc F_{\hat\ell}[\mu_2](y)\big|
		\ll \big|\mc F_{\hat\ell}[\nu](y)\big|
		\sim \frac{\frac{\pi\rho}{2}}{\cos\frac{\pi\rho}{2}}y^{\rho-1}.
	\]
	Since $f$ is regularly varying with index $\beta-1=2\hat\ell$ and $\widehat p$
	is a polynomial of degree at most $2\hat\ell-2$, relation \eqref{L136}
	follows in this case; see \eqref{L130}.

	Now let us consider the case when $\hat\ell=\kappa+1$.
	Then $\ell(\mu_2)\le\kappa<\hat\ell$, and hence 
	\Cref{L163} and \cref{L158}\,(ii) imply that
	\[
		\big|\RC_{\hat\ell}[\mu_2](y)\big|
		\le \big|\mc F_{\hat\ell}[\mu_2](y)\big|
		\sim cy^{2\hat\ell-2}
	\]
	with some $c>0$.  
	Without loss of generality assume that $\ell(\mu_1^+)=\hat\ell$
	(the case $\ell(\mu_1^-)=\hat\ell$ is similar).
	It follows from assumption (e) and \cref{L167} that
	\[
		f(y) \gtrsim y^{\beta-1}\int_1^y\frac{\mu_1^+([0,s))}{s^{\beta+1}}\DD s
		\sim \frac{1}{\beta}\big|\mc F_{\hat\ell}[\mu_1^+](y)\big|
		\gg y^{2\hat\ell-2},
	\]
	which yields \eqref{L136} also in this case.
\item
	Combining \eqref{L132}, \eqref{L134}, \eqref{L135} and \eqref{L136}
	we arrive at
	\[
		\frac{1}{i}\lim_{y\to\infty}\frac{q(iy)}{f(y)}
		= (-1)^\kappa\Bigl(\frac{\pi}{2}+i\bigl(\eta_+-\eta_-\bigr)\Bigr),
	\]
	which, by \cref{L91}, shows the remaining assertions in (i).
\item
	Let us now prove item (ii).
	Assume that $f$ and $\omega\notin i\bb R$ are such that \eqref{L32} holds.
	We obtain from \cref{L16} and \cref{L108}\,(i) that
	\begin{equation}\label{L43}
		\lim_{t\to\infty}\frac{\mu((0,t))}{tf(t)}
		= \lim_{t\to\infty}\frac{\mu((-t,0))}{tf(t)}
		= \frac{\omega}{\pi(\alpha+1)}\cos\bigl(\arg\bigl((-1)^m\omega\bigr)\bigr)
	\end{equation}
	where $\alpha=\beta-1=2m$ with $m\in\bb N_0$ 
	and $\big|\arg\bigl((-1)^m\omega\bigr)\big|\le\frac{\pi}{2}$.
	Since, by assumption, $\omega\notin i\bb R$, we have $\Re((-1)^m\omega)>0$
	and hence
	\[
		|\omega|\cos\bigl(\arg\bigl((-1)^m\omega\bigr)\bigr) = \Re\bigl((-1)^m\omega)\bigr)
		= (-1)^m\Re\omega.
	\]
	In particular, the limits in \eqref{L43} are non-zero, and therefore \eqref{L80} holds
	with $\zeta=1$ and
	\begin{equation}\label{L44}
		\lim_{t\to\infty}\frac{\mu((-t,t))}{tf(t)} = (-1)^m\frac{2}{\pi\beta}\Re\omega.
	\end{equation}
	This implies that $t\mapsto\mu((-t,t))$ is regularly varying with index $\beta$,
	which, in turn, yields that $\beta=2\kappa+1$ and hence $m=\kappa$.
	Choosing $f$ as in \eqref{L41} we obtain from \eqref{L44} 
	that $\Re\omega=(-1)^\kappa\frac{\pi}{2}$.
\end{Steps}
\end{proof}

\medskip\noindent
To illustrate the Abelian direction of \cref{L40}, let us consider the following example.

\begin{example}\label{L51}
	Let $\mu$ be a measure on $\bb R$ that satisfies
	\begin{equation}\label{L55}
	\begin{aligned}
		\mu([0,t)) &= t\log t + a_+t + \BigO(t^\gamma),
		\\[1ex]
		\mu((-t,0)) &= t\log t + a_-t + \BigO(t^\gamma)
	\end{aligned}
	\end{equation}
	as $t\to+\infty$ with $a_+,a_-\ge 0$ and $\gamma<1$.
	We can write $\mu$ as $\mu=\mu_0+\mu_1+\mu_2$ where $\mu_0,\mu_1,\mu_2$
	are measures satisfying
	\begin{alignat*}{2}
		& \mu_0([-1,1]) = 0; \qquad &
		& \mu_0([1,t)) = \mu_0((-t,-1)) = t\log t, \qquad t>1;
		\\[1ex]
		& \mu_1([0,t)) = a_+t, \qquad && \mu_1((-t,0)) = a_-t, \qquad t>0;
		\\[1ex]
		& \mu_2((-t,t)) = \BigO(t^\gamma), \qquad && t\to+\infty.
	\end{alignat*}
	Since $\int_{\bb R}(1+|t|)^{-1}\D\mu(t)=\infty$ (see \cref{L159}), assumptions (a)--(d)
	in \cref{L40} are fulfilled with $\beta=1$.  For (e) we consider
	\[
		\eta_\pm = \lim_{t\to\infty}
		\biggl(\raisebox{0.5ex}{$\displaystyle\int_1^t\frac{a_\pm s}{s^2}\DD s$}
		\bigg/\raisebox{-0.5ex}{$\displaystyle\frac{2t\log t}{t}$}\biggr)
		= \frac{a_\pm}{2},
	\]
	which shows that (e) is also satisfied.
	Now let $q(z)=c_0+\widetilde C[\mu](z)$ with $c_0\in\bb R$.
	Then, by \cref{L40}\,(i), relation \eqref{L32} holds with
	\[
		f(r) \sim \frac{\mu_0((-r,r))}{r} = 2\log r, \qquad
		\omega = \frac{\pi}{2}+i\frac{a_+-a_-}{2},
	\]
	i.e.\
	\[
		\lim_{r\to\infty}\frac{q(rz)}{\log r} = -a_++a_-+i\pi
	\]
	locally uniform for $z\in\bb C^+$.
	An example of a function $q$ with a measure $\mu$ as in \eqref{L55} is
	\[
		q(z) = (a_--a_++i\pi)\log(z+i) + \pi(a_++1)i, \qquad z\in\bb C^+,
	\]
	where the $\BigO(t^\gamma)$ term is actually $\BigO(\log t)$.
\end{example}

% #########################################################
%
%
%\addtocontents{toc}{\contentsline{section}{\textbf{Appendix}}{}}
\appendix
\makeatletter
\DeclareRobustCommand{\@seccntformat}[1]{%
  \def\temp@@a{#1}%
  \def\temp@@b{section}%
  \ifx\temp@@a\temp@@b
  Appendix\ \csname the#1\endcsname.\quad%
  \else
  \csname the#1\endcsname\quad%
  \fi
}
\makeatother
\renewcommand{\thelemma}{\Alph{section}.\arabic{lemma}}
%
%
% #########################################################

%
%
%
\section[Regularly varying functions and some theorems of Karamata]{Regularly varying functions and \newline some theorems of Karamata}
\label{L115}

\noindent
In this appendix we provide some classical results about regularly varying functions
in slightly extended or rounded-off formulations.
A very good source for the theory of regular variation is \cite{bingham.goldie.teugels:1989};
this is our standard reference.

Recall the definition of regular variation in Karamata's sense.

\begin{definition}\label{L37}
	A function $f\DF[x_0,\infty)\to(0,\infty)$ with $x_0>0$ is called \emph{regularly varying}
	with \emph{index} $\alpha\in\bb R$ if it is measurable and
	\[
		\forall\lambda>0\DP \lim_{r\to\infty}\frac{f(\lambda r)}{f(r)} = \lambda^\alpha.
	\]
	A regularly varying function $f$ with index $0$ is also called \emph{slowly varying}.
\end{definition}

\medskip\noindent
Examples include functions $f$ behaving for large $r$ like
\[
	f(r) = r^\alpha\cdot\big(\log r\big)^{\beta_1}\cdot\big(\log\log r\big)^{\beta_2}
	\cdot\ldots\cdot
	\big(\underbrace{\log\cdots\log}_{\text{\footnotesize$m$\textsuperscript{th} iterate}}r\big)^{\beta_m},
\]
where $\alpha,\beta_1,\ldots,\beta_m\in\bb R$.  Other examples are
$f(r) = r^\alpha e^{(\log r)^\beta}$ with $\beta\in(0,1)$,
or $f(r) = r^\alpha e^{\frac{\log r}{\log\log r}}$;
see \cite[\S 1.3]{bingham.goldie.teugels:1989}.
In many respects, regularly varying functions of index $\alpha$ behave like the
power function $r^\alpha$.
For instance, if $f$ is regularly varying with index $\alpha$, then
\begin{equation}\label{L130}
	f(r) \gg r^{\alpha-\varepsilon}, \qquad
	f(r) \ll r^{\alpha+\varepsilon}
	\qquad\text{as} \ r\to\infty
\end{equation}
for every $\varepsilon>0$, which follows, e.g.\ from the Potter bounds;
see \cite[Theorem~1.5.6\,(iii)]{bingham.goldie.teugels:1989}.

Another property that reflects the power-like behaviour is a fundamental result 
by J.~Karamata about primitives of regularly varying functions.
We state a comprehensive formulation collecting what is proved
in \cite[Section~1.5.6]{bingham.goldie.teugels:1989}.
More precisely, item (i) in \cref{L36} follows from
\cite[Theorem~1.5.11\,(i) and Proposition~1.5.9a]{bingham.goldie.teugels:1989};
item (ii) follows from
\cite[Theorem~1.5.11\,(ii) and Proposition~1.5.9b]{bingham.goldie.teugels:1989}.

\begin{theorem}[Karamata]\label{L36}
	Let $x_0>0$ and let $f:[x_0,\infty)\to(0,\infty)$ be measurable and locally bounded.
	Further, assume that $f$ is regularly varying with index $\alpha\in\bb R$.
	\begin{Enumerate}
	\item
		Suppose that $\alpha+1\ge 0$.  Then the function
		$x\mapsto\int_{x_0}^x f(t)\DD t$
		is regularly varying with index $\alpha+1$, and
		\[
			\lim_{x\to\infty}
			\biggl(
			\raisebox{3pt}{$\displaystyle xf(x)$}
			\bigg/
			\raisebox{-2pt}{$\displaystyle\int_{x_0}^x f(t)\DD t$}
			\biggr)
			= \alpha+1.
		\]
	\item
		Suppose that $\int_{x_0}^\infty f(t)\DD t<\infty$.
		Then $\alpha+1\le 0$, the function
		$x\mapsto\int_x^\infty f(t)\DD t$
		is regularly varying with index $\alpha+1$, and
		\[
			\lim_{x\to\infty}\biggl(
			\raisebox{3pt}{$\displaystyle xf(x)$}
			\bigg/
			\raisebox{-2pt}{$\displaystyle\int_x^\infty f(t)\DD t$}
			\biggr)
			= -(\alpha+1).
		\]
	\end{Enumerate}
\end{theorem}

\medskip\noindent
In the following we often use integration by parts in its proper
measure-theoretic form as stated in the following lemma.

\begin{lemma}\label{L34}
	Let $-\infty<a<b\leq\infty$, let $\mu$ and $\nu$ be measures on $[a,b)$.
	Then
	\begin{equation}\label{L103}
		\int_{[a,b)}\mu([a,t))\DD\nu(t) = \int_{[a,b)}\nu((t,b))\DD\mu(t).
	\end{equation}
	If these integrals are finite, then $\lim_{x\to b}\mu([a,x))\nu([x,b))=0$.
\end{lemma}

\begin{proof}
	If $\nu([a,b))=\infty$, then either both sides are zero (when $\mu$ is the zero measure)
	or both sides are infinite (otherwise).
	In the case when $\nu$ is a finite measure, relation \eqref{L103} follows
	from Fubini's theorem.
	To show the last assertion, assume that both sides of \eqref{L103} are finite
	and let $x\in(a,b)$.  Then \eqref{L103} applied to $[a,x)$ instead of $[a,b)$ yields
	\begin{align*}
		\int_{[a,x)}\mu([a,t))\DD\nu(t)
		&= \int_{[a,x)}\nu((t,x))\DD\mu(t)
		\\
		&= \int_{[a,x)}\nu((t,b))\DD\mu(t) - \mu([a,x))\nu([x,b)).
	\end{align*}
	Letting $x\to b$ we obtain the claimed limit relation.
\end{proof}

\medskip\noindent
From \cref{L36} we obtain the next proposition, which relates integrals with respect
to a measure to integrals involving the corresponding distribution function.

\begin{proposition}\label{L33}
	Let $\sigma$ be a measure on $[1,\infty)$, which is not the zero measure.
	Assume that the distribution function $t\mapsto\sigma([1,t))$ is
	regularly varying with index $\alpha$
	\textup{(}since $t\mapsto\sigma([1,t))$ is non-decreasing, we have $\alpha\ge 0$\textup{)}.
	\begin{Enumerate}
	\item
		Let $\gamma\in\bb R\setminus\{0\}$.  Then
		\begin{equation}\label{L58}
			\lim_{x\to\infty}\biggl(
			\raisebox{3pt}{$\displaystyle\int_{[1,x)} t^\gamma\DD\sigma(t)$}
			\bigg/
			\raisebox{-2pt}{$\displaystyle\int_1^x t^{\gamma-1}\sigma([1,t))\DD t$}
			\biggr)
			= \max\{\alpha,-\gamma\}
			=
			\begin{cases}
				\alpha, & \alpha+\gamma\ge 0,
				\\[2mm]
				|\gamma|, & \alpha+\gamma<0.
			\end{cases}
		\end{equation}
		The function $x\mapsto\int_1^x t^{\gamma-1}\sigma([1,t))\DD t$ is
		regularly varying with index $\alpha+\gamma$ if $\alpha+\gamma\ge 0$,
		and non-decreasing and bounded \textup{(}in particular, slowly varying\textup{)}
		if $\alpha+\gamma<0$; if $\alpha\neq 0$ or $\gamma\leq 0$, 
		the same holds for $x\mapsto\int_{[1,x)}t^\gamma\DD\sigma(t)$.

		In particular, if $\int_{[1,\infty)} t^\gamma\DD\sigma(t)=\infty$, then $\alpha+\gamma\geq 0$.
		Moreover, if $\int_{[1,\infty)} t^\gamma\DD\sigma(t)<\infty$ and $\sigma$
		is an infinite measure, then $\alpha+\gamma\le0$.
	\item
		Let $\rho>0$ and assume that $\int_{[1,\infty)}t^{-\rho}\D\sigma(t)<\infty$.
		Then $\sigma([1,x))\ll x^\rho$, and hence $\alpha\leq\rho$.
		Further, we have
		\[
			\lim_{x\to\infty}\biggl(
			\raisebox{3pt}{$\displaystyle\int_{[x,\infty)}\frac{\D\sigma(t)}{t^\rho}$}
			\bigg/
			\raisebox{-2pt}{$\displaystyle\int_x^\infty\frac{\sigma([1,t))}{t^{\rho+1}}\DD t$}
			\biggr)
			= \alpha.
		\]
		The function $x\mapsto\int_x^\infty t^{-(\rho+1)}\sigma([1,t))\DD t$
		is regularly varying with index $\alpha-\rho$; if $\alpha\neq 0$, 
		the same holds for $x\mapsto\int_x^\infty t^{-\rho}\D\sigma(t)$.
	\end{Enumerate}
\end{proposition}

\begin{proof}
\hfill
\begin{Steps}
\item
	For the proof of \Enumref{1} we integrate by parts
	(using the measure $\DD\nu(t)\DE t^{\gamma-1}\DD t$ in \cref{L34})
	to obtain
	\begin{align}
		\int_{[1,x)}\sigma([1,t))t^{\gamma-1}\DD t
		&= \int_{[1,x)}\Bigl(\int_t^x s^{\gamma-1}\DD s\Bigr)\DD\sigma(t)
		\nonumber\\[1ex]
		&= \frac 1\gamma x^\gamma\sigma([1,x))
		- \frac{1}{\gamma}\int_{[1,x)}t^\gamma\DD\sigma(t),
		\label{L104}
	\end{align}
	and therefore
	\begin{equation}\label{L133}
		\raisebox{3pt}{$\displaystyle \int_{[1,x)}t^\gamma\DD\sigma(t)$}
		\bigg/
		\raisebox{-2pt}{$\displaystyle \int_{[1,x)}t^{\gamma-1}\sigma([1,t))\DD t$}
		=
		\biggl(\raisebox{3pt}{$\displaystyle x^\gamma\sigma([1,x))$}
		\bigg/
		\raisebox{-2pt}{$\displaystyle \int_{[1,x)}t^{\gamma-1}\sigma([1,t))\DD t$}
		\biggr)-\gamma.
	\end{equation}

\item
	First we consider the case when $\alpha+\gamma\ge 0$.
	We can apply \cref{L36}\,(i) to obtain that the integral
	$\int_{[1,x)}\sigma([1,t))t^{\gamma-1}\DD t$ is regularly varying with
	index $\alpha+\gamma$, and that the quotient on the right-hand side
	of \eqref{L133} tends to $\alpha+\gamma$.  The asserted limit relation follows.
	In particular, if $\alpha\ne 0$, also the
	integral $\int_1^x t^\gamma\DD\sigma(t)$ is regularly varying
	with index $\alpha+\gamma$.

\item
	Now assume that $\alpha+\gamma<0$. Since $\alpha\ge 0$, we have $\gamma<0$
	and hence $\max\{\alpha,-\gamma\}=-\gamma=|\gamma|$.
	The integral $\int_1^\infty t^{\gamma-1}\sigma([1,t))\DD t$ converges,
	and hence $\lim_{x\to\infty}x^\gamma\sigma([1,x))=0$ and 
	by \eqref{L104}
	\[
		\int_{[1,\infty)}t^\gamma\DD\sigma(t)
		= |\gamma|\int_{[1,\infty)}t^{\gamma-1}\sigma([1,t))\DD t
		< \infty.
	\]
	In particular, the function $x\mapsto\int_{[1,x)}t^\gamma\DD\sigma(t)$ is
	slowly varying and the asserted limit relation holds.

\item
	For the last statement in (i) assume that $\int_{[0,\infty)}t^\gamma\DD\sigma(t)<\infty$,
	that $\sigma$ is an infinite measure, and suppose that $\alpha+\gamma>0$.
	If $\alpha>0$, then $\int_{[1,x)}t^\gamma\DD\sigma(t)$ is regularly varying
	with positive index $\alpha+\gamma$ and hence unbounded, a contradiction.
	If $\alpha=0$, then $\gamma>0$, and therefore
	$\int_{[1,x)}t^\gamma\DD\sigma(t)\ge\sigma([1,x))\to\infty$ as $x\to\infty$,
	again a contradiction.
\item
	For the proof of \Enumref{2} we argue in a similar way.  Integrate by parts
	(using the measure $\D\nu(t)=\frac{\D t}{t^{\rho+1}}$) to obtain
	\[
		\int_1^\infty\sigma([1,t))\frac{\D t}{t^{\rho+1}}
		= \frac 1\rho\int_{[1,\infty)}\frac{\D\sigma(t)}{t^\rho},
		\qquad
		\lim_{x\to\infty}\frac{\sigma([1,x))}{x^\rho} = 0.
	\]
	The second relation shows, in particular, that $\alpha\le\rho$.
	We integrate by parts again to obtain
	\begin{align*}
		\int_x^\infty\sigma([1,t))\frac{\D t}{t^{\rho+1}}
		&= \int_x^\infty\sigma([1,x))\frac{\D t}{t^{\rho+1}}
		+ \int_x^\infty\sigma([x,t))\frac{\D t}{t^{\rho+1}}
		\\[1ex]
		&= \frac{1}{\rho}\cdot\frac{\sigma([1,x))}{x^\rho}
		+\frac{1}{\rho}\int_{[x,\infty)}\frac{\D\sigma(t)}{t^\rho},
	\end{align*}
	and hence
	\[
		\raisebox{3pt}{$\displaystyle \int_{[x,\infty)}\frac{\D\sigma(t)}{t^\rho}$}
		\bigg/
		\raisebox{-2pt}{$\displaystyle \int_x^\infty\frac{\sigma([1,t))}{t^{\rho+1}}\DD t$}
		= \rho - \bigg(
		\raisebox{3pt}{$\displaystyle\frac{\sigma([1,x))}{x^\rho}$}
		\bigg/
		\raisebox{-2pt}{$\displaystyle \int_x^\infty\frac{\sigma([1,t))}{t^{\rho+1}}\DD t$}
		\bigg).
	\]
	By \cref{L36}\,(ii) the integral $\int_{[1,x)}t^{-(\rho+1)}\sigma([1,t))\DD t$
	is regularly varying with index $\alpha-\rho$, and the quotient on the
	right-hand side tends to $\rho-\alpha$.  The assertions made in (ii) follow.
\end{Steps}
\end{proof}

\medskip\noindent
The following example shows that, when $\alpha=0$ and $\gamma>0$ in \cref{L33}\,(i), 
the function $x\mapsto\int_{[1,x)}t^\gamma\DD\sigma(t)$ may fail to
be regularly varying.
Instead of $t^\gamma$ we consider an arbitrary function $g$ with $g(t)\sim t^\gamma$
as $t\to\infty$.  This more general example is used in \cref{L202}.

\begin{example}\label{L69}
	Define the discrete measure $\sigma$ supported on $\{e^k\DS k\in\bb N\}$ 
	with the following point masses:
	\[
		\sigma(\{e^k\}) = 1, \qquad k\in\bb N.
	\]
	The distribution function satisfies $\sigma([1,t))=0$ if $t\le e$
	and $\sigma([1,t))=n$ if $t\in(e^n,e^{n+1}]$ for $n\in\bb N$.
	The relation $t\in(e^n,e^{n+1}]$ is
	equivalent to $\log t-1\le n<\log t$, and hence $\sigma([1,t))\sim\log t$ 
	as $t\to\infty$.
	This shows that the distribution function $t\mapsto\sigma([1,t))$ 
	is slowly varying.
	
	Now let $g:[1,\infty)\to(0,\infty)$ be a function such that
	$g(t)\sim t^\gamma$ as $t\to\infty$ with $\gamma>0$, and define the function
	\[
		f(x) \DE \int_{[1,x)}g(t)\DD\sigma(t), \qquad x\in(1,\infty).
	\]
	For $x\in(e^{n-1},e^n]$ we have $f(x) = \sum_{k=1}^{n-1}g(e^k)$.
	It follows easily that, as $n\to\infty$,
	\begin{equation}\label{L203}
		\sum_{k=1}^{n-1}g(e^k) \sim \sum_{k=1}^{n-1}e^{\gamma k}
		= \frac{e^{\gamma n}-e^\gamma}{e^\gamma-1} \sim \frac{e^{\gamma n}}{e^\gamma-1}.
	\end{equation}
	Now choose $\lambda=\sqrt{e}$. % let $\lambda\in(1,c)$ and choose $\xi\in\bigl(1,\frac{c}{\lambda}\bigr)$.
	From \eqref{L203} we obtain
	\begin{align*}
		\frac{f(\lambda e^m)}{f(e^m)} &= \frac{\sum_{k=1}^m g(e^k)}{\sum_{k=1}^{m-1}g(e^k)}
		\to e^\gamma
		\qquad \text{as}\;\;m\to\infty,
		\\[1ex]
		\frac{f(\lambda e^{m+\frac12})}{f(e^{m+\frac12})} 
		&= \frac{\sum_{k=1}^m g(e^k)}{\sum_{k=1}^m g(e^k)}
		= 1,
	\end{align*}
	which implies that $f$ is not regularly varying.
	
	Note that $f(x)\asymp x^\gamma$ as $x\to\infty$, which can be seen from \eqref{L203}.
\end{example}

\medskip\noindent
The next topic we discuss is the asymptotic behaviour of Stieltjes transforms:
let $\mu$ be a measure on $[0,\infty)$ which 
satisfies $\int_{[0,\infty)}(1+t)^{-1}\DD\mu(t)<\infty$;
then the Stieltjes transform of $\mu$ is defined by
\begin{equation}\label{L22}
	\ms S[\mu](x)\DE \int_{[0,\infty)}\frac{\D\mu(t)}{t+x},\qquad x>0.
\end{equation}
As is common practice in the literature, we use a sign convention that is different
from the one in \eqref{L300} (i.e.\ different from the convention for the Cauchy transform) 
in order to obtain functions that are defined on the positive half-line.

\begin{remark}
\label{L26}
	Let $\mu$ be a measure on $[0,\infty)$
	which satisfies $\int_{[0,\infty)}(1+t)^{-1}\DD\mu(t)<\infty$. 
	Then, by the dominated and monotone convergence theorems, 
	\[
		\lim_{x\to\infty}\ms S[\mu](x) = 0, \qquad
		\lim_{x\to\infty}x\ms S[\mu](x) = \mu\bigl([0,\infty)\bigr).
	\]
	In particular, if $\mu$ is not the zero measure, then, as $x\to\infty$, 
	\begin{alignat*}{2}
		& \frac{1}{x} \asymp \ms S[\mu](x) \ll 1 \quad&& \text{if $\mu$ is finite},
		\\
		& \frac{1}{x} \ll \ms S[\mu](x) \ll 1 \quad&& \text{if $\mu$ is infinite}.
	\end{alignat*}
\end{remark}

\medskip\noindent
Karamata's theorem about the Stieltjes transform \cite{karamata:1931a,karamata:1931} characterises regular variation of 
the Stieltjes transform and gives precise information about the size of $\ms S[\mu](x)$ also when $\mu$ is infinite. 
We use it in a formulation that includes a boundary case;
this is often excluded, e.g.\ in \cite[Theorem~1.7.4]{bingham.goldie.teugels:1989}
or \cite{seneta:1976}.

\begin{theorem}[Karamata]\label{L1}
	Let $\mu$ be a measure on $[0,\infty)$, which is not the
	zero measure and satisfies $\int_{[0,\infty)}(1+t)^{-1}\DD\mu(t)<\infty$.
	Then the following two statements are equivalent:
	\begin{Enumerate}
	\item
		the distribution function $t\mapsto\mu([0,t))$ is regularly varying
		with index $\alpha$;
	\item
		$\ms S[\mu]$ is regularly varying with index $\alpha-1$.
	\end{Enumerate}
	If \textup{(i)} and \textup{(ii)} hold, then $\alpha\in[0,1]$ and
	\begin{equation}\label{L105}
		\ms S[\mu](x) \sim C_\alpha\int_x^\infty \frac{\mu([0,t))}{t^2}\DD t,
		\qquad x\to\infty.
	\end{equation}
	with
	\begin{equation}\label{L129}
		C_\alpha \DE
		\begin{cases}
			\dfrac{\pi\alpha(1-\alpha)}{\sin(\pi\alpha)}, & \alpha\in(0,1);
			\\[1ex]
			1, & \alpha\in\{0,1\}.
		\end{cases}
	\end{equation}
	The integral in \eqref{L105} is finite for every $x>0$.
\end{theorem}

\medskip\noindent
As in the usual presentations in textbooks we follow the lines
of \cite{karamata:1931a} and reduce the problem to the Laplace--Stieltjes transform.
Recall that the Laplace--Stieltjes transform of a positive measure $\nu$ on $[0,\infty)$
is the function $\ms L[\nu]\DF\bb R\to[0,\infty]$ defined as 
\[
	\ms L[\nu](x)\DE\int_{[0,\infty)}e^{-xt}\DD\nu(t),
	\qquad x\in\bb R.
\]
In the proof of \cref{L1} we also need the concept of a regularly varying function at $0$:
a function $g:(0,x_0]\to(0,\infty)$ with $x_0>0$ is called 
\emph{regularly varying at $0$ with index $\beta$}
if $x\mapsto g\bigl(\frac{1}{x}\bigr)$ is regularly varying with index $-\beta$
in the sense of \cref{L37}.

\begin{proof}[Proof of \cref{L1}]
\hfill
	\begin{Steps}
	\item
		The relation with the Stieltjes transform is established by Fubini's theorem:
		for $x>0$ we have
		\begin{align}
			\ms S[\mu](x) &= \int_{[0,\infty)}\frac 1{t+x}\DD\mu(t)
			= \int_{[0,\infty)}\bigg(\int_{[0,\infty)}
			e^{-(t+x)s}\DD s\bigg)\DD\mu(t)
			\nonumber\\[1ex]
			&= \int_{[0,\infty)}e^{-xs}\bigg(\int_{[0,\infty)}
			e^{-ts}\DD\mu(t)\bigg)\DD s
			%= \ms L\big[\ms L[\nu](s)\D s\big](x)
			= \int_{[0,\infty)}e^{-xs}\,\ms L[\mu](s)\DD s.
			\label{L150}
		\end{align}
		Let $\sigma$ be the measure on $[0,\infty)$ with density $\ms L[\mu]$, i.e.\
		\[
			\sigma([0,t)) \DE \int_0^t \ms L[\mu](s)\DD s, \qquad t>0.
		\]
		The latter integral is finite since, again by Fubini's theorem,
		\[
			\int_0^t \ms L[\mu](s)\DD s
			= \int_0^t \int_{[0,\infty)}e^{-sr}\DD\mu(r)\DD s
			= \int_{[0,\infty)}\frac{1-e^{-tr}}{r}\DD\mu(r);
		\]
		the finiteness of the last integral follows easily from the
		assumption $\int_{[0,\infty)}(1+t)^{-1}\D\mu(t)<\infty$.
		Now \eqref{L150} can be written as $\ms S[\mu]=\ms L[\sigma]$.
	\item
		Assume first that $x\mapsto\mu([0,x))$ is regularly varying with index $\alpha$.
		Since $x\mapsto\mu([0,x))$ is non-decreasing, $\alpha\ge0$;
		by \cref{L33}\,(ii) we have $\alpha\le1$.
		It follows from \cite[Theorem~1.7.1]{bingham.goldie.teugels:1989}
		in the form of \cite[(1.7.3)]{bingham.goldie.teugels:1989} that
		\begin{equation}\label{L128}
			\ms L[\mu]\Bigl(\frac 1x\Bigr)
			\sim \Gamma(1+\alpha)\mu([0,x)),
			\qquad x\to\infty.
		\end{equation}
		In particular, the function $x\mapsto\ms L[\mu](\frac 1x)$ is regularly varying
		with index $\alpha$.
		Note that in \cite{bingham.goldie.teugels:1989} the right-continuous
		distribution function $x\mapsto\mu([0,x])$ is used.
		However, the latter is regularly varying if and only if
		$x\mapsto\mu([0,x))$ is regularly varying and
		$\mu([0,x])\sim\mu([0,x))$, $x\to\infty$ if these functions
		are regularly varying.

		Before we apply the Laplace--Stieltjes transform a second time,
		we need the asymptotic behaviour of
		$x\mapsto\sigma\bigl(\bigl[0,\frac{1}{x}\bigr)\bigr)$ as $x\to\infty$.
		From \eqref{L128} we obtain
		\[
			\sigma\bigl(\bigl[0,\tfrac{1}{x}\bigr)\bigr)
			= \int_0^{\frac 1x}\ms L[\mu](s)\DD s
			= \int_x^\infty \ms L[\mu]\Bigl(\frac 1t\Bigr)\frac 1{t^2}\DD t
			\sim \Gamma(1+\alpha)\int_x^\infty\frac{\mu([0,t))}{t^2}\DD t,
			\qquad x\to\infty.
		\]
		It follows from \cref{L33}\,(ii) that the function
		$x\mapsto\sigma\bigl(\bigl[0,\frac{1}{x}\bigr)\bigr)$ is regularly varying
		with index $\alpha-1$,
		and hence, $t\mapsto\sigma([0,t))$ is regularly varying at $0$
		with index $1-\alpha$.
		Now \cite[Theorem~1.7.1$'$]{bingham.goldie.teugels:1989} implies that
		\begin{align*}
			\ms S[\mu](x) &= \ms L[\sigma](x)
			\sim \Gamma(2-\alpha)\sigma\bigl(\bigl[0,\tfrac{1}{x}\bigr)\bigr)
			\\[1ex]
			&\sim \Gamma(1+\alpha)\Gamma(2-\alpha)
			\int_{[x,\infty)}\frac{\mu([0,t))}{t^2}\DD t,
			\qquad x\to\infty,
		\end{align*}
		and that $\ms S[\mu]$ is regularly varying (at infinity) with index $\alpha-1$.
		It follows from the reflection formula for the Gamma function that
		$\Gamma(1+\alpha)\Gamma(2-\alpha)=\alpha(1-\alpha)\Gamma(\alpha)\Gamma(1-\alpha)=\frac{\pi(1-\alpha)\alpha}{\sin(\pi\alpha)}$
		when $\alpha\in(0,1)$.
	\item
		Conversely, assume that (ii) holds.
		Again by \cite[Theorem~1.7.1$'$]{bingham.goldie.teugels:1989},
		the function $t\mapsto\sigma([0,t))$ is regularly varying at $0$ with index $1-\alpha$.
		Since $\ms L[\mu]$ is non-increasing, we can apply
		\cite[Theorem~1.7.2b]{bingham.goldie.teugels:1989} to deduce that $\ms L[\mu]$
		is regularly varying at $0$ with index $-\alpha$.
		Finally, we obtain from \cite[Theorem~1.7.1]{bingham.goldie.teugels:1989}
		that (i) holds.
	\end{Steps}
\end{proof}

\begin{remark}\label{L3}
	In the case when $\alpha<1$ we can use \cref{L36}\,(ii) to rewrite
	the right-hand side of \eqref{L105} to obtain the standard formulation
	as in \cite[Theorem~1.7.4]{bingham.goldie.teugels:1989};
	namely, under the assumption of \cref{L1} we have
	\[
		\ms S[\mu](x) \sim \frac{\pi\alpha}{\sin(\pi\alpha)}
		\cdot\frac{\mu([0,x))}{x},
		\qquad x\to\infty,
	\]
	where the first fraction on the right-hand side is understood as $1$ if $\alpha=0$.
\end{remark}

\medskip\noindent
In the proofs of our main results we often need the following elementary facts 
where we compare the Stieltjes transforms of two measures.

\begin{lemma}\label{L54}
	Let $\mu_1$ and $\mu_2$ be measures on $[0,\infty)$
	such that $\mu_2$ is not the zero measure and that
	$\int_{[0,\infty)}(1+t)^{-1}\D\mu_i(t)<\infty$ for $i\in\{1,2\}$.  
	\begin{Enumerate}
	\item
		If $\mu_1([0,t))\le\mu_2([0,t))$ for all $t\in(0,\infty)$, then
		$\ms S[\mu_1](x) \le \ms S[\mu_2](x)$ for all $x\in(0,\infty)$.
	\item
		We have
		\begin{equation}\label{L106}
			\limsup_{x\to\infty}\frac{\ms S[\mu_1](x)}{\ms S[\mu_2](x)}
			\le \limsup_{r\to\infty}\frac{\mu_1([0,r))}{\mu_2([0,r))}.
		\end{equation}
	\end{Enumerate}
\end{lemma}

\begin{proof}
	We use \cref{L34} with the measure $\D\nu(t)\DE\frac{\D t}{(t+x)^2}$
	for $x>0$ to rewrite the Stieltjes transforms:
	\begin{equation}\label{L126}
		\ms S[\mu_i](x) = \int_{[0,\infty)}\frac{1}{t+x}\DD\mu_i(t)
		= \int_0^\infty\frac{\mu_i([0,t))}{(t+x)^2}\DD t,
		\qquad x>0,\; i\in\{1,2\},
	\end{equation}
	which immediately yields the assertion in (i).
	
	Let us now prove the statement in (ii).
	If $\mu_1=0$ or the right-hand side of \eqref{L106} is infinite,
	then there is nothing to prove.  Hence we assume that $\mu_1$ is not the zero measure
	and the right-hand side of \eqref{L106} is finite.
	Let $r_0>0$ be such that $\mu_i([0,r_0))>0$ for $i\in\{1,2\}$.
	From \eqref{L126} we obtain,
	for $x>0$ (and the asymptotic relations as $x\to\infty$),
	\begin{align*}
		\ms S[\mu_1](x)
		&= \int_0^\infty\frac{\mu_1([0,t))}{(t+x)^2}\DD t
		= \underbrace{\int_0^{r_0}\frac{\mu_1([0,t))}{(t+x)^2}\DD t}_{\le\frac{1}{x^2}\cdot r_0\mu_1([0,r_0))}
		+ \underbrace{\int_{r_0}^\infty\frac{\mu_1([0,t))}{(t+x)^2}\DD t}_{\ge\mu_1([0,r_0))\cdot\frac{1}{r_0+x}}
		\\[1ex]
		&\sim \int_{r_0}^\infty\frac{\mu_1([0,t))}{(t+x)^2}\DD t
		\le \sup_{t\ge r_0}\frac{\mu_1([0,t))}{\mu_2([0,t))}
		\cdot\int_{r_0}^\infty\frac{\mu_2([0,t))}{(t+x)^2}\DD t
		\\[1ex]
		&\sim \sup_{t\ge r_0}\frac{\mu_1([0,t))}{\mu_2([0,t))}
		\cdot\int_0^\infty\frac{\mu_2([0,t))}{(t+x)^2}\DD t
		= \sup_{t\ge r_0}\frac{\mu_1([0,t))}{\mu_2([0,t))}\cdot\ms S[\mu_2](x).
	\end{align*}
	With this we can deduce that
	$\limsup\limits_{x\to\infty}\frac{\ms S[\mu_1](x)}{\ms S[\mu_2](x)}\le\sup\limits_{t\ge r_0}\frac{\mu_1([0,t))}{\mu_2([0,t))}$.
	Since $r_0$ was arbitrary, the assertion follows.
\end{proof}

\medskip\noindent
We also need a comparison result when we integrate powers with respect
to two different measures.

\begin{lemma}\label{L59}
	Let $\nu_1,\nu_2$ be measures on $[0,\infty)$
	and let $\gamma\in\bb R\setminus\{0\}$.
	\begin{Enumerate}
	\item
		Assume that $\gamma<0$ and that $\nu_1([0,x))\lesssim\nu_2([0,x))$ as $x\to\infty$.
		Then
		\begin{equation}\label{L79}
			\int_{[1,x)}t^\gamma\DD\nu_1(t) \lesssim \int_{[1,x)}t^\gamma\DD\nu_2(t),
			\qquad x\to\infty.
		\end{equation}
		If, in addition, $\int_{[1,\infty)}t^\gamma\DD\nu_2(t)=\infty$, then
		\begin{equation}\label{L66}
			\limsup_{x\to\infty}\frac{\int_{[1,x)}t^\gamma\DD\nu_1(t)}{\int_{[1,x)}t^\gamma\DD\nu_2(t)}
			\le \limsup_{x\to\infty}\frac{\nu_1([0,x))}{\nu_2([0,x))}.
		\end{equation}
	\item
		Assume that $\gamma>0$ and that $x\mapsto\nu_2([0,x))$ is regularly varying
		with index $\alpha>0$.  Then
		\begin{equation}\label{L67}
			\limsup_{x\to\infty}\frac{\int_{[1,x)}t^\gamma\DD\nu_1(t)}{\int_{[1,x)}t^\gamma\DD\nu_2(t)}
			\le \frac{\alpha+\gamma}{\alpha}
			\limsup_{x\to\infty}\frac{\nu_1([0,x))}{\nu_2([0,x))}.
		\end{equation}
	\item
		Assume that $\gamma>0$ and that $x\mapsto\nu_i([0,x))$ are regularly varying
		with index $\alpha_i>0$ for $i\in\{1,2\}$.  Then
		\[
			\frac{\int_{[1,x)}t^\gamma\DD\nu_1(t)}{\int_{[1,x)}t^\gamma\DD\nu_2(t)}
			\sim \frac{\alpha_1(\alpha_2+\gamma)}{\alpha_2(\alpha_1+\gamma)}
			\cdot\frac{\nu_1([0,x))}{\nu_2([0,x))},
			\qquad x\to\infty.
		\]
	\end{Enumerate}
\end{lemma}

\begin{proof}
	If the limit superior on the right-hand side of \eqref{L67}
	is $+\infty$, then there is nothing to prove for (ii).
	Hence, for the proof (ii) we assume that the right-hand side 
	of \eqref{L67} is finite.
	In (i) the right-hand side of \eqref{L66} is finite by assumption.
	Let $M>0$ and $x_0>1$ be such that
	\begin{equation}\label{L62}
		\forall x\ge x_0\DP \frac{\nu_1([0,x))}{\nu_2([0,x))} \le M.
	\end{equation}

	(i) Assume that $\gamma<0$.  For $i\in\{1,2\}$, integration by parts yields
	\begin{align*}
		\int_{[1,x)}t^\gamma\DD\nu_i(t)
		&= x^\gamma\nu_i\bigl([0,x)\bigr) - \nu_i\bigl([0,1)\bigr)
		- \gamma\int_1^x t^{\gamma-1}\nu_i([0,t))\DD t
		\\[1ex]
		&= x^\gamma\nu_i\bigl([0,x)\bigr) - \nu_i\bigl([0,1)\bigr)
		+ |\gamma|\int_1^{x_0}t^{\gamma-1}\nu_i([0,t))\DD t
		+ |\gamma|\int_{x_0}^x t^{\gamma-1}\nu_i([0,t))\DD t
		\\[1ex]
		&= x^\gamma\nu_i\bigl([0,x)\bigr)
		+ |\gamma|\int_{x_0}^x t^{\gamma-1}\nu_i([0,t))\DD t + c_i
	\end{align*}
	with some $c_i\in\bb R$.  Together with \eqref{L62} we obtain, for $x\ge x_0$,
	\begin{align*}
		\int_{[1,x)}t^\gamma\DD\nu_1(t)
		&= x^\gamma\nu_1\bigl([0,x)\bigr)
		+ |\gamma|\int_{x_0}^x t^{\gamma-1}\nu_1([0,t))\DD t + c_1
		\\[1ex]
		&\le Mx^\gamma\nu_2\bigl([0,x)\bigr)
		+ M|\gamma|\int_{x_0}^x t^{\gamma-1}\nu_2([0,t))\DD t + c_1
		\\[1ex]
		&= M\int_{[1,x)}t^\gamma\DD\nu_2(t) - Mc_2 + c_1.
	\end{align*}
	This proves \eqref{L79}.  
	Now assume that $\int_{[1,\infty)}t^\gamma\DD\nu_2(t)=\infty$.
	Then
	\[
		\frac{\int_{[1,x)}t^\gamma\DD\nu_1(t)}{\int_{[1,x)}t^\gamma\DD\nu_2(t)}
		\le M + \frac{c_1-Mc_2}{\int_{[1,x)}t^\gamma\DD\nu_2(t)}
		\to M, \qquad x\to\infty,
	\]
	from which we can deduce that
	\[
		\limsup_{x\to\infty}
		\frac{\int_{[1,x)}t^\gamma\DD\nu_1(t)}{\int_{[1,x)}t^\gamma\DD\nu_2(t)}
		\le \sup_{t\ge x_0}\frac{\nu_1([0,t))}{\nu_2([0,t))}.
	\]
	Since $x_0$ was arbitrary, the inequality in \eqref{L66} follows.

	(ii)
	Now we assume that $\gamma>0$ and that $x\mapsto\nu_2([0,x))$ is regularly
	varying with index $\alpha$.
	The latter, together with \cref{L33}\,(i) and \cref{L36}\,(i), implies
	\begin{equation}\label{L124}
		\int_{[1,x)}t^\gamma\DD\nu_2(t)
		\sim \frac{\alpha}{\alpha+\gamma}x^\gamma\nu_2([1,x)),
		\qquad x\to\infty.
	\end{equation}
	The monotonicity of $t\mapsto t^\gamma$ and \eqref{L62} yield, for $x\ge x_0$,
	\begin{align}
		\int_{[1,x)}t^\gamma\DD\nu_1(t) &\le x^\gamma\nu_1\bigl([1,x)\bigr)
		\le x^\gamma\nu_1\bigl([0,x)\bigr)
		\le Mx^\gamma\nu_2\bigl([0,x)\bigr)
		\nonumber\\[1ex]
		&= Mx^\gamma\nu_2\bigl([1,x)\bigr)\biggl(1+\frac{\nu_2([0,1))}{\nu_2([1,x))}\biggr).
		\label{L116}
	\end{align}
	Since, by assumption, $\alpha>0$, we have $\nu_2([1,x))\to\infty$.
	Together with \eqref{L124} and \eqref{L116}, this shows that
	\[
		\limsup_{x\to\infty}
		\frac{\int_{[1,x)}t^\gamma\DD\nu_1(t)}{\int_{[1,x)}t^\gamma\DD\nu_2(t)}
		\le M\frac{\alpha+\gamma}{\alpha}.
	\]
	Now the assertion follows as in the proof of (i).

	(iii)
	In the same way as in the proof of (ii), we obtain from \cref{L33}\,(i) 
	and \cref{L36}\,(i) that
	\[
		\int_{[1,x)}t^\gamma\DD\nu_i(t)
		\sim \frac{\alpha_i}{\alpha_i+\gamma}x^\gamma\nu_i([1,x)),
		\qquad x\to\infty,
	\]
	for $i\in\{1,2\}$.  From this the result follows since $\nu_i([1,x))\sim\nu_i([0,x))$.
\end{proof}

\medskip\noindent
The next proposition contains asymptotic results about Stieltjes transforms
of certain measures; it plays a key role in the proofs of some
of the main results of the paper.

\begin{proposition}\label{L122}
	Let $\nu$ be a measure on $[0,\infty)$.
	Further, let $h:[0,\infty)\to[0,\infty)$ be a continuous function
	such that $h(t)>0$ for $t>0$ and $h(t)\sim t^\gamma$, $t\to\infty$, 
	with $\gamma\in\bb R\setminus\{0\}$.
	Assume that
	\begin{equation}\label{L118}
		\int_{[0,\infty)}h(t)\DD\nu(t) = \infty
		\qquad\text{and}\qquad
		\int_{[0,\infty)}\frac{h(t)}{1+t}\DD\nu(t) < \infty
	\end{equation}
	and define the measure $\sigma$ on $[0,\infty)$ by 
	\[
		\D\sigma(t)=h(t)\DD\nu(t), \qquad t\in[0,\infty).
	\]
	Then the Stieltjes transform $\ms S[\sigma]$ is well defined and
	satisfies
	\begin{equation}\label{L189}
		\frac{1}{x} \ll \ms S[\sigma](x) \ll 1, \qquad x\to\infty.
	\end{equation}
	Now let $\alpha\ge0$ and consider the following two statements:
	\begin{Enumeratealph}
	\item
		$t\mapsto\nu([0,t))$ is regularly varying with index $\alpha$;
	\item
		$\ms S[\sigma]$ is regularly varying with index $\alpha+\gamma-1$.
	\end{Enumeratealph}
	Then the following relations are true.
	\begin{Enumerate}
	\item
		Assume that \textup{(a)} holds.  Then $\alpha+\gamma\in[0,1]$.
		
		If $\alpha>0$, then \textup{(b)} holds and, as $x\to\infty$,
		\begin{equation}\label{L127}
			\ms S[\sigma](x) \sim
			\begin{cases}
				\displaystyle \frac{\pi\alpha}{\sin(\pi(\alpha+\gamma))}\cdot x^{\gamma-1}\nu([0,x)),
				& \alpha+\gamma\in(0,1),
				\\[3ex]
				\displaystyle \frac{\alpha}{x}\int_1^x t^{\gamma-1}\nu([0,t))\DD t,
				& \alpha+\gamma=0,
				\\[3ex]
				\displaystyle \alpha\int_x^\infty t^{\gamma-2}\nu([0,t))\DD t,
				& \alpha+\gamma=1.
			\end{cases}
		\end{equation}
		If $\alpha=0$ and $\gamma\in(0,1)$, then
		\begin{equation}\label{L120}
			\ms S[\sigma](x) \ll x^{\gamma-1}\nu([0,x)), \qquad x\to\infty.
		\end{equation}
	\item
		If \textup{(b)} holds and $\alpha+\gamma>0$, then also \textup{(a)} holds.
	\end{Enumerate}
\end{proposition}

\begin{proof}
	The relations in \eqref{L118} imply that the Stieltjes transform $\ms S[\sigma]$
	is well defined and that the measure $\sigma$ is infinite.
	Hence \eqref{L189} follows from \cref{L26}.

	(i)
	Assume that (a) holds.
	It follows from \eqref{L118} that $\int_{[1,\infty)}t^\gamma\DD\nu(t)=\infty$
	and $\int_{[1,\infty)}t^{\gamma-1}\DD\nu(t)<\infty$, which,
	by \cref{L33}\,(i), yields that $\alpha+\gamma\ge0$ and $\alpha+\gamma-1\le0$.

	Let us first consider the case when $\alpha>0$.
	The first relation in \eqref{L118} and \cref{L33}\,(i) imply that, as $x\to\infty$,
	\begin{equation}\label{L125}
		\sigma([0,x)) \sim \int_{[1,x)}t^\gamma\DD\nu(t)
		\sim \alpha\int_1^x t^{\gamma-1}\nu([0,t))\DD t
	\end{equation}
	and that $x\mapsto\sigma([0,x))$ is regularly varying with index $\alpha+\gamma$.
	Hence we can apply \cref{L1} to deduce that (b) holds and that
	\begin{equation}\label{L119}
		\ms S[\sigma](x) \sim C_{\alpha+\gamma}
		\int_x^\infty \frac{\sigma([0,t))}{t^2}\DD t,
	\end{equation}
	where $C_{\alpha+\gamma}$ is defined in \eqref{L129}.
	If $\alpha+\gamma<1$, then \cref{L36}\,(ii) and \eqref{L125} yield
	\[
		\ms S[\sigma](x) \sim \frac{C_{\alpha+\gamma}}{1-\alpha-\gamma}
		\cdot\frac{\sigma([0,x))}{x}
		\sim \alpha\frac{C_{\alpha+\gamma}}{1-\alpha-\gamma}
		\cdot\frac{1}{x}\int_1^x t^{\gamma-1}\nu([0,t))\DD t.
	\]
	In the case when $\alpha+\gamma=0$, the first fraction in front of the integral
	is equal to $1$.
	If $\alpha+\gamma\in(0,1)$, we can apply \cref{L36}\,(i) to arrive at
	\[
		\ms S[\sigma](x) \sim \frac{\alpha C_{\alpha+\gamma}}{1-\alpha-\gamma}
		\cdot\frac{1}{x}\cdot\frac{1}{\alpha+\gamma}\cdot x^\gamma\nu([0,x))
		= \frac{\pi\alpha}{\sin(\pi(\alpha+\gamma))}\cdot x^{\gamma-1}\nu([0,x)),
	\]
	which proves \eqref{L127} in the first two cases.
	
	Now assume that $\alpha+\gamma=1$.  We obtain from \eqref{L125} 
	and \cref{L36}\,(i) that
	\[
		\sigma([0,x)) \sim \alpha x^\gamma\nu([0,x)),
	\]
	which, together with \eqref{L119}, implies that
	\[
		\ms S[\sigma](x) \sim \alpha C_1\int_x^\infty t^{\gamma-2}\nu([0,t))\DD t,
	\]
	which proves \eqref{L127} in the third case.

	Let us now assume that $\alpha=0$ and $\gamma\in(0,1)$.
	It follows again from \cref{L33}\,(i) that
	\begin{equation}\label{L123}
		\sigma([0,x)) \sim \int_{[1,x)}t^\gamma\DD\nu(t)
		\ll \int_1^x t^{\gamma-1}\nu([0,t))\DD t.
	\end{equation}
	Define the measure $\mu$ on $[0,\infty)$
	by $\D\mu(t)=t^{\gamma-1}\nu([0,t))\DD t$.  Then
	\[
		\int_{[1,\infty)}\frac{1}{t}\DD\mu(t)
		= \int_{[1,\infty)}t^{\gamma-2}\nu([0,t))\DD t < \infty
	\]
	since the integrand in the second integral is regularly varying
	with index $\gamma-2<-1$.
	From \cref{L36}\,(i) we obtain that $x\mapsto\mu([0,x))$
	is regularly varying with index $\gamma$, and \eqref{L123} implies that
	\[
		\mu([0,x)) = \int_0^x t^{\gamma-1}\nu([0,t))\DD t
		\gg \sigma([0,x)).
	\]
	Hence we can apply \cref{L54} and \cref{L3} to deduce that
	\begin{align*}
		\ms S[\sigma](x) &\ll \ms S[\mu](x)
		\sim \frac{\pi\gamma}{\sin(\pi\gamma)}\cdot\frac{\mu([0,x))}{x}
		\\[1ex]
		&= \frac{\pi\gamma}{\sin(\pi\gamma)}\cdot\frac{1}{x}\int_0^x t^{\gamma-1}\nu([0,t))\DD t
		\sim \frac{\pi}{\sin(\pi\gamma)}\cdot\frac{1}{x}x^\gamma\nu([0,x)),
	\end{align*}
	which proves \eqref{L120}.

	(ii)
	Assume that (b) holds and that $\alpha+\gamma>0$.  It follows from \cref{L1}
	that $x\mapsto\sigma([0,x))$ is regularly varying with index $\alpha+\gamma$.
	Since $\D\nu(t)=\frac{1}{h(t)}\DD\sigma(t)$, $t\in(0,\infty)$, we obtain 
	from \cref{L33}\,(i) that (a) holds.
\end{proof}

\medskip\noindent
The following example shows that the implication (b)\,$\Rightarrow$\,(a) is not
valid in general if $\alpha+\gamma=0$.
The example is also used in \cref{L208,L212}.

\begin{example}\label{L202}
	Let $h:[0,\infty)\to[0,\infty)$ be a continuous function such that
	$h(t)>0$ for $t>0$ and $h(t)\sim t^\gamma$ with $\gamma<0$.
	Further, let $\sigma$ be the discrete measure on $[0,\infty)$ as in \cref{L69},
	i.e.\ with point masses $\sigma(\{e^k\})=1$, $k\in\bb N$.
	Let $\nu$ be the measure on $[0,\infty)$ such that $\D\sigma(t)=h(t)\DD\nu(t)$,
	$t\in(0,\infty)$ and $\nu(\{0\})=0$; it is also discrete and has point 
	masses 
	\[
		\nu(\{e^k\})=\frac{1}{h(e^k)}, \qquad k\in\bb N.
	\]
	Since
	\begin{align*}
		\int_{[0,\infty)}h(t)\DD\nu(t) &= \int_{[0,\infty)}\hspace*{1ex}\D\sigma(t) = \infty,
		\\[1ex]
		\int_{[0,\infty)}\frac{h(t)}{1+t}\DD\nu(t) 
		&= \int_{[0,\infty)}\frac{1}{1+t}\DD\sigma(t)
		= \sum_{k=1}^\infty\frac{1}{1+e^k} < \infty,
	\end{align*}
	condition \eqref{L118} is satisfied.
	\Cref{L69} shows that the distribution function $t\mapsto\sigma([0,t))$ is 
	slowly varying.  Hence, by \cref{L1}, $\ms S[\sigma]$ is regularly varying
	with index $-1$, i.e.\ \textup{(b)} in \cref{L122} is satisfied with $\alpha+\gamma-1=-1$.
	Note that $\sigma([0,t))\sim\log t$ and therefore 
	\begin{equation}\label{L209}
		\ms S[\sigma](x)\sim\frac{\log x}{x},
		\qquad x\to\infty,
	\end{equation}
	by \cref{L3}.
	On the other hand, \cref{L69} with $g$ such that $g(t)=\frac{1}{h(t)}$, $t\in[1,\infty)$,
	also implies that 
	\[
		x\mapsto\nu([0,x)) = \int_{[0,x)}g(t)\DD\sigma(t) = f(x),
	\]
	with $f$ from \cref{L69},
	is not regularly varying, i.e.\ \textup{(a)} in \cref{L122} is not satisfied.
	Hence the implication \textup{(b)}\,$\Rightarrow$\,\textup{(a)} does not hold in 
	general when $\alpha+\gamma=0$.
\end{example}

\medskip\noindent
Finally, we need a comparison result for Stieltjes transforms of measures
as in \cref{L122}.

\begin{lemma}\label{L121}
	Let $\nu_1,\nu_2$ be measures on $[0,\infty)$, let $h:[0,\infty)\to[0,\infty)$
	be a continuous function such that $h(t)>0$ for $t>0$ and $h(t)\sim t^\gamma$, $t\to\infty$,
	with some $\gamma\in\bb R\setminus\{0\}$.
	Assume that
	\begin{equation}\label{L192}
		\int_{[0,\infty)}h(t)\DD\nu_2(t) = \infty
		\qquad\text{and}\qquad
		\int_{[0,\infty)}\frac{h(t)}{1+t}\DD\nu_2(t) < \infty
	\end{equation}
	and define the measures $\sigma_i$ on $[0,\infty)$ by
	$\D\sigma_i(t)=h(t)\DD\nu_i(t)$, $t\in[0,\infty)$, $i\in\{1,2\}$.
	Further, assume that $t\mapsto\nu_2([0,t))$ is regularly varying
	with strictly positive index and that the limit
	\begin{equation}\label{L193}
		\lim_{t\to\infty}\frac{\nu_1([0,t))}{\nu_2([0,t))}
	\end{equation}
	exists in $[0,\infty)$.  Then $\ms S[\sigma_1]$ is well defined and
	\begin{equation}\label{L194}
		\lim_{x\to\infty}\frac{\ms S[\sigma_1](x)}{\ms S[\sigma_2](x)}
		= \lim_{t\to\infty}\frac{\nu_1([0,t))}{\nu_2([0,t))}.
	\end{equation}
\end{lemma}

\begin{proof}
	Since $\int_{[1,\infty)}t^{\gamma-1}\DD\nu_2(t)<\infty$ by \eqref{L192},
	we can use the existence of the limit in \eqref{L193} and 
	\cref{L59}\,(i) (when $\gamma<1$) or \cref{L59}\,(ii) (when $\gamma>1$) 
	to deduce that $\int_{[1,\infty)}t^{\gamma-1}\DD\nu_1(t)<\infty$
	(for $\gamma=1$ this follows directly).
	Hence $\ms S[\sigma_1]$ is well defined.
	
	Denote the limit in \eqref{L193} by $c$.
	Let us first consider the case when $c>0$.  Then $t\mapsto\nu_1([0,t))$
	is regularly varying with the same index as $t\mapsto\nu_2([0,t))$,
	and $\int_{[1,\infty)}t^\gamma\DD\nu_1(t)=\infty$
	(again by \cref{L59}).
	We now obtain from \cref{L59}\,(i) (when $\gamma<0$) or \cref{L59}\,(iii) (when $\gamma>0$) 
	that
	\[
		\lim_{x\to\infty}\frac{\sigma_1([0,x))}{\sigma_2([0,x))}
		= \lim_{x\to\infty}\frac{\int_{[1,x)}t^\gamma\DD\nu_1(t)}{\int_{[1,x)}t^\gamma\DD\nu_2(t)}
		= c.
	\]
	By \cref{L54}, this implies that \eqref{L194} holds.
	
	Now assume that $c=0$.  From \cref{L59}\,(i) and (ii) we can deduce that
	\[
		\sigma_1([0,x)) \asymp \int_{[1,x)}t^\gamma\DD\nu_1(t)
		\ll \int_{[1,x)}t^\gamma\DD\nu_2(t)
		\sim \sigma_2([0,x)),
	\]
	which, together with \cref{L54}, yields \eqref{L194}.
\end{proof}

%---------
%   FINISH
%---------

{\footnotesize
\begin{flushleft}
	M.~Langer \\
	Department of Mathematics and Statistics \\
	University of Strathclyde \\
	26 Richmond Street \\
	Glasgow G1 1XH \\
	UNITED KINGDOM \\
	email: \texttt{m.langer@strath.ac.uk} \\[5mm]
\end{flushleft}
\begin{flushleft}
	H.\,Woracek\\
	Institute for Analysis and Scientific Computing\\
	Vienna University of Technology\\
	Wiedner Hauptstra{\ss}e\ 8--10/101\\
	1040 Wien\\
	AUSTRIA\\
	email: \texttt{harald.woracek@tuwien.ac.at}\\[5mm]
\end{flushleft}
}

\end{document}